\RequirePackage{fix-cm}
\pdfoutput=1
\documentclass[leqno, fleqn, centertags, 12pt]{article}

\usepackage{latexsym}
\usepackage{amsmath}

\usepackage[T1]{fontenc}
\usepackage{fourier}
\usepackage[mathbf]{euler}
\newcommand{\mathbold}[1]{\mathbf{#1}}

\usepackage{comment}
\usepackage{fullpage}
\usepackage{tikz-cd}

\newcommand{\dis}{\displaystyle}
\newcommand{\old}[1]{\oldstylenums{#1}}
\newcommand{\oold}{\oldstylenums{0}}

\makeatletter
\renewcommand*\@arabic[1]{\oldstylenums{\number#1}}
\renewcommand\makeatother

\newcommand{\id}{\mathop{\mbox{id}\fff}\nolimits}

\newcommand{\mmod}{\mathop{\mbox{Mod}}\hspace{0.075em}}
\newcommand{\mm}{\mathop{\mbox{Mod}}\hspace{0.075em}}
\newcommand{\ttt}{\mathop{\mathcal I}\hspace*{0.05em}}
\newcommand{\ppp}{\mathop{\mathcal P}\hspace*{0.05em}}

\newcommand{\zero}{_{\dff \oldstylenums{0}\fff}}
\newcommand{\one}{_{\dff \oldstylenums{1}\fff}}
\newcommand{\two}{_{\dff \oldstylenums{2}\fff}}
\newcommand{\three}{_{\dff \oldstylenums{3}\fff}}
\newcommand{\four}{_{\dff \oldstylenums{4}\fff}}

\newcommand{\zzz}{\mathbb{Z}}

\newcommand{\toto}{\longrightarrow}
\newcommand{\ttoo}{\hspace*{0.2em}\longrightarrow\hspace*{0.2em}}

\newcommand{\tors}{\ttt(S)}
\newcommand{\mms}{\mmod(S)}

\def\sss{\hspace{0.05em}\ }
\def\dss{\hspace{0.1em}\ }
\def\trs{\hspace{0.15em}\ }
\def\qss{\hspace{0.2em}\ }
\def\pss{\hspace{0.3em}\ }
\def\oss{\hspace{0.4em}\ }

\def\halfff{\hspace*{0.025em}}
\def\fff{\hspace*{0.05em}}
\def\dff{\hspace*{0.1em}}
\def\trf{\hspace*{0.15em}}
\def\qff{\hspace*{0.2em}}
\def\pff{\hspace*{0.3em}}
\def\off{\hspace*{0.4em}}

\def\ttff{{\hspace*{-0.05em}--\hspace*{0.15em}}}

\def\ffdot{\hspace*{-0.1em}.\hspace*{0.2em}\ }
\def\dfdot{\hspace*{-0.2em}.\hspace*{0.4em}\ }

\def\dfcom{\hspace*{-0.2em},\hspace*{0.4em}\ }

\newcommand{\hnsp}{\hspace*{-0.05em}}
\newcommand{\nsp}{\hspace*{-0.1em}}
\newcommand{\nnsp}{\hspace*{-0.15em}}
\newcommand{\snsp}{\hspace*{-0.175em}}
\newcommand{\dnsp}{\hspace*{-0.2em}}
\newcommand{\tnsp}{\hspace*{-0.3em}}

\newcommand{\ccomp}{\mathsf{c}}

\newcommand{\cutc}{\nsp/\tnsp/\hspace{-0.08em}c}

\newcommand{\hclass}[1]{[\dff #1\dff]}

\makeatletter
\renewcommand{\@makefntext}[1]{\vspace*{0.5ex}
\parindent=0em
\hspace*{-0.4em}
\hbox to 0.4em{\hss\@makefnmark}\hspace*{0.4em}{#1}
}
\makeatother

\newcommand{\proof}{\vspace{0.5\bigskipamount}{\textbf{{\emph{Proof}.}}\hspace*{0.7em}}}
\newcommand{\prooftitle}[1]{\vspace*{0.5\bigskipamount}{\textbf{{\emph{#1}.}}\hspace{0.7em}}}
\newcommand{\eproof}{ $\blacksquare$}

\newcounter{mysectionnumber}
\newcommand{\mysection}[2]{
\setcounter{equation}{0}\refstepcounter{mysectionnumber}
\section*{ \textnormal{{\old{\themysectionnumber}.} {#1}}}\label{#2}}

\newcounter{myparnum}[mysectionnumber]
\newcommand{\mypar}[2]{\refstepcounter{myparnum}{\vspace{\medskipamount}\textbf{{\themyparnum. #1\label{#2}}}\hspace{0.5em}}}
\renewcommand{\themyparnum}{\old{\themysectionnumber}.\old{\arabic{myparnum}}}

\newcounter{mylemmanum}[myparnum]
\numberwithin{equation}{section}

\newcommand{\myitpar}[1]{\vspace{\medskipamount}\textbf{\textit{#1}}\hspace*{0.5em}}
\newcommand{\myit}[1]{\textbf{\textit{#1}}\hspace{0.0em}}

\newcounter{myapparnum}[mysectionnumber]
\newcommand{\myappar}[2]{
\refstepcounter{myapparnum}
\vspace{\bigskipamount}
\textbf{{\themyapparnum. #1}\label{#2}}
\hspace{0.5em}}
\renewcommand{\themyapparnum}{A\fff. 
\arabic{myapparnum}}

\newcounter{myappendnumber}
\newcommand{\myappend}[2]{\setcounter{footnote}{0}
\setcounter{myapparnum}{0}
\refstepcounter{myappendnumber}
\section*{ \textnormal{Appendix\fff.\oss 
{#1}}}\label{#2}}

\title{Elements\qss of\qss Torelli\qss topology\halfff:\\
II.\qss The\qss extension\qss problem}
\author{Nikolai\dss V.\dss Ivanov}
\date{}

\begin{document}

\setlength{\baselineskip}{12pt plus 0pt minus 0pt}
\setlength{\parskip}{12pt plus 0pt minus 0pt}
\setlength{\abovedisplayskip}{12pt plus 0pt minus 0pt}
\setlength{\belowdisplayskip}{12pt plus 0pt minus 0pt}

\newskip\smallskipamount \smallskipamount=3pt plus 0pt minus 0pt
\newskip\medskipamount   \medskipamount  =6pt plus 0pt minus 0pt
\newskip\bigskipamount   \bigskipamount =12pt plus 0pt minus 0pt

\maketitle

\vspace*{12ex}

\myit{\hspace*{0em}\large Contents}\vspace*{1ex} \vspace*{\bigskipamount}\\ 
\myit{\phantom{1}\old{1}.}\hspace*{0.5em} Introduction \hspace*{0.5em} \vspace*{0.25ex}\\
\myit{\phantom{1}\old{2}.}\hspace*{0.5em} Surfaces, circles, and twists  \hspace*{0.5em} \vspace*{0.25ex}\\
\myit{\phantom{1}\old{3}.}\hspace*{0.5em} Pure diffeomorphisms and reduction systems  \hspace*{0.5em} \vspace*{0.25ex}\\
\myit{\phantom{1}\old{4}.}\hspace*{0.5em} The sirens of Torelli topology  \hspace*{0.5em} \vspace*{0.25ex}\\
\myit{\phantom{1}\old{5}.}\hspace*{0.5em} Weakly Torelli diffeomorphisms  \hspace*{0.5em} \vspace*{0.25ex}\\
\myit{\phantom{1}\old{6}.}\hspace*{0.5em} Extensions by the identity and the symmetry property  \hspace*{0.5em} \vspace*{0.25ex}\\
\myit{\phantom{1}\old{7}.}\hspace*{0.5em} General extensions  \hspace*{0.5em} \vspace*{0.25ex}\\
\myit{\phantom{1}\old{8}.}\hspace*{0.5em} Multi-twists about the boundary  \hspace*{0.5em} \vspace*{0.25ex}\\
\myit{\phantom{1}\old{9}.}\hspace*{0.5em} Torelli groups of surfaces with boundary\fff? \hspace*{0.5em} 
\vspace*{1.5ex}\\
\myit{Appendix.}\hspace*{0.5em} 
A converse to Theorem\qss \ref{weakly-torelli} \hspace*{0.5em} \vspace*{1.5ex}\\
\myit{References}\hspace*{0.5em}  \hspace*{0.5em}  \vspace*{0.25ex}

\footnotetext{\hspace*{-0.65em}\copyright\ 
Nikolai\dss V.\qss Ivanov,\oss \old{2017}.\oss 
Neither the work reported in this paper,\qss 
nor its preparation were 
supported by any governmental 
or non-governmental agency,\qss 
foundation,\qss 
or institution.}

\vspace{12ex}

\renewcommand{\baselinestretch}{1}
\selectfont

\newpage
\mysection{Introduction}{introduction}

\vspace*{\medskipamount}
\myitpar{Torelli groups and diffeomorphisms.} The\qss 
Teichm\"{u}ller modular group\qss 
$\mms$\qss of an orientable surface\dss $S$\dss is defined as the group of isotopy classes of
orientation-preserving diffeomorphisms of\dss $S$\nnsp.\oss
Both diffeomorphisms and isotopies are required to preserve the boundary\qss $\partial S$\qss
only set-wise.\oss
If\dss $S$\dss is closed,\oss the subgroup of elements\qss $\mmod(S)$\qss
acting trivially on\qss $H\one(S\fff,\pff \zzz)$\qss
is called the\qss \emph{Torelli\dss group}\qss of\dss $S$\dss
and is denoted by\qss $\tors$\dnsp.\oss
We will denote the homology group\qss $H\one(S\fff,\pff \zzz)$\qss
simply by\qss $H\one(S)$\dnsp.\oss
A diffeomorphism of a closed surface\dss $S$\dss is called a\qss \emph{Torelli diffeomorphism}\qss
if it acts trivially on\qss $H\one(S)$\dnsp.

There are good reasons to hesitate with defining Torelli groups 
and Torelli diffeomorphisms for surfaces with boundary.\oss
Requiring diffeomorphisms to act trivially
on\dss $H\one(S)$\dss in the case when\qss $\partial S\qff \neq\qff \emptyset$\trs 
seems to be fairly naive.\oss

\myitpar{The extension problem.} This paper is devoted to the\qss 
\emph{extension problem}\qss in Torelli to\-po\-lo\-gy.\oss 
By this we understand the following circle of questions.
Let\dss $Q$\dss be a subsurface of a connected
closed orientable surface\dss $S$\dnsp.\oss
We are interested in extensions of diffeomorphisms of\dss $Q$\dss 
to diffeomorphisms of\dss $S$\dss acting trivially on\dss $H\one(S)$\dnsp,\oss
i.e.\qss in extensions which are Torelli diffeomorphisms.\oss
As usual,\oss the case of connected subsurfaces\dss $Q$\dss is the most important one.

The theory of Teichm\"{u}ller modular groups suggests to consider
only diffeomorphisms of\dss $Q$\dss equal to the identity on\qss $\partial Q$\dnsp,\oss
and that the appropriate isotopies of such diffeomorphisms are the isotopies fixed on\qss $\partial Q$\dnsp.\oss
At the same times this allows to simultaneously address the question of when
the extension\dss $\varphi\backslash S$\dss of a diffeomorphism\dss 
$\varphi$\dss of\dss $Q$\dss by the identity to 
a diffeomorphism of\dss $S$\dss is a Torelli diffeomorphism.\oss
The latter question was already addressed by A. Putman,\oss
whose paper\qss \cite{p}\qss was one of the main sources 
of inspiration for the present paper\halfff.\oss 

Closely related is the question of when there exists a multi-twist diffeomorphism\dss $\tau$\dss
about\qss $\partial Q$\qss 
such that\qss $\tau\qff\circ\qff (\varphi\backslash S\halfff)$\qss
is a Torelli diffeomorphism.\oss
For some applications this is more important
than\dss $\varphi\backslash S$\dss itself\dss being\sss a Torelli diffeomorphism.\oss

\myitpar{The ambient surface\sss $S$\sss as a structure on\dss $Q$\dnsp.}\qss 
It is hardly surprising that the answers to such questions 
depend on the embedding of\dss $Q$\dss
into the ambient closed surface\dss $S$\dnsp.\oss
But they depend this embedding
only in a mild and controlled way discovered by A. Putman\qss \cite{p}.\oss
In contrast with\qss \cite{p},\oss we prefer not to spell out in advance
what structure on the surface\dss $Q$\dss is induced by its the embedding in\dss $S$\nnsp.\oss 
We treat the embedding of\dss $Q$\dss in\dss $S$\dss 
as a part of the structure of\dss $Q$\snsp,\oss
keeping the mind open to using whatever part of this structure may be needed.\oss

Let\qss $c\qff =\qff \partial Q$\qss 
and let\qss $\ccomp Q$\qss be the closure of the complement\qss $S\smallsetminus Q$\nnsp.\oss
Clearly,\pss $\ccomp Q$\dss
is a subsurface of\dss $S$\dss with the boundary\dss $c$\nnsp.\oss
If\dss $Q$\dss is connected,\oss
then the answers to these questions depend only on the partition of\qss 
$c$\dss
into the boundaries\dss $\partial P$\dss 
of components\dss $P$\dss of\qss $\ccomp Q$\dnsp.

\myitpar{Weakly Torelli diffeomorphisms.}
Our answers are stated in terms of homology groups related to\dss $Q$\qss
and the ambient surface\dss $S$\dnsp.\oss
We will use the singular homology theory and work with singular chains when
this is either necessary or more natural than working with homology classes.\oss

An obvious necessary condition for the existence of a Torelli diffeomorphism of\dss
$S$ extending a diffeomorphism\dss $\varphi$\dss of\dss $Q$\dss
is the following\fff:\oss for every cycle\dss $\sigma$\dss in\dss $Q$\dss
the cycle\dss $\varphi_{\fff *\fff}(\sigma)$\dss 
should be homologous to\dss $\sigma$\dss in\dss $S$\nnsp.\oss
We call diffeomorphisms\dss $\varphi$\dss of\dss $Q$\dss
fixed on\dss $c$\dss and having
this property\qss \emph{weakly Torelli diffeomorphisms}\pss of\dss $Q$\nnsp.\oss

\myitpar{Homology theory.}
We use the standard notations\qss $C_{\fff *}(\fff\bullet\fff)$\qss and\dss $\partial$\dss 
for the groups of singular chains and boundary maps.\oss
Given a topological space\dss $X$\dss and a subspace\dss $A$\dss of\dss $X$\nnsp,\oss
we denote by\dss $Z\one(X\fff,\pff A)$\dss the group of relative cycles of\dss
$(X\fff,\pff A)$\dnsp,\oss i.e.\qss the group of chains\qss
$\alpha\qff \in\qff C_{\fff *\fff}(X)$\qss
such that\qss
$\partial \alpha\qff \in\qff C_{\fff *\fff}(A)$\dnsp.\oss
By\dss $\hclass{\sigma}$\sss we denote 
the homology class of a cycle\qss $\sigma$\nnsp.\oss

\myitpar{The action on relative cycles.}
Let\qss $\mathbold{Z}\one(Q\fff,\pff c)$\qss be the group of 
relative cycles\qss $\alpha\qff \in\qff Z\one(Q\fff,\pff c)$\qss
such that the boundary\dss $\partial \alpha$\dss 
is a boundary in\dss $\ccomp Q$\dss also.\oss

Let\dss $\varphi$\dss
be a weakly Torelli diffeomorphism of\dss $Q$\nnsp.\oss
It turns out that understanding the action of\dss $\varphi$\dss
on relative cycles\qss 
$\alpha\qff \in\qff \mathbold{Z}\one(Q\fff,\pff c)$\qss 
is the key to the extension problem for\dss $\varphi$\nnsp.\oss
If\dss $\alpha$\dss is an arbitrary relative cycle of\qss $(Q\fff,\pff c)$\nnsp,\oss
then\qss $\varphi_{\fff *\fff}(\alpha)\qff -\qff \alpha$\qss is a cycle in\dss $Q$\dss
because\dss $\varphi$\dss is fixed on\dss $c$\nnsp.\oss
The following theorem is the first main result of this paper\halfff.\oss

\myitpar{Theorem.}
\emph{If\qss $\varphi$\qss is a weakly Torelli diffeomorphism,\oss
then for every\qss $\alpha\qff \in\qff \mathbold{Z}\one(Q\fff,\pff c)$\qss
the cycle\qss
$\varphi_{\dff *}(\alpha)\qff -\qff \alpha$\qss 
is homologous in\dss $Q$\dss to a cycle in\dss $c$\nnsp.\oss}

\vspace*{6pt}
See Theorem\qss \ref{weakly-torelli}.\oss
It turns out that the homology classes of cycles in\dss $c$\dss
homologous to the cycle\qss
$\varphi_{\dff *}(\alpha)\qff -\qff \alpha$\qss 
in\dss $Q$\dss to a big extent depend only on\dss $\partial \alpha$\dss
and even only on the homology class\qss $\hclass{\partial \alpha}\qff \in\qff H\zero(c)$\dnsp.\oss
In order to state this in a precise and efficient form,\oss
one need to introduce two homology groups reflecting the way
the manifold\dss $c$\dss is situated in\dss $S$\nnsp.\oss

Let\qss $\mathbold{K}\zero(c)$\qss be the intersection of the kernels 
of the inclusion homomorphisms\qss\vspace*{3pt}
\[
\quad
H\zero(c)\ttoo H\zero(Q)
\hspace*{1.5em}\mbox{ and }\hspace*{1.5em}
H\zero(c)\ttoo H\zero(\ccomp Q)\dff.
\]

\vspace*{-9pt}
A relative cycle\qss $\alpha\qff \in\qff Z\one(Q\fff,\pff c)$\qss
belongs to\qss $\mathbold{Z}\one(Q\fff,\pff c)$\qss
if and only if\qss
$\hclass{\partial \alpha}\qff \in\qff \mathbold{K}\zero(c)$\dnsp.\qff\oss
Let\qss
$\mathbold{H}\one(c)$\qss be the quotient of\dss $H\one(c)$\dss
by the kernel of the inclusion homomorphism\qss\vspace*{3pt}
\[
\quad
\iota
\qff \colon\qff
H\one(c)\ttoo H\one(S)\dff.
\]

\vspace*{-9pt}
The group\qss $\mathbold{H}\one(c)$\qss is canonically isomorphic to the image of\dss
$\iota$\nnsp.\oss 
The advantage of the group\qss $\mathbold{H}\one(c)$\qss
is that it lives on\dss $c$\dss
instead of\dss $S$\nnsp.\qff\oss
If\dss $Q$\dss is connected,\oss
then the groups\qss $\mathbold{K}\zero(c)$\qss
and\qss $\mathbold{H}\one(c)$\qss
depend only on the partition of\qss 
$c$\dss
into the boundaries\dss $\partial P$\dss 
of components\dss $P$\dss of\qss $\ccomp Q$\dnsp.

\myitpar{Observation.}
\emph{Suppose that\qss $\varphi$\qss is a weakly Torelli diffeomorphism
and\qss $\alpha\qff \in\qff \mathbold{Z}\one(Q\fff,\pff c)$\dnsp.\oss
Let\qss $\gamma\qff \in\qff C\one(c)$\qss be a cycle homologous in\dss $Q$\dss to 
the cycle\qss
$\varphi_{\dff *}(\alpha)\qff -\qff \alpha$\nnsp.\oss
Then the image of the homology class\qss
$\hclass{\gamma}\qff \in\qff H\one(c)$\qss
in the group\qss $\mathbold{H}\one(c)$\qss
depends only on\qss
$\hclass{\partial \alpha}\qff \in\qff H\zero(c)$\dnsp.\oss}

\vspace*{6pt}
See Section\qss \ref{weakly-torelli-section}.\oss
The image
of\dss $\hclass{\gamma}$\dss in the group\dss $\mathbold{H}\one(c)$\dss 
is called the\qss \emph{difference class}\qss of\dss $\alpha$\dss
and is denoted by\qss $\Delta_{\dff \varphi}\dff(\alpha)$\dnsp.\oss
The difference classes\qss 
$\Delta_{\dff \varphi}\dff(\alpha)$\qss lead to a homomorphism\vspace*{3pt}
\[
\quad
\delta_{\fff \varphi}
\qff \colon\qff
\mathbold{K}\zero(c)\ttoo \mathbold{H}\one(c)
\]

\vspace*{-9pt}
called the\dss \emph{$\delta$\dnsp-difference map}\pss of\dss $\varphi$\nnsp.\oss
The\qss $\delta$\dnsp-difference map of\dss $\varphi$\dss 
turns out to be the obstruction for the extension\qss
$\varphi\backslash S$\qss of\dss $\varphi$\dss by the identity
to be a Torelli diffeomorphism of\dss $S$\nnsp.\oss

\myitpar{Theorem.}\oss \emph{The extension\qss $\varphi\nsp\backslash S$\qss
of\dss $\varphi$\sss to a diffeomorphism of\dss $S$\dss by the identity        
is a Torelli diffeomorphism if and only if\qss $\varphi$\dss 
is a weakly Torelli diffeomorphism
and\qss $\delta_{\fff \varphi}\off =\off \oold$\nnsp.\oss}

\vspace*{\medskipamount}
See Theorem\qss \ref{extension-by-identity}.\oss
This theorem provides a characterization of diffeomorphisms\dss $\varphi$\dss
such that\qss $\varphi\nsp\backslash S$\qss is a Torelli diffeomorphism 
in terms of the action of\dss $\varphi$\dss on\dss $Q$\dss
and some minor and independent of\dss $\varphi$\dss 
information about the embedding of\dss $Q$\dss into\dss $S$\dnsp.\oss
For connected subsurfaces\dss $Q$\dss it is closely related to a
theorem of A.\qss Putman.\oss
See A.\qss Putman\oss \cite{p},\oss Theorem\qss 3.3.\oss

\myitpar{Symmetric homomorphisms.}
The intersection pairing\oss
$\displaystyle
H\zero(c)\dff \times H\one(c) \ttoo \zzz$\oss
leads to a ca\-non\-i\-cal pairing\oss
$\displaystyle
\mathbold{K}\zero(c)\dff \times \mathbold{H}\one(c) \ttoo \zzz$\nsp,\oss
denoted by\qss
$\displaystyle
\langle\qff \bullet\fff,\pff \bullet \qff\rangle_{\dff c}$\dnsp.\oss
A map\oss
$\displaystyle
\delta
\dff \colon\dff
\mathbold{K}\zero(c)\ttoo \mathbold{H}\one(c)$\oss
is called\qss \emph{symmetric}\pss if\oss\vspace*{3pt}
\[
\quad
\langle\qff a\halfff\fff,\off \delta(b) \qff\rangle_{\dff c}
\off =\off
\langle\qff b\halfff,\off \delta(a) \qff\rangle_{\dff c}
\]

\vspace*{-9pt}
for all\qss $a\fff,\pff b\qff\in\qff  \mathbold{K}\zero(c)$\dnsp.\oss
It turns out that the\dss $\delta$\dnsp-difference map\dss $\delta_{\fff \varphi}$\dss 
is symmetric for every weakly Torelli
diffeomorphism\dss $\varphi$\nnsp.\oss
See Theorem\qss \ref{symmetry-theorem}.

\myitpar{Completely reducible homomorphisms.}
Suppose now that the subsurface\dss $Q$\dss is connected.\oss
Then the decomposition of\qss $c\qff =\qff \partial Q$\dss
into the disjoint union of the boundaries\dss $\partial P$\dss 
of the components\dss $P$\dss of\qss $\ccomp Q$\qss leads to
direct sum decompositions of the groups\qss
$\displaystyle
\mathbold{K}\zero(c)$\qss
and\qss
$\displaystyle
\mathbold{H}\one(c)$\dnsp.\oss
A homomorphism\oss
$\displaystyle
\mathbold{K}\zero(c)\ttoo \mathbold{H}\one(c)$\oss is said to be\qss
\emph{completely reducible}\pss if\dss it\dss respects these direct sum decompositions.\oss

\myitpar{Theorem.}\oss \emph{The diffeomorphism\dss $\varphi$\dss can be extended
to a Torelli diffeomorphism of\dss $S$\dss if and only if\dss $\varphi$\dss is peripheral
and the\dss $\delta$\dnsp-difference map\dss $\delta_{\fff \varphi}$\dss is completely reducible.\oss}

\vspace{\medskipamount}
See Theorem\qss \ref{char-ext}.\oss 
This theorem is complemented by the following
theorem\qss (which is,\oss in fact\halfff,\oss used in its proof\fff)\qss
about the realization of completely reducible symmetric homomorphisms\oss
$\displaystyle
\mathbold{K}\zero(c)\ttoo \mathbold{H}\one(c)$\oss 
as\dss $\delta$\dnsp-difference maps.\oss

\myitpar{Theorem.}\oss \emph{If\oss 
$\displaystyle
\delta
\dff \colon\dff 
\mathbold{K}\zero(c)\ttoo \mathbold{H}\one(c)$\oss
is a completely reducible and symmetric,\oss
then\qss $\delta\qff =\qff \delta_{\fff \varphi}$ for some diffeomorphism\dss $\varphi$\dss
of\qss $Q$\qss which is equal to the identity on\qss $\partial Q$\qss
and can be extended to a Torelli diffeomorphism of\dss $S$\nnsp.\oss}

\vspace*{\medskipamount}
See Theorem\qss \ref{realization-of-delta}.\oss
The above characterizations of diffeomorphisms\dss $\varphi$\dss such that\qss
$\varphi\nsp\backslash S$\qss is a Torelli diffeomorphisms and such that\dss $\varphi$\dss
can be extended to a Torelli diffeomorphism of\dss $S$\dss imply,\oss
together with the last theorem,\oss
that neither of these properties can be characterized entirely in terms
of action of\dss $\varphi$\dss on\dss $Q$\dnsp.\oss
These characterizations do depend on the embedding of\dss $Q$\dss in\dss $S$\dnsp,\oss
although in a minor and controlled way.\oss
For diffeomorphisms\dss $\varphi$\dss such that\qss 
$\varphi\nsp\backslash S$\qss is a Torelli diffeomorphism
this was observed by A.\qss Putman\oss \cite{p}.

\myitpar{Diagonal maps.} A homomorphism\oss
$\displaystyle
H\zero(c)\ttoo \mathbold{H}\one(c)$\oss
is called\qss \emph{diagonal}\qss if
for every component\dss $C$\dss of\dss $c$\dss
it maps the homology class of any $\oold$\dnsp-simplex in\qss $C$\qss 
to the image in\dss $\mathbold{H}\one(c)$\dss of an integer multiple of the fundamental class\qss
$\hclass{C}$\nsp.\oss

\myitpar{Theorem.} \emph{There exists a multi-twist diffeomorphism\dss $\tau$\dss 
about\qss $\partial Q$\qss such that\qss $\tau\qff \circ\qff (\varphi\nsp\backslash S\halfff)$\qss
is a Torelli diffeomorphism if and only if\pss
$\varphi$\pss a weakly Torelli diffeomorphism and its difference map\qss
$\delta_{\fff \varphi}$\qss is equal to
the restriction to\qss $\mathbold{K}\zero(c)$\qss of a diagonal map\oss
$\displaystyle
H\zero(c)\ttoo \mathbold{H}\one(c)$\dnsp.\oss}

\vspace*{\medskipamount}
See Theorem\qss \ref{twist-correctable}.\oss
It turns out that if each component of\qss $\ccomp Q$\qss 
has\qss $\leqslant\qff 3$\qss boundary components,\oss
then every completely reducible symmetric map is equal
to the restriction of a diagonal map.\oss
See Theorem\qss \ref{3-diagonal}.\oss
This leads to the following theorem.

\myitpar{Theorem.} \emph{Suppose that every component of\qss $\ccomp Q$\qss 
has\oss $\leqslant\qff 3$\qss boundary components.\oss
If\qss $\varphi$\qss can be extended to a Torelli diffeomorphism of\dss $S$\nnsp,\oss 
then there exists a multi-twist diffeomorphism\dss $\tau$\dss
about\qss $\partial Q$\qss such that\qss $\tau\qff \circ\qff (\varphi\nsp\backslash S\halfff)$\qss
is a Torelli diffeomorphism.\oss}

\vspace*{\medskipamount}
See Theorem\qss \ref{3-components-identity}.\oss
As the examples in Section\qss \ref{extensions}\qss show,\oss
the assumption that every component of\qss $\ccomp Q$\qss 
has\qss $\leqslant\qff 3$\qss boundary components cannot be 
relaxed even to\qss $\leqslant\qff 4$\nnsp.\oss

\myitpar{Outline of the paper\halfff.} 
Section\dss \ref{prelim}\dss
is devoted mostly to fixing the meaning of such common terms as
a subsurface,\qss a circle,\qss a twist\halfff.\oss
Section\dss \ref{pure-reduction}\dss is devoted to a review of reduction
systems and of their basic properties.\oss
This material is used only in Section\qss \ref{extensions}\fff.\oss
Section\qss \ref{extensions}\qss is devoted
partly to the motivation behind this paper\halfff,\qss
partly to the dangers one may encounter while venturing in the Torelli sea
after gaining experience in the Teichm\"{u}ller modular groups waters.\oss
The next four sections are devoted to 
the detailed statements and proofs of the results outlined above.\oss
Section\dss \ref{torelli-groups}\dss is devoted to speculations about
how one should define Torelli groups for surfaces with non-empty boundary.\oss
While Section\dss \ref{torelli-groups} depends on the previous sections,\oss
it is more of a continuation of Section\dss \ref{extensions}\dss
than of Sections\dss \ref{weakly-torelli-section}\trf{\ttff}\ref{multi-twists-about-the-boundary}.\oss
The Appendix is devoted to a converse of the first theorem above.

\mysection{Surfaces,\qss circles,\qss and\qss twists}{prelim}

\vspace*{\medskipamount}
\myitpar{Surfaces.} By a\qss \emph{surface}\qss we understand a compact orientable 
$\oldstylenums{2}$\dnsp-manifold 
with a possibly empty boundary.\oss
We denote by\dss $S$\dss a fixed connected oriented surface.\oss
By a\qss \emph{subsurface}\qss of $S$ we understand a codimension\dss $\oold$\dss submanifold\qss
$Q$\qss of\qss $S$\qss such that each component of\qss $\partial Q$\qss is either equal to
a component of\qss $\partial S$\dfcom or disjoint from\qss $\partial S$\dfdot 
For a subsurface $Q$ of $S$ we will denote by $\ccomp Q$ 
the closure of its set-theoretic complement $S\smallsetminus\nsp Q$\dfdot
Clearly,\oss $\ccomp Q$\qss is also a subsurface of\qss $S$\halfff\dfdot

\myitpar{Circles.} A\dss \emph{circle}\qss on\dss $S$\dss is defined as
a submanifold of\dss $S$\dss diffeomorphic to 
the standard circle\dss $S^{\fff \oldstylenums{1}}$\dss and disjoint from\dss $\partial S$\dfdot
A circle in\dss $S$\dss is said to be\qss \emph{non-peripheral}\qss if does not bound 
an annulus together with a component of the boundary\dss $\partial S$\dnsp,\oss 
and\qss \emph{non-trivial}\qss 
if\halfff,\oss in addition,\oss
it does not bound a disc in\dss $S$\dfdot

A circle $D$ in $S$ 
is called\qss \emph{separating}\qss if\qss
$S\smallsetminus D$\dss is not connected.\oss 
Then\dss $S\smallsetminus D$\dss consist of two components.\oss 
The closures of these components are subsurfaces of $S$ having $D$ as a boundary component.\oss
Their other boundary components are included in\qss $\partial S$\dnsp.\oss 

A\qss \emph{bounding pair of circles}\qss 
on a connected\qss \emph{closed}\qss surface $S$ is as an unordered pair\dss 
$C\fff,\pff C'$\dss of disjoint non-isotopic circles in $S$
such that both\sss $C$\sss and\sss $C'$\sss are non-separating\halfff,\oss 
but $S\smallsetminus (C\cup C')$ is not connected.\oss
In this case\qss $S\smallsetminus (C\cup C')$\qss consist of two components.\qss

\myitpar{Twist diffeomorphisms and Dehn twists.} Let $A$ be an annulus,\oss
i.e. a surface diffeomorphic to 
$S^{\fff \oldstylenums{1}}
\times 
[\fff {\fff \oldstylenums{0}},\dff {\fff \oldstylenums{1}}\fff]$\nnsp.\oss
As is well known,\oss the group of diffeomorphisms of\dss $A$\halfff\dfcom
fixed in a neighborhood of\dss $\partial A$\dss and considered up to isotopies 
fixed in a neighborhood\dss $\partial A$\halfff\dfcom
is an infinite cyclic group.\oss
A diffeomorphism of\dss $A$\dss is called a\qss \emph{twist diffeomorphism}\qss of\dss $A$\dss
if it is fixed in a neighborhood of\dss $\partial A$\dss and its isotopy class 
is a generator of this group.\oss 
 
If an annulus\dss $A$\dss is a subsurface of a surface\dss $S$\dnsp,\oss
then any twist diffeomorphism of $A$ can be extended by 
the identity to a diffeomorphism of\qss $S$\dnsp.\oss
Such extensions are called\qss \emph{twist diffeomorphisms}\qss of\qss $S$\dnsp,\oss 
and their isotopy classes are called\qss \emph{Dehn twists}.\oss
The Dehn twist resulting from a twist diffeomorphism of an annulus\dss $A$\dss
in\dss $S$\dss is said to be a\qss \emph{Dehn twist about}\qss a circle\dss $C$\dss in\dss $S$\dss
if\trs $C$\trs is contained in\dss $A$\dss as a deformation retract\halfff.\oss

\myitpar{Left and right twists.} Since\dss $S$\dss is oriented,\oss
every annulus\dss $A$\dss contained in\dss $S$\dss is also oriented.\oss
The orientation of an annulus\dss $A$\dss allows 
to choose a preferred isotopy class of twist diffeomorphisms of\dss $A$\dnsp.\oss
The twist diffeomorphisms in this isotopy class and their extensions to
twist diffeomorphisms of\dss $S$\dss are called the\qss \emph{left twist diffeomorphisms},\oss
and other twist diffeomorphisms are called the\qss \emph{right twist diffeomorphisms}.\oss

The isotopy classes of the left or right twist diffeomorphisms about\dss $C$\dss are called,\oss
respectively,\oss \emph{left}\pss or\pss \emph{right Dehn twists}\qss about\dss $C$\dnsp.\oss
They are uniquely determined by\dss $C$\dnsp.\oss
The left Dehn twist about\dss $C$\dss is denoted by\pss $t_{\dff C}$\nsp.\oss
The right Dehn twist about\dss $C$\dss is the inverse of the left one
and hence is equal to\pss $t_{\dff C}^{\dff -\dff \oldstylenums{1}}$\nsp.\oss
Let\qss $G$\qss be a diffeomorphism of\qss $S$\qss and 
let\qss $g\dff\in\dff \mms$\qss be its isotopy class.\oss
Then\oss
\[
\quad
g\qff t_{\dff C}\qff g^{\fff -\dff \oldstylenums{1}}
\off =\off 
t_{\dff G\fff(C)}\fff.
\]
In particular\halfff,\oss if\oss $G\fff(C)\off =\off C$\nnsp,\oss
or\qss if\qss $G\fff(C)$\qss is isotopic to\qss $C$\snsp,\oss
then\qss $t_{\dff C}$\qss and\qss $g$\qss commute.\oss

\myitpar{The action of Dehn twists on homology.}\qss 
Since\dss $S$\dss is oriented,\oss there is a canonical 
skew-symmetric pairing on $H\one(S)$\dnsp,\oss 
known as the\qss \emph{intersection\dss pairing}.\oss
We denote it by
\[
\quad
(a\fff,\qff b)
\off \longmapsto\off 
\langle\qff a\fff,\qff b \qff\rangle\dff.
\]
Let $C$ be a circle on\dss $S$\dfdot
Let us choose an orientation of\dss $C$\nnsp.\oss 
Let\dss $\hclass{C}$\dss
be the image in\dss $H\one(S)$\dss of the fundamental class of\dss $C$\dss
with this orientation.\oss
Then\qss $t_{\dff C}$\qss acts on\qss $H\one(S)$\qss by the formula\oss\vspace*{3pt}
\begin{equation*}
\quad
\bigl(t_{\dff C}\bigr)_*\fff(a)
\off =\off 
a\qff +\qff  
\langle\qff a\fff,\qff \hclass{C} \qff\rangle\qff\hclass{C}\dff,
\end{equation*}

\vspace*{-9pt}
and the powers of\qss $t_{\fff C}$\qss act by the formula\oss\vspace*{3pt} 
\begin{equation*}
\quad
\left(t_{\dff C}^{\dff m_{\phantom{O}}}\hspace*{-0.5em}\right)_*\fff(a)\off 
=\off a\qff +\qff m\dff 
\langle\qff a\fff,\qff \hclass{C} \qff\rangle\qff\hclass{C}\dff.
\end{equation*}

\vspace*{-9pt}
Changing the orientation of $C$ replaces $\hclass{C}$ by $- \fff\hclass{C}$
and hence does not change the right hand sides of these formulas.\oss
The formula for the action of\qss $t_{\dff C}$\qss implies that\qss 
$t_{\fff C}$\qss belongs to the Torelli group\dss $\tors$\dss
if and only if\oss 
$\hclass{C}\off =\off \oold$\nnsp,\oss
i.e.\qss if\dss and\sss only\sss if\qss $C$\qss is a separating circle.\oss

\myitpar{Multi-twists.} Let\dss $c$\dss be a one-dimensional closed submanifold of\qss $S$\dfdot
A\qss \emph{Dehn multi-twist about}\qss about\dss $c$\dss 
is defined as a product\qss $t$\qss of the form\vspace*{3pt}
\begin{equation*}
\quad
t\off =\off \prod\nolimits_{\fff O}\qff t_{\dff O}^{\dff m_{\dff O}},
\end{equation*}

\vspace*{-9pt}
where\qss $O$\qss runs over all components of\dss $c$\dss
and\dss $m_{\dff O}$\dss are integers.\oss
Every Dehn multi-twist about\dss $c$\dss can be represented by a diffeomorphism of\dss $S$\dss
equal to the identity outside of a subsurface of\dss $S$\dss 
containing\dss $c$\dss as a deformation retract\halfff.\oss
Such diffeomorphisms are called\qss \emph{multi-twist diffeomorphisms}\qss about\dss $c$\dnsp.\oss

Suppose that\qss $C\fff,\pff D$\qss is a bounding pair of circles.\oss 
Then\dss $C$\dss and\dss $D$\dss can be oriented 
in such a way that\oss 
$\hclass{C}\off =\off \hclass{D}$\oss
and hence\qss
$t_{\dff C}$\qss and\qss $t_{\dff D}$\qss induce the same maps on\qss $H\one(S)$\dnsp.\oss
Therefore the multi-twists\oss 
$\displaystyle
t_{\dff C}\qff t_{\dff D}^{\dff -\dff \oldstylenums{1}}\fff,\off\off 
t_{\dff D}\qff t_{\dff C}^{\dff -\dff \oldstylenums{1}}$\qss
belong to\qss $\tors$\dnsp.\oss
Both of them are called the\pss 
\emph{Dehn--Johnson twists}\pss about the bounding pair\qss $C\fff,\pff D$\dnsp.\oss

\mysection{Pure\qss diffeomorphisms\qss and\qss reduction\qss systems}{pure-reduction}

\vspace*{\bigskipamount}
This section is devoted to an overview of the main notions related to
reduction systems of diffeomorphisms of a surface\dss $S$\dss and
elements of\qss $\mms$\dnsp.\oss
This material is needed only for the next section,\oss 
which is devoted to the motivation behind the main results of this paper\halfff,\oss
but is not needed for the proofs of these results.

\myitpar{Cutting surfaces and diffeomorphisms.}
Let\dss $c$\dss be one-dimensional closed submanifold of a surface $S$\nnsp.\oss 
We denote by\qss $S\cutc$\qss
the result of cutting\qss $S$\qss along\qss $c$\dfdot 
The components of\qss $S\cutc$\qss are called\qss \emph{parts}\qss
into which $c$ divides\qss $S$\nnsp.\oss
The canonical map\oss 
\[
\quad
p\cutc\dff\colon\dff S\cutc\off\toto\off S
\] 
induces a diffeomorphism\oss
$\displaystyle
(p\cutc)^{\dff -\dff \oldstylenums{1}}\dff(S\smallsetminus c)\ttoo S\smallsetminus c$\nnsp,\oss 
treated as an identification.\oss
If\qss ${p\cutc}$\qss is injective on a component\dss $Q$\dss of\oss $S\cutc$\nnsp,\qff\oss 
then the induced map\qss 
$\displaystyle
Q\ttoo p\cutc\dff(Q)$\qss 
is also treated as an identification
and\dss $Q$\dss is treated as a subsurface of\dss $S$\nnsp.\oss 

Any diffeomorphism\dss $\psi\dff \colon\dff S\toto S$\dss such that\qss 
$\psi\fff(c)\qff =\qff c$\qss induces a diffeomorphism 
\[
\quad
\psi\cutc\dff \colon\dff S\cutc \ttoo S\cutc\dff.
\] 
If\pss $\psi\cutc$\qss leaves a component\dss $Q$\dss 
of\qss $S\cutc$\qss 
invariant\halfff,\oss 
then\qss $\psi\cutc$\qss induces a diffeomorphism
\[
\quad
\psi_Q\dff\colon Q\ttoo Q\fff,
\]
called the\pss \emph{restriction}\pss of\qss $\psi$\qss to\qss $Q$\nnsp.\oss

\myitpar{Systems of circles and reduction systems.} 
A one-dimensional closed submanifold\dss $c$\dss 
of a surface\dss $S$\dss is called a\qss 
\emph{system of circles on}\dss $S$\qss if the components of\dss $c$\dss are all
non-trivial circles on\dss $S$\dss and are pair-wise non-isotopic.\oss

A system of circles\dss $c$\dss on\dss $S$\dss 
is called a\qss \emph{reduction system}\qss for a diffeomorphism\dss $\psi$\dss of\qss $S$\qss
if\pss $\psi\fff(c)\qff =\qff c$\nnsp.\qff\oss
A system of circles\dss $c$\dss on\dss $S$\dss 
is called a\qss \emph{reduction system}\qss for an element\qss $f\dff\in\dff \mms$\oss
if\qss $c$\trs is a reduction system for some diffeomorphism\dss $\psi$\dss 
in the isotopy class\qss $f$\snsp,\oss
i.e.\qss if\dss $f$\dss can be represented by a diffeomorphism\dss
$\psi$\dss of\trs $S$\trs such that\pss $\psi\fff(c)\qff =\qff c$\nnsp.\oss

\myitpar{Reducible and pseudo-Anosov elements.}
A non-trivial\qss element\qss $f\dff\in \mmod(S)$\qss is said to be\qss \emph{reducible}\qss 
if\dss there exists a non-empty reduction system for\qss $f$\halfff\dnsp,\oss
and\qss \emph{irreducible}\qss otherwise.\qss
An irreducible element of infinite order is called a\qss \emph{pseudo-Anosov element}.\oss
Thurston's theory\qss \cite{thurston}\qss provides a lot of information
about pseudo-Anosov elements.\oss

\myitpar{Pure diffeomorphisms and elements.} 
Let\dss $\psi$\dss be a diffeomorphism  of\dss $S$\nnsp.\oss
A system of circles\dss $c$\dss is said to be a\qss
\emph{pure reduction system}\qss for\dss $\psi$\dss if\dss $c$\dss 
is a reduction system for\dss $\psi$\dss 
and the following four conditions hold.\vspace*{-9pt}
\begin{itemize}
\item[(a)] \dnsp$\psi$\dss is orientation-preserving\fff.\oss 
\item[({\fff}b)] Every component of\qss $S\cutc$\qss is invariant under\trs $\psi\cutc$\nnsp.\oss  
\item[(c)] \dnsp$\psi$\dss is equal to the identity 
           in a neighborhood of\dss $c\dff \cup\dff \partial S$\nnsp.\oss
\item[(d)] For each component\trs $Q$\trs of\qss $S\cutc$\qss the isotopy class of 
the restriction\qss $\psi_Q\dff \colon\dff Q\toto Q$\qss
is either pseudo-Anosov,\oss or\dss contains\qss $\id_Q$\snsp.\oss
\end{itemize}

\vspace*{-9pt}
A diffeomorphism\dss $\psi$\dss of\dss $S$\dss is said to be\pss \emph{pure}\oss 
if\qss $\psi$\qss admits a pure reduction system.\oss
An isotopy class\dss $f\dff \in\dff \mmod(S)$\dss is\dss said\dss to\dss be\pss \emph{pure}\fff\oss 
if\fff\qss $f$\qss contains\dss a\dss pure\dss diffeomorphism.\oss
A system of circles\dss $c$\dss is said to be a\qss
\emph{pure reduction system}\qss for an element\qss $f\dff\in\dff \mms$\qss
if\dss $c$\dss is a pure reduction system for some diffeomorphism\dss $\psi$\dss
in the isotopy class\qss $f$\snsp.\oss
The\qss \emph{isotopy extension theorem}\qss implies that the property of being a pure
reduction system for\qss $f$\qss depends only 
on the isotopy class of the submanifold\dss $c$\dss in\qss $S$\dnsp.\oss

\myitpar{The subgroups\qss $\ttt_{\dnsp m\fff}(S)$\dnsp.}
Given an integer\dss $m$\snsp,\oss the group\dss $\ttt_{\dnsp m\fff}(S)$\qss 
is defined as the subgroup of\qss $\mms$\qss
consisting of the isotopy classes of diffeomorphisms acting trivially on\qss
$H\one(S\fff,\pff \zzz/m\dff\zzz)$\dnsp.\oss
Clearly,\qss $\tors$\qss is contained in\dss 
$\ttt_{\dnsp m\fff}(S)$\pss for every\qss $m\dff \in\dff \zzz$\snsp.\oss
In the present paper the groups\dss $\ttt_{\dnsp m\fff}(S)$\qss
play only a technical role and only in Section\qss \old{\ref{extensions}}\qss
devoted to the motivation.\oss

\mypar{Theorem.}{pure-torelli} \emph{If\pss $m\qff \geqslant\qff \old{3}$\snsp,\qff\oss
then all elements of\pss $\ttt_{\dnsp m\fff}(S)$\pss are pure.}

\mypar{Lemma.}{pure-abc} \emph{Suppose that\dss $c$\dss 
is a reduction system for a diffeomorphism\qss $\psi$\qss
representing an element of\qss $\ttt_{\dnsp m\fff}(S)$\dnsp,\oss 
where\pss $m\qff \geqslant\qff \old{3}$\nnsp.\qff\oss
Then\qss $\psi$\pss leaves every component of\qss $c$\qss invariant\halfff,\oss
and\qss $\psi\cutc$\pss leaves every component of\oss $S\cutc$\qss invariant\halfff.\oss}

\prooftitle{Proofs}\qss See\qss Theorem\qss \old{1}.\old{7}\qss and\qss Theorem\qss \old{1}.\old{2}\qss 
of\oss \cite{i-book}\oss respectively.\oss  \eproof

\myitpar{Minimal pure reduction systems.} A system of circles\dss $c$\dss is said to be a\qss
\emph{minimal pure reduction system}\qss for an element\qss $f\dff\in\dff \mms$\qss if\dss
$c$\dss is a pure reduction system for\qss $f$\snsp,\oss but no proper subsystem of\dss $c$\dss is.\oss
Any pure reduction system for\qss $f$\qss contains a minimal pure reduction system,\oss
and,\oss in fact\halfff,\oss it is unique.\oss
Moreover\halfff,\oss up to isotopy such a minimal pure reduction system 
depends only on\qss $f$\snsp.\oss 

In fact\halfff,\oss the set of the isotopy classes of 
components of a minimal pure reduction system for\qss $f$\qss
is nothing else but the canonical reduction system of\qss $f$\qss 
in the sense of\qss \cite{i-book}.\oss 
This easily follows from the results of\oss \cite{i-book}\halfff,\oss Chapter \old{7}.\oss
See also\oss \cite{im},\oss Section\qss \old{3}.\oss
Since the canonical reduction system of\qss $f$\qss 
is defined invariantly in terms of\qss $f$\snsp,\oss
it depends only on\qss $f$\qss and hence the same is true 
for the minimal pure reduction systems for\qss $f$\snsp.\oss

\mysection{The\qss sirens\qss of\qss Torelli\qss topology}{extensions}

\vspace*{\medskipamount}
\myitpar{Extensions by the identity.} 
Suppose that\qss $Q$\qss is a subsurface of\dss $S$\nnsp.\off\oss 
If\dss $\varphi$\dss is a diffeomorphism of\dss $Q$\dss equal to the identity 
in a neighborhood of the boundary\dss $\partial Q$\nnsp,\oss
then\dss $\varphi$\dss canonically extends to a diffeomorphism of\dss $S$\dss
equal to the identity on the complementary surface\dss $\ccomp Q$\nnsp.\oss 
We will denote this extension by\dss $\varphi\backslash S$\ffdot 

More generally\nsp,\oss let\dss $c$\dss be a one-dimensional closed submanifold of\dss $S$\dss
and let\dss $Q$\dss be a component of\dss $S\cutc$\dfdot
The canonical map\qss $p\cutc\dff \colon\dff S\cutc \toto S$\qss 
is injective on\qss $Q\dff \smallsetminus\dff \partial Q$\snsp,\oss
but may map two different components of\dss $\partial Q$\dss 
onto the same component of\dss $c$\nnsp.\oss
The image\qss $p\cutc\fff(Q)$\trs is always a subsurface of\dss $S$\nnsp.\oss
If\dss $\varphi$\dss is a diffeomorphism of\dss $Q$\dss 
equal to the identity in a neighborhood of\dss $\partial Q$\nnsp,\oss
then\dss $\varphi$\dss induces a diffeomorphism of the image\dss $p\cutc\fff(Q)$\dnsp.\oss
The latter\halfff,\oss in turn,\oss canonically extends to a diffeomorphism of\dss $S$\dss
equal to the identity on\qss $\ccomp\dff \bigl(\fff p\cutc(Q)\bigr)$\dnsp.\oss 
We will denote this extension also by $\varphi\backslash S$\nnsp.\oss

\myitpar{A construction of abelian subgroups of Teichm\"{u}ller modular groups.}
While the extensions by the identity are rarely mentioned explicitly\nsp,\oss
they play a crucial role in the background of
the theory of Teichm\"{u}ller modular groups.\oss
As an example,\oss let us outline their role in 
the classification of abelian subgroups
of Teichm\"{u}ller modular groups.\oss

Let\dss $c$\dss be a one-dimensional closed submanifold of\dss $S$\snsp.\oss
Suppose that for each component\dss $Q$\dss of\qss $S\cutc$\qss
we are given a diffeomorphism\qss $\varphi(Q)\dff \colon\dff Q\toto Q$\qss
equal to identity in a neighborhood of the boundary\dss $\partial Q$\dnsp.\oss
Extensions\qss $\varphi(Q)\backslash S$\qss of these diffeomorphisms by the identity to 
diffeomorphisms\qss $S\dff\toto\dff S$\qss
obviously commute.\oss
Therefore,\oss the isotopy classes of these extensions
generate an abelian subgroup of\qss $\mms$\dnsp.\oss

It turns out that every abelian subgroup of\qss $\mmod(S)$\qss consisting
of pure elements is contained in an abelian subgroup of this form.\oss
Together with Theorem\qss \ref{pure-torelli}\qss this implies that
for every subgroup\qss $\mathfrak{G}$\qss of\qss $\mmod(S)$\qss 
the intersection\qss $\mathfrak{G}\dff\cap\dff \ttt_{\dnsp m\fff}(S)$\dnsp,\oss
which is a subgroup of finite index in\qss $\mathfrak{G}$\dnsp,\oss
is contained in a subgroup of this form.\oss

\myitpar{A better construction of abelian subgroups.}
It takes into account the special role of Dehn multi-twists.\oss
The extensions\qss $\varphi(Q)\backslash S$\qss 
are equal to the identity in a neighborhood of\dss $c$\dnsp.\oss 
Therefore,\oss these extensions commute with twist diffeomorphisms
representing Dehn twists about components of\dss $c$\dnsp,\oss
and hence the isotopy classes of these extensions together with multi-twists
about\dss $c$\dss generate an abelian subgroup\qss $\mathfrak{A}$\qss of\qss $\mmod(S)$\dnsp.\oss

Adding Dehn multi-twists generators can be replaced by adding to\dss $c$\dss
components isotopic to the already present ones.\oss
But if we allow adding Dehn multi-twists,\oss 
we may assume that\dss $c$\dss is a system of circles\qss
and that for each\dss $Q$\dss the diffeomorphism\qss 
$\varphi(Q)$\qss is either isotopic to\qss $\id_Q$\dnsp,\oss
or belongs to a pseudo-Anosov isotopy class.\oss

The resulting class of the abelian subgroups is the same as the original one.\oss
But now the constructed subgroup admits a more detailed description.\oss
Without changing the subgroup,\oss
we can put aside diffeomorphisms\qss
$\varphi(Q)$\qss which are isotopic to the identity.\oss
The isotopy classes of extensions of remaining diffeomorphisms
together with Dehn twists about the components of\dss $c$\dss 
turn out to be the free generators of\qss $\mathfrak{A}$\dnsp.\oss

\myitpar{Pasting.} This construction of abelian subgroups 
is a\qss \emph{pasting argument}\qss par excellence.\oss
We start with a system of circles\dss $c$\dss on our surface\dss $S$\dnsp.\oss
It subdivides\dss $S$\dss into several pieces.\oss
For each piece we consider an object located on it\qss 
(a diffeomorphism of the piece in our example)\qss
and paste them together into a new object on the whole surface\dss $S$\dss
(an abelian group in our example).\oss
Usually one needs to add something located,\oss at least morally,\oss 
on the cutting submanifold\dss $c$\dss ({\dff}Dehn multi-twists about\dss $c$\nsp).\oss

\myitpar{Commutants\dss and\dss bicommutants.} 
The\qss \emph{commutant}\qss $X'$\qss of a subset\qss 
$X$\qss of a group\qss $\mathfrak{G}$\qss
is the subgroup of\qss $\mathfrak{G}$\qss 
consisting of elements commuting with all\qss $g\dff\in\dff X$\dnsp.\oss
The\qss \emph{bicommutant}\qss of\qss $X$\qss is the commutant of the commutant of\qss $X$\dnsp,\oss
i.e.\dss the subgroup\qss $X''$\dnsp.\oss
The\qss \emph{commutant}\qss $g'$\qss and the\qss \emph{bicommutant}\qss $g''$\qss 
of an element\qss $g\dff\in\dff \mathfrak{G}$\qss are the defined
as the commutant and the bicommutant\dss of\dss 
the one-element subset\qss $\{\dff g\dff\}$\dnsp.\oss

\myitpar{Commutants\dss and\dss bicommutants\dss in\qss \emph{$\mmod(S)$}\dnsp.}
Let us fix a subgroup\qss $\Gamma$\qss of finite index in\qss $\mmod(S)$\qss
consisting of pure elements.\oss
For example,\oss by Theorem\qss \ref{pure-torelli}\qss we can take as\qss $\Gamma$\oss
the subgroup\dss $\ttt_{\dnsp m\fff}(S)$\qss for any\qss $m\qff \geqslant\qff 3$\nnsp.\oss

Suppose that\dss $f\dff \in\dff \Gamma$\dss and\dss $c$\dss 
is a pure reduction system for\dss $f$\nnsp.\oss
Then\dss $f$\dss can be represented by a diffeomorphism\qss 
$\varphi\dff \colon\dff S\qff\toto\qff S$\qss such that\dss $c$\dss 
is a pure reduction system for\dss $\varphi$\nnsp.\oss
Let us consider for every component\dss $Q$\dss of\pss $S\cutc$\pss the restriction\qss 
\[
\quad
\varphi_Q\dff \colon\dff Q\ttoo Q\dff,
\]
i.e.\qss the diffeomorphism\qss $Q\qff\toto\qff Q$\qss induced by\dss $\varphi\cutc$\nnsp.\oss
The restriction\dss $\varphi_Q$\dss is equal to the identity 
in a neighborhood of\qss $\partial Q$\snsp,\oss
and hence the extension\qss $\varphi_Q\nnsp\backslash S$\qss of\qss 
$\varphi_Q$\qss to a diffeomorphism of\dss $S$\dss is defined.\oss 

By applying the above construction to\dss $c$\dss 
and diffeomorphisms\qss $\varphi_Q\nnsp\backslash S$\nnsp,\oss 
we get an abelian subgroup of\qss $\mmod(S)$\dnsp.\oss
Let\qss $\mathfrak{B}(f)$\qss be the intersection of this subgroup with\qss $\Gamma$\dnsp.\oss
It turns out that\qss $\mathfrak{B}\fff(f)$\qss is a subgroup of finite index in 
the bicommutant\qss $f''$\qss of\qss $f$\qss in\qss $\Gamma$\dnsp.\oss 
The bicommutant\qss $f''$\qss itself is equal to\qss $\mathfrak{B}\fff(g)$\qss
for an element\qss $g$\qss closely related to\qss $f$\dnsp.\oss

The\qss (somewhat disguised)\qss commutants,\oss bicommutants,\oss
and groups\qss $\mathfrak{B}\fff(f)$\qss are the key tool  
for algebraic characterizations of Dehn twists and related elements 
of\qss $\mmod(S)$\qss
and for the classification of automorphisms of $\mmod(S)$\dnsp.\oss 
See\oss \cite{i-aut},\oss \cite{i-book}.\oss
Ad hoc analogues of these construction for Torelli groups
were used by B. Farb and the author\oss \cite{fi-announ}\oss
in order to give an algebraic characterization of 
Dehn and Dehn--Johnson twists in\qss $\tors$\dnsp.\oss

\myitpar{Cut-and-paste.}\qss The construction 
of the groups\qss $\mathfrak{B}\fff(f)$\qss is a\qss 
\emph{cut-and-paste argument}\qss par excellence.\oss
We started with a pure isotopy class\dss $f$\nnsp.\oss
It has a representative\dss $\varphi$\dss which can be 
cut into sufficiently well understood pieces\qss $\varphi_Q$\nnsp.\oss

Namely\nsp,\oss the isotopy class\dss $f_Q$\dss of each piece\trs 
$\varphi_Q$\trs is either pseudo-Anosov,\oss
or contains the identity.\oss
In both cases the commutant and the bicommutant of\dss $f_Q$\dss 
in\dss $\mm(Q)$\dss are well understood.\oss
The case when\dss $f_Q$\dss contains identity is trivial.\oss
If\dss $f_Q$\dss is pseudo-Anosov,\oss then both\dss $f_Q'$\dss and\dss $f_Q''$\dss
contain an infinite cyclic group generated by a pseudo-Anosov element
as a subgroup of finite index.\oss
See\oss \cite{mc-nc}\oss or\oss \cite{i-book}.\oss

After cutting\dss $\varphi$\dss into pieces\qss $\varphi_Q$\nsp,\oss 
we considered the extensions\qss $\varphi_Q\nnsp\backslash S$\qss 
of these pieces to\qss $S$\qss and pasted them together 
into the group\qss $\mathfrak{B}\fff(f)$\dnsp.\oss

\myitpar{A siren song.}\qss The sirens sing that one can 
cut and paste in Torelli groups also.\oss
A. Putman begins his paper\oss \cite{p}\oss as follows.\oss
\begin{quote}
We introduce machinery to allow ``cut-and-paste''-style 
inductive arguments in the Torelli subgroup\qss \ldots .\qss  
In the past these arguments have been problematic because restricting 
the Torelli group to subsurfaces gives 
different groups depending on how the subsurfaces are embedded.\oss
We define a category \ldots
\end{quote}
While Putman's machinery allows to use ``cut-and-paste'' arguments
in some situations,\oss they remain problematic in other ones.\oss
The main obstacle is the lack of a suitable definition of Torelli groups
for surfaces with boundary.\oss
By this reason after cutting\dss $\varphi$\dss into restrictions\dss $\varphi_Q$\dss
we are not in the theory of\dss Torelli groups anymore.\qss

In fact{\halfff},\oss the heart of Putman's paper\oss \cite{p}\oss is a definition of
Torelli groups for surfaces with boundary.\oss
But he does not address the question if 
the isotopy classes of restrictions\dss 
$\varphi_Q$\dss belong to his version of\dss Torelli groups.\oss 
The siren sings that they do belong,\oss 
and one only has to read\oss \cite{p}\qss carefully \ldots

\myitpar{A bicommutant siren song.}\qss 
Suppose that\dss $S$\dss is a closed surface and\qss $f\dff\in\dff\tors$\dnsp.\oss
It is tempting to think that one construct a subgroup of finite index 
in the bicommutant\qss $f''$\qss of\qss $f$\qss in\qss $\tors$\qss by
adapting the construction of\qss $\mathfrak{B}\fff(f)$\dnsp.\oss 

The element\dss $f$\dss is a pure element by Theorem\qss \ref{pure-torelli}.\oss
Let $c$ be a pure reduction system for a diffeomorphism\dss $\varphi$\dss of\qss $S$\qss
representing\dss $f$\nnsp.\oss
Let us consider{\halfff},\oss 
for each component\dss $Q$\dss of\dss $S\cutc$\nsp,\qff\oss 
the restriction\dss $\varphi_Q$\dss  
and the extension\qss $\varphi_Q\nnsp\backslash S$\pss of\qss $\varphi_Q$\snsp.\oss 

The sirens sing that the group generated by the
Dehn multi-twists about\dss $c$\dss belonging to\dss $\tors$\dss 
together with the isotopy classes of extensions\qss 
$\varphi_Q\nnsp\backslash S$\pss is a subgroup of finite index 
in the bicommutant\qss $f''$\qss of\qss $f$\qss in\qss $\tors$\dnsp.\oss

\myitpar{An extension siren song.}\qss
The sirens sing that for each component\dss $Q$\dss of\qss $S\cutc$\qss
the isotopy class\qss $f_Q\nsp\backslash S$\qss
of the extension\qss $\varphi_Q\nnsp\backslash S$\qss
belongs to\dss $\tors$\dnsp,\oss
at least up to Dehn multi-twists about\dss $c$\dnsp.\oss
Otherwise the bicommutant song does not make much sense.\oss

Since only the component\dss $Q$\dss of\dss $S\cutc$\dss
affects\oss $f_Q\nsp\backslash S$\nnsp,\oss
one may forget about other components.\oss
This leads to a plain extension version of the bicommutant song\halfff.\oss

Let\qss $f\dff\in\dff \tors$\dnsp,\oss
and let\qss $\varphi$\qss be a representative of\qss $f$\qss
leaving a subsurface\dss $Q$\dss of\dss $S$\dss invariant{\halfff}.\oss
Without any loss of generality one can assume that\qss 
$\varphi$\qss is equal to the identity in a neighborhood of\dss $\partial Q$\dss in\dss $S$\nnsp.\oss
Consider the restriction\dss $\varphi_Q$\dss and its extension\dss $\varphi_Q\nnsp\backslash S$\snsp.\oss

Let\oss $f_Q\nsp\backslash S$\oss be the isotopy class of\qss $\varphi_Q\nnsp\backslash S$\nnsp.\oss
The sirens sing that\oss $f_Q\nsp\backslash S$\oss belongs to\qss 
$\tors$\qss up to Dehn multi-twist about\qss $\partial Q$\nnsp.\oss
In other words,\oss
\[
\quad
t\dff \cdot\dff \bigl(f_Q\nsp\backslash S\bigr)\qff\in\qff \tors
\] 
for some Dehn multi-twist $t$ about $\partial Q$\nnsp.\qff\oss
Let us return to the reality.\oss

\myitpar{Example.} Suppose that $S$ is a closed surface
and $Q$ is a subsurface of $S$ such that both $Q$ 
and its complementary surface $\ccomp Q$ are connected.\oss
Suppose that the boundary $\partial Q$ consists of four components
$C\one\fff,\pff C\two\fff,\pff C\three\fff,\pff C\four$\nnsp.\oss
Let\dss $C$\dss be a circle bounding in $Q$ a disc with two holes 
together with\dss $C\one$\dss and\dss $C\two$\nnsp,\oss
and\qss $D$\qss be a circle bounding in\qss $\ccomp Q$\qss a disc with two holes 
together with the same circles\dss $C\one$\dss and\dss $C\two$\nnsp.\oss
The union of these two discs with two holes is a genus $\old{1}$ surface
with the boundary $C\cup D$\dfcom
and hence\qss\vspace*{3pt}
\[
\quad
f\off =\off t_{\dff C}\qff t_{\dff D}^{\dff -\dff \oldstylenums{1}}
\]

\vspace*{-9pt}
is a Dehn-Johnson twist{\halfff}.\oss
In particular{\halfff},
$f$\qss
belongs to the Torelli group $\tors$\dfdot

Let\qss $\tau_{\dff C}\dff,\off \tau_{\dff D}$\qss be twist diffeomorphisms
of\dss $S$\dss representing\qss $t_{\dff C}\dff,\off t_{\dff D}$\qss 
and equal to the identity on\qss $\ccomp Q\dff,\off Q$\qss
re\-spec\-tive\-ly\nsp.\qff\oss
By the definition,\oss
$\displaystyle 
f
\off =\off 
t_{\dff C}\qff t_{\dff D}^{\dff -\dff \oldstylenums{1}}$\oss
is the isotopy class of\vspace*{3pt}
\[
\quad
\varphi
\off =\off 
\tau_{\dff C}\qff \circ\qff \tau_{\dff D}^{\dff -\dff \oldstylenums{1}}\dff.
\]

\vspace*{-9pt} 
Since $\tau_{\fff D}$\dss is equal to identity on\dss $Q$\nsp,\oss
the restriction\dss $\varphi_Q$\dss is equal to the restriction\dss $(\tau_{\dff C})_{\fff Q}$\dfdot
Since\qss $\tau_{\dff C}$\qss is equal to the identity on $\ccomp Q$\dfcom
the extension\qss $\tau_{\dff C}\hnsp\backslash S$\qss is equal to\qss $\tau_{\dff C}$\nsp.\off\oss
Hence\qss\vspace*{3pt} 
\[
\quad
\varphi_Q\hnsp\backslash S
\off =\off 
\tau_{\dff C}
\]

\vspace*{-9pt}
and the isotopy class\qss $f_Q\nsp\backslash S$\qss of\qss $\varphi_Q\hnsp\backslash S$\qss
is equal to\qss $t_{\fff C}$\nnsp.\oss 
Since\qss $C$\qss is non-separating in\qss $S$\dnsp,\oss\vspace*{3pt}
\[
\quad
f_Q\nsp\backslash S
\off =\off 
t_{\fff C}\pff \not\in\pff \tors\dff.
\]

This ruins the most optimistic extensions hopes.\qss
The rest is ruined by the fact that\vspace*{3pt}
\[
\quad
t\cdot \bigl(f_Q\nsp\backslash S\bigr)
\off =\off 
t\cdot t_{\fff C}
\pff \not\in\pff \tors
\] 

\vspace*{-9pt}
for any Dehn multi-twist\qss $t$\qss about\qss $\partial Q$\dnsp.\oss
Let us prove this.\oss

To begin with,\qss let us orient the circles\qss
$C\one\fff,\pff \ldots\fff,\pff C\four$\pss 
as components of\pss $\partial Q$\nnsp.\oss
The only relation between the homology classes\qss 
$\hclass{C\one}\fff,\pff \ldots,\pff \hclass{C\four}$\qss is\vspace*{3pt}
\[
\quad
\hclass{C\one}\qff +\qff \hclass{C\two}\qff +\qff 
\hclass{C\three}\qff +\qff \hclass{C\four}
\off =\off 
\oold\dff.
\]

\vspace*{-9pt}
In particular{\halfff},\oss any three of the classes\qss 
$\hclass{C\one}\fff,\pff \ldots\fff,\pff \hclass{C\four}$\qss 
are linearly independent.\oss
Let us orient\dss $C$\dss and\dss $D$\dss in such a way that\qss 
$\displaystyle \hclass{C}
\off =\off 
\hclass{D}
\off =\off 
\hclass{C\one}\qff +\qff \hclass{C\two}$\nnsp.\oss

There is a circle\dss $A$\dss in\dss $S$\dss such that\dss 
$A$\dss is disjoint from\dss $C\two$\dss and\dss $C\four$\dss and
intersects each of three circles\qss $C$\nnsp,\oss $C\one$\nnsp,\oss 
and\dss $C\three$\dss transversely at one point.\oss
One can orient\dss $A$\dss in such a way that the intersection numbers of
the homology class\dss $a\qff =\qff \hclass{A}$\dss are\vspace*{6pt}
\[
\quad
\langle\qff a\fff,\pff \hclass{C} \qff\rangle
\off\phantom{\two} =\off 
\langle\qff a\fff,\pff \hclass{C\one} \qff\rangle
\off =\off \old{1}\dff,
\quad
\langle\qff a\fff,\pff \hclass{C\three} \qff\rangle
\off =\off 
-\qff \old{1}\dff,
\]

\vspace*{-33pt}
\[
\quad
\langle\qff a\fff,\pff \hclass{C\two} \qff\rangle
\off =\off 
\langle\qff a\fff,\pff \hclass{C\four} \qff\rangle
\off =\off 
\oold\dff.
\]

\vspace*{-6pt}
Every Dehn multi-twist\dss $t$\dss about\qss $\partial Q$\qss has the form\vspace*{3pt}
\[
\quad
t
\off =\off 
\prod\nolimits_{\dff i\qff =\qff \old{1}}^{\dff \old{4}} 
\left(\fff t_{\dff  C_{\dff i}}\right)^{m_{\dff i}},
\]

\vspace*{-9pt}
where\qss $m\one\fff,\pff \ldots,\pff m\four$\qss are integers.\oss
The product\qss $t\cdot t_{\dff C}$\qss acts on\dss $a$\dss as follows\halfff:\vspace*{6pt}
\[
\quad
\bigl(\fff t\cdot t_C\halfff\bigr)_{\fff *}\dff(a)
\off =\off 
a\qff +\qff \hclass{C}\qff +\qff m\one\fff \hclass{C\one}\qff -\qff m\three\fff \hclass{C\three} 
\]

\vspace*{-33pt}
\[
\quad
\phantom{\bigl(t\cdot t_C\bigr)_{\fff *}\dff(a)
\off }
=\off a\qff +\qff \bigl(\fff\hclass{C\one}\qff +\qff \hclass{C\two}\fff\bigr)
\qff +\qff 
m\one\fff \hclass{C\one}\qff -\qff m\three\fff\hclass{C\three}
\]

\vspace*{-33pt}
\[
\quad
\phantom{\bigl(t\cdot t_C\bigr)_{\fff *}\dff(a)\off }
=\off a\qff +\qff (m\one\qff +\qff \old{1})\fff \hclass{C\one}\qff 
+\qff \hclass{C\two}\qff -\qff m\three\fff \hclass{C\three}\off
\]

\vspace*{-33.5pt}
\[
\quad
\phantom{\bigl(t\cdot t_C\bigr)_{\fff *}\dff(a)\off }
\neq\off a\dff,
\]

\vspace*{-6.5pt}
because the classes\qss 
$\hclass{C\one}\fff,\oss \hclass{C\two}\fff,\oss  \hclass{C\three}$\qss
are linearly linearly independent.\oss
It follows that\vspace*{3pt} 
\[
\quad
t\cdot (f_Q\nsp\backslash S)
\off =\off 
t\cdot t_{\dff C}
\qff \not\in\qff
\tors
\] 

\vspace*{-9pt}
for any Dehn multi-twist\dss $t$\dss about\qss $\partial Q$\snsp.\oss

\myitpar{Example: a modification.} In the above example\dss $\partial Q$\dss  
is not a pure reduction system for\dss $f$\dnsp.\qff\oss
Let us modify the example in such a way that it will be.\qff\oss
Let\qss $m\dff\in\dff\zzz$\qss
and\dss let\vspace*{3pt}
\[
\quad
f^{\dff m}
\off =\off  
t_{\dff C}^{\qff m}\cdot t_{\dff D}^{\qff -\dff m}
\off \in\off
\tors
\]

\vspace{-9pt}
be the isotopy class of\oss
$\displaystyle
\varphi^{\dff m}
\off =\off 
\tau_{\dff C}^{\qff m}\qff \circ\qff \tau_{\dff D}^{\qff -\dff m}$\nsp.\oss
Then the isotopy class\pss 
$\displaystyle 
f_Q^{\qff m}\dff\backslash\halfff S$\oss 
of\oss 
$\displaystyle 
\varphi_Q^{\qff m}\dff\backslash\halfff S$\oss
is equal to\qss 
$\displaystyle 
t_{\fff C}^{\dff m}$\nsp.\off\oss
As above,\oss\vspace*{3pt} 
\begin{equation*}
\quad
t\cdot \bigl(f_Q^{\qff m}\dff\backslash\halfff S\bigr)
\off =\off 
t\cdot t_{\dff C}^{\qff m}
\off \not\in \off
\tors
\end{equation*}

\vspace*{-9pt} 
for any Dehn multi-twist\dss $t$\dss about\qss $\partial Q$\dfdot\vspace*{3pt}

Let\qss $\psi$\qss be a diffeomorphism of\qss $S$\qss equal to the identity
in a neighborhood of\qss $\ccomp Q$\dnsp.\oss
Let\qss\vspace*{2pt}
\[
\quad
\otherkappa
\off =\off 
\varphi^{\dff m}\qff \circ\qff \psi\dff.
\]

\vspace*{-10pt}
Then\oss
$\displaystyle
\otherkappa_{\dff Q}
\off =\off  
T_{\fff C}^{\qff m}\qff \circ\qff \psi_Q$\oss
and\oss
$\displaystyle
\otherkappa_{\dff Q}\nsp\backslash\fff S
\off =\off  
T_{\dff C}^{\qff m}\qff \circ\qff \psi$\nnsp.\qff\oss
Let\qss\vspace*{3pt}
\[
\quad
g\fff,\qquad\off 
g_Q\fff,\qquad 
k\fff,\qquad 
k_{\dff Q}\fff,\qquad\qff 
k_{\dff Q}\nsp\backslash\fff S
\]

\vspace*{-9pt}
be the isotopy classes of{\halfff},\oss respectively,\oss
\[
\quad
\psi\fff,\qquad 
\psi_Q\fff,\qquad 
\otherkappa\fff,\qquad 
\otherkappa_{\dff Q}\fff,\qquad 
\otherkappa_{\dff Q}\nsp\backslash\fff S\fff.
\]

\vspace*{-9pt}
Then\off\oss 
$\displaystyle
k
\off =\off 
f^{\dff m}\dff \cdot\dff g$\nnsp,\off\oss
$\displaystyle
k_{\dff Q}
\off =\off 
t_{\dff C}^{\qff m}\dff \cdot\dff g_{\fff Q}$\nsp,\off\oss
$\displaystyle
k_{\dff Q}\nsp\backslash\fff S
\off =\off 
t_{\dff C}^{\qff m}\cdot g$\nnsp,\off\oss
and\qss if\qff\oss 
$g\dff \in\dff \tors$\dnsp,\off\oss 
then\vspace*{6pt}
\[
\quad
k
\off =\off 
f^{\dff m}\dff \cdot\dff g
\off \in\off 
\tors
\hspace*{1.5em}\mbox{ and }\hspace*{1.5em}
t\dff \cdot\dff \bigl(\fff k_{\dff Q}\nsp\backslash\fff S \fff\bigr)
\off =\off 
t\fff\cdot\fff t_{\dff C}^{\qff m}\dff \cdot\dff g
\off \not\in\off
\tors
\]

\vspace*{-6pt}
for any Dehn multi-twist\sss $t$\sss about\sss $\partial Q$\dnsp.\off\oss
If{\halfff},\oss in addition,\oss 
$t_{\dff C}^{\qff m}\dff \cdot\dff g_{\fff Q}
\dff \in\dff
\mmod(Q)$\oss 
is a pseudo-Anosov element{\halfff},\oss
then\qss $k_{\dff Q}$\qss is irreducible.\oss\vspace*{3pt}

Let\qss $\eta$\qss be a diffeomorphism of\qss $S$\qss equal to the identity
in a neighborhood of\qss $Q$\dnsp.\oss
Let\qss\vspace*{3pt}
\[
\quad
\lambda\off =\off \otherkappa\qff \circ\qff \eta\fff.
\]

\vspace*{-9pt}
Then\oss
$\displaystyle
\lambda_{\dff Q}
\off =\off  
\otherkappa_{\dff Q}$\dnsp,\off\oss
$\displaystyle
\lambda_{\dff Q}\nsp\backslash\fff S
\off =\off  
\otherkappa_{\dff Q}\nsp\backslash\fff S$\nnsp,\off\oss
and\oss
$\displaystyle
\lambda_{\dff \ccomp\fff Q}
\off =\off  
\tau_{\dff D}^{\dff -\dff m}\qff \circ\qff \eta_{\dff \ccomp\fff Q}$\nsp.\qff\oss
Let\qss\vspace*{3pt}
\[
\quad
e\fff,\qquad\off 
e_{\dff \ccomp \fff Q}\fff,\qquad 
l\fff,\qquad 
l_{\dff Q}\fff,\qquad\qff 
l_{\dff \ccomp \fff Q}
\]

\vspace*{-9pt}
be the isotopy classes of{\halfff},\oss respectively,\oss\vspace*{3pt}
\[
\quad
\eta\fff,\qquad 
\eta_{\dff \ccomp \fff Q}\fff,\qquad 
\lambda\fff,\qquad 
\lambda_{\fff Q}\fff,\qquad 
\lambda_{\dff \ccomp \fff Q}\fff.
\]

\vspace{-9pt}
Then\off\oss 
$\displaystyle
l
\off =\off 
k\dff \cdot\dff e$\dnsp,\off\oss
$\displaystyle
l_{\dff Q}
\off =\off 
k_{\dff Q}$\nsp,\off\oss
$\displaystyle
l_{\trf \ccomp \fff Q}
\off =\off 
t_{\dff D}^{\dff -\dff m}\dff \cdot\dff e_{\dff \ccomp\fff Q}$\snsp,\off\oss
and\qss if\qff\oss 
$g\fff,\pff e\dff \in\dff \tors$\dnsp,\off\oss 
then\vspace*{6pt}
\[
\quad
l
\off =\off 
k\dff \cdot\dff e\off\in\off \tors
\hspace*{1.5em}\mbox{ and }\hspace*{1.5em}
t\dff \cdot\dff \bigl(\fff l_{\dff Q}\nsp\backslash\fff S \fff\bigr)
\off =\off 
t\dff \cdot\dff \bigl(\fff k_{\dff Q}\nsp\backslash\fff S \fff\bigr)
\off \not\in \off
\tors
\]

\vspace{-6pt}
for any Dehn multi-twist\sss $t$\sss about\sss $\partial Q$\nnsp.\off\oss
If{\halfff},\oss in addition,\oss 
$t_{\dff D}^{\dff -\dff m}\dff \cdot\dff e_{\dff \ccomp\fff Q}\dff\in\dff\mmod(\ccomp Q)$\oss 
is a pseudo-Anosov element{\halfff},\oss
then\qss $l_{\trf \ccomp \fff Q}$\qss is irreducible.\oss\vspace*{6pt} 

To sum up,\oss if\oss 
$g\fff,\pff e
\qff \in\qff \tors$\oss 
and if the isotopy classes\oss
$t_{\dff C}^{\qff m}\dff \cdot\dff g_{\fff Q}$\nnsp,\off\oss
$t_{\dff D}^{\dff -\dff m}\dff \cdot\dff e_{\dff \ccomp\fff Q}$\oss
are pseudo-Anosov,\oss then\vspace*{3pt}
\[
\quad
l
\off \in\off 
\tors\fff,\qquad
t\dff \cdot\dff \bigl(l_{\dff Q}\nsp\backslash\fff S\bigr)
\off \not\in \off
\tors\fff,
\]

\vspace*{-6pt}
and by Lemma\qss \ref{pure-abc}\qss $\partial Q$\qss 
is contained in any pure reduction system for\dss $l$\snsp.\oss
Since the isotopy classes of\qss $\lambda_{\dff Q}$\qss and\qss $\lambda_{\dff \ccomp Q}$\qss
are pseudo-Anosov,\oss $\partial Q$\qss is actually a pure reduction system for\dss $l$\snsp.\oss
It remains to construct diffeomorphisms\qss $\psi\fff,\pff \eta$\qss with such properties.\vspace*{3pt}

\myitpar{Construction of\pss $\psi\fff,\pff \eta$\nnsp.}
We will construct\qss $\psi\fff,\pff \eta$\qss under the assumption 
that the genus of both\dss $Q$\dss and\dss $\ccomp Q$\dss is at least\dss $\old{1}$\nnsp.\oss
Then $Q$ contains circles bounding in $Q$ a surface 
of genus\qss $\geqslant\qff \old{1}$\qss with one boundary component\halfff.\oss
Such circles are 
bounding circles in\qss both\qss $Q$\qss and\qss $S$\dnsp.\oss
A well known construction of Thurston
leads to pseudo-Anosov elements of\qss $\mmod(Q)$\qss
of the form\qss $t\one^a\dff \cdot\dff t\two^b$,\off\oss
where both\qss $t\one$\qss and\qss $t\two$\qss are Dehn twists about such bounding circles.\oss

Such an element\qss $t\one^a\dff \cdot\dff t\two^b$\qss can be represented
by a composition of the form\qss $\tau\one^{\dff a}\qff \circ\qff \tau\two^{\dff b}$\nnsp,\oss
where\qss $\tau\one\dff,\qss \tau\two$\qss are twist diffeomorphisms of\qss $Q$\qss
about bounding circles.\oss
Consider now\qss $\tau\one\dff,\qss \tau\two$\qss as twist diffeomorphisms of\dss $S$\dss
about the same circles.\oss 
Then\vspace*{3pt} 
\[
\quad
\psi\off =\off \tau\one^{\dff a}\qff \circ\qff \tau\two^{\dff b}
\]

\vspace*{-9pt}
is a diffeomorphism of\qss $S$\qss equal to the identity on\qss $\ccomp Q$\dnsp.\oss
Since Dehn twists about bounding circles belong to\qss $\tors$\dnsp,\oss
the isotopy class\qss $g$\qss of\qss $\psi$\qss belongs to\qss $\tors$\dnsp.\oss
The isotopy class\qss $g_{\fff Q}$\qss of the restriction\qss $\psi_Q$\qss
is equal to\pss $t\one^a\dff \cdot\dff t\two^b$\pss and hence\qss $g_{\fff Q}$\qss
is pseudo-Anosov.\oss
It follows that\oss $t_{\dff C}^{\qff m}\dff \cdot\dff (g_{\fff Q})^n$\oss 
is pseudo-Anosov for all sufficiently large\qss $m\dff,\qss n$\nsp.\oss
See\oss \cite{i-book},\oss the proof of Lemma\qss 5.8\qss 
in the case when\dss $g$\dss from this Lemma is reducible\qss (not our\dss $g$\dnsp).\dfdot

By the same arguments,\qss there is a diffeomorphism\qss 
$\eta$\qss of\qss $S$\qss such that\qss 
$\eta$\qss is equal to the identity on\qss $Q$\dnsp,\oss
its isotopy class\qss $e$\qss belongs to\qss $\tors$\dnsp,\oss
and the isotopy class\qss $e_{\dff \ccomp\fff Q}$\qss 
of its restriction\qss $\eta_{\dff \ccomp\fff Q}$\qss
is a pseudo-Anosov element of\qss $\mmod(\ccomp Q)$\dnsp.\oss
By the same reasons as above,\oss
$t_{\dff D}^{\dff -\dff m}\dff \cdot\dff (e_{\dff \ccomp Q})^n$\oss
is pseudo-Anosov for all sufficiently large\qss $m\dff,\qss n$\nsp.\oss

Let\qss $m\dff,\qss n$\qss be two sufficiently large integers
and replace\qss $\psi\dff,\off \eta$\qss by their\qss $n$\dnsp-th powers.\oss
Then\qss $m$\qss and the new\qss $\psi\dff,\off \eta$\qss have all the required properties.\qss

\mysection{Weakly\qss Torelli\qss diffeomorphisms}{weakly-torelli-section}

\vspace*{\medskipamount}
\myitpar{Notations related to the singular homology theory.} For a topological space\qss $X$\qss
we denote by\qss $C_{\dff n\fff}(X)$\qss the group of singular\qss $n$\dnsp-chains in\qss $X$\nsp,\oss
and by\qss $\partial$\qss the boundary map of the singular chain complex\qss $C_{\dff *\fff}(X)$\dnsp.\oss
For a continuous map\qss $F\dff\colon X\qff\toto\qff Y$\qss
we denote by\qss $F_{\fff *\fff}$\qss the induced maps\qss
$C_{\dff *\fff}(X)\ttoo C_{\dff *\fff}(Y)$\qss and\qss 
$H_{\dff *\fff}(X)\ttoo H_{\dff *\fff}(Y)$\dnsp.\oss

For a cycle\dss $\alpha$\dss
we denote by\qss $\hclass{\alpha}$\qss the homology class.\oss
The space in question usually is clear from the context\halfff,\oss
but if it is important to distinguish between homology classes of\dss $\alpha$\dss
in different spaces,\oss 
we denote by\qss $\hclass{\alpha}_{\fff X}\dss \in\dff H_{\fff *}(X)$\qss 
the homology class of\dss $\alpha$\dss in\qss $X$\nnsp.\oss

By\qss $Z_{\dff n\fff}(X)$\qss we denote the group of\dss $n$\dnsp-cycles in\dss $X$\nnsp.\oss
For subspace\qss $A$\qss of\qss $X$\qss
we denote by\qss $Z_{\dff n\fff}(X\fff,\pff A)$\qss the group of relative\dss $n$\dnsp-cycles of
the pair\qss $(X\fff,\pff A)$\dnsp,\oss i.e.\qss the group of
singular chains\qss $\alpha\dff\in\dff C_{\dff n\fff}(X)$\qss in\qss $X$\qss 
such that\qss $\partial\alpha\dff\in\dff C_{\dff n\dff -\dff 1\fff}(A)$\dnsp.\oss
Then\oss\vspace*{3pt}
\begin{equation*}
\quad
H_{\dff n\fff}(X\fff,\pff A)\off 
=\off Z_{\dff n\fff}(X\fff,\pff A)\bigl/\bigl(\fff\partial C_{\dff n\dff +\dff 1\fff}(X)\qff 
+\qff C_{\dff n\fff}(A)\fff\bigr)\fff.
\end{equation*}

\vspace*{-9pt}
Let\oss
$\displaystyle
i_{\dff *}\dff\colon H_{\dff n\fff}(A)\ttoo H_{\dff n\fff}(X)$\oss
and\oss
$\displaystyle
\partial_*\dff\colon H_{\fff n\qff +\qff 1\fff}(X\fff,\pff A)\ttoo H_{\dff n\fff}(A)$\oss
be the maps induced by the inclusion\oss
$\displaystyle
i\dff\colon A\ttoo X$\oss
and the boundary map of the homology sequence 
of the pair\qss $(X\fff,\pff A)$\qss respectively.\off\oss
The exactness of the homology sequence of the pair\qss $(X\fff,\pff A)$\qss 
implies\dss that the kernel of\dss $i_{\dff *}$\dss is equal to the image of\dss $\partial_*$\nnsp.\oss
Let\oss\vspace*{3pt} 
\[
\quad
K_{\dff n\fff}(A\to X)
\off =\off
\mbox{kernel of}\hspace*{1em}
 i_{\fff *}\dff\colon\dff H_{\dff n\fff}\dff(A)\ttoo H_{\dff n\fff}(X)\dff,
\]

\vspace*{-36pt}
\[
\quad
H_{\dff n\fff}(A\nnsp:\fff X) 
\off =\off 
H_{\dff n\fff}(A)\bigl/\partial_*\dff H_{\dff n\dff +\dff 1\fff}(X\fff,\pff A)
\off =\off 
H_{\dff n\fff}(A)\bigl/\dff K_{\dff n\fff}(A\to X)\dff.
\]

\vspace*{-9pt}
Obviously,\oss the canonical homomorphism\oss 
$H_{\dff n\fff}(A\nnsp:\fff X)\toto H_{\dff n\fff}(X)$\oss
is an isomorphism onto the image of\dss $i_{\fff *}$\nnsp.\oss
We will often identify\qss $H_{\dff n\fff}(A\nnsp:\fff X)$\qss with this image.

\myitpar{The framework.} 
Let us fix for the rest of the paper a closed connected surface\dss $S$\dss and 
a subsurface\dss $Q$\dss of\dss $S$\dnsp.\off\oss
Let\oss $c\off =\off \partial Q$\dnsp.\off\oss
The subsurface\dss $Q$\dss is not required to be connected,\oss
although the case of a connected~subsurface is the most important one.\off\oss
Let\vspace*{3pt} 
\[
\quad
\mathbold{K}\zero(c)
\off =\off
K\zero(c\to Q)\qff \cap\qff K\zero(c\to \ccomp Q)\dff,
\]

\vspace*{-36pt}
\[
\quad
\mathbold{Z}\one(Q\fff,\pff c)
\off =\off
\{\qff \alpha\dff\in\dff Z\one(Q\fff,\pff c) \pff\mid\off
\hclass{\partial\alpha}\qff \in\qff \mathbold{K}\zero(c)\pff\}\qff,
\]

\vspace{-36pt}
\[
\quad
\mathbold{H}\one(c)
\off =\off
H\one(c\nnsp:\fff S)
\off =\off
H\one(c)\bigl/\dff K\one(c\to S)\dff.
\]

\vspace{-9pt}
Since\oss
$\alpha\qff \in\qff Z\one(Q\fff,\pff c)$\oss
implies that\oss
$\hclass{\partial \alpha}\qff \in\qff K\zero(c\to Q)$\dnsp,\oss
the condition\oss
$\hclass{\partial\alpha}\qff \in\qff \mathbold{K}\zero(c)$\oss
in the definition of\qss
$\mathbold{Z}\one(Q\fff,\pff c)$\qss
can be replaced by\oss
$\hclass{\partial\alpha}\qff \in\qff K\zero(c\to \ccomp Q)$\dnsp.\oss

\myitpar{Weakly Torelli diffeomorphisms.}
Let\qss $\varphi\dff \colon\dff Q\toto Q$\qss
be a diffeomorphism fixed in a neighborhood of\dss $c$\nnsp.\off\oss
If\qss 
$\alpha\qff \in\qff Z\one(Q\fff,\pff c)$\dnsp,\oss 
then\oss 
\[
\quad
\partial\dff \varphi_{\fff *\fff}(\alpha)
\off =\off
\varphi_{\fff *\fff}(\partial\alpha)
\off =\off 
\partial \alpha
\]
and hence\qss
$\varphi_{\fff *\fff}(\alpha)\qff -\qff \alpha$\qss is a cycle in\dss $Q$\nnsp.\oss

The diffeomorphism\dss $\varphi$\dss 
is said to be\qss \emph{weakly Torelli}\oss if\dss for every cycle\dss
$\sigma$\dss in\dss $Q$\dss
the cycle\qss 
$\varphi_{\fff *\fff}(\sigma)\qff -\qff \sigma$\qss
is homologous to\dss $\old{0}$\dss in\dss $S$\dss
or\halfff,\oss what is the same,\oss 
$\dis
\hclass{\varphi_{\fff *\fff}(\sigma)\qff -\qff \sigma\dff}_{\dff S}
\off =\off 
\oold$\nnsp.\oss
If\dss $\varphi$\dss admits an extension to a Torelli diffeomorphism\qss 
$\psi\dff \colon\dff S\toto S$\nnsp,\oss then
\[
\quad
\hclass{\varphi_{\fff *\fff}(\sigma)\qff -\qff \sigma}_{\dff S}
\off =\off
\hclass{\psi_{\fff *\fff}(\sigma)\qff -\qff \sigma}_{\dff S}
\off =\off
\psi_{\fff *\fff}\hclass{\sigma}_{\dff S}\qff -\qff \hclass{\sigma}_{\dff S}
\off =\off
\oold
\]
for every\qss $\sigma\qff \in\qff Z\one(Q)$\qss 
and hence\dss $\varphi$\dss is a weakly Torelli diffeomorphism.\oss
It turns out that the property of being a weakly Torelli diffeomorphism has strong implications
for the cycles\qss
$\varphi_{\fff *\fff}(\alpha)\qff -\qff \alpha$\qss
for all relative cycles\qss 
$\alpha\qff \in\qff Z\one(Q\fff,\pff c)$\dnsp.\oss

\myitpar{The Mayer--Vietoris equivalence.}
The surface\qss $S$\qss is equal to the union of subsurfaces\qss $Q$\qss and\qss $\ccomp Q$\dnsp,\oss
with the intersection\oss 
$Q\dff\cap\dff \ccomp Q
\off =\off 
\partial Q
\off =\off 
\partial \ccomp Q
\off =\off 
c$\nnsp.\oss
By a well known result{\halfff},\oss usually hidden in the proof of the Mayer--Vietoris theorem,\oss 
the inclusion\vspace*{0pt}
\begin{equation*}
\quad
C_{\dff *\fff}(Q)\qff +\qff C_{\dff *\fff}(\ccomp Q)\off \toto\qff C_{\dff *\fff}(S)
\end{equation*}

\vspace*{-12pt}
is a homology equivalence,\oss 
i.e. it induces an isomorphism in homology\nsp.\oss

\myitpar{The Mayer--Vietoris boundary map.}
The Mayer--Vietoris boundary map\pss\vspace*{3pt}
\[
\quad
\mathcal{D}\dff \colon\dff H\one(S)\ttoo H\zero(c)
\] 

\vspace*{-9pt}
is defined as follows.\oss
By the Mayer--Vietoris equivalence,\oss
every homology class\qss $a\dff \in\dff H\one(S)$\qss had the form\qss
$a\qff =\qff \hclass{\alpha\qff +\qff \beta}$\qss
for some\oss $\alpha\dff\in\dff C\one(Q)$\oss and\oss
$\beta\dff\in\dff C\one(\ccomp Q)$\dnsp.\off\oss 
It is well known that\qss 
$\hclass{\fff \partial\alpha \fff}\qff \in\qff H\zero(c)$\qss
depends only on\qss $a$\nsp.\oss
By the definition,\oss
$\displaystyle
\mathcal{D}\dff(\halfff a\halfff)\off =\off \hclass{\fff \partial\alpha \fff}$\snsp.\oss 

The image of the Mayer--Vietoris boundary map\qss $\mathcal{D}$\qss
is equal to\qss $\mathbold{K}\zero(c)$\dnsp.\oss
Indeed,\oss the group\qss $\mathbold{K}\zero(c)$\qss 
is the group of homology classes of\dss $\oold$\dnsp-chains in\dss $c$\dss
which bound in both\dss $Q$\dss and\dss $\ccomp Q$\nnsp.\oss
Obviously\halfff,\oss the same is true for the image of\qss 
$\mathcal{D}$\nnsp.\oss

\myitpar{The intersection pairing.}
The orientation\dss $S$\dss
induces an orientation of\dss $Q$\snsp,\oss
which,\oss in turn,\oss induces an orientation of\dss $c\qff =\qff \partial Q$\nnsp.\oss
The Poincar\'{e} duality leads to a bilinear map\oss\vspace*{3pt}
\begin{equation}
\label{zero-one-intersection}
\quad
\langle\halfff \bullet\fff,\pff \bullet \fff\rangle_{\dff c}\dff\colon\dff
H\zero(c)\dff \times\dff H\one(c) \ttoo \zzz\dff
\end{equation}

\vspace*{-9pt}
called the\qss \emph{intersection pairing}\halfff.\oss
By the Poincar\'{e} duality or
by the following description this pairing is non-degenerate.\qff\oss
For each component\qss $C$\qss of\qss $c$\dnsp,\off\oss
let\qss $o\fff(C)\qff \in\qff H\zero(c)$\qss be the homology class of
any singular\qss $\old{0}$\dnsp-simplex in\qss $C$\dnsp,\oss
and let\qss $\hclass{C}$\qss be the fundamental class of\qss $C$\dnsp.\qff\oss 
Then\oss\vspace*{3pt}
\begin{equation}
\label{zero-one-duality}
\quad
\langle\qff o\fff(C)\fff,\pff \hclass{D} \qff\rangle_{\dff c}
\off =\off \old{1}
\quad\mbox{if}\quad
C\off =\off D\fff,
\hspace*{1.5em}\mbox{and}\hspace*{1.5em}
\langle\qff o\fff(C)\fff,\pff \hclass{D} \qff\rangle_{\dff c}
\off =\off \old{0}
\quad\mbox{if}\quad
C\off \neq\off D\dff.
\end{equation}

\vspace*{-9pt}
Let\oss
$\dis
\iota\dff \colon\dff 
H\one(c)\toto H\one(S)$\oss 
be the inclusion ho\-mo\-mor\-phism.\oss
The pairing\qss 
$\langle\halfff \bullet\fff,\pff \bullet \fff\rangle_{\dff c}$\qss 
is related
to the standard intersection pairing\vspace*{3pt}
\[
\quad
\langle\halfff \bullet\fff,\pff \bullet \fff\rangle\dff\colon\dff
H\one(S)\dff \times H\one(S) \ttoo \zzz\dff
\]

\vspace*{-9pt}
and the Mayer--Vietoris boundary map\qss 
$\mathcal{D}$\qss
by the formula\oss\vspace*{3pt}
\begin{equation}
\label{mv-pairing}
\quad
\langle\dff a\qff\fff,\pff\dff \iota\dff(b) \qff\rangle
\off\off =\off\off
\langle\pff \mathcal{D}\fff(a)\fff,\pff\fff b \qff\rangle_{\dff c}\dff,
\end{equation}

\vspace*{-9pt}
where\qss $a\dff\in\dff H\one(S)$\qss
and\qss
$b\dff\in\dff H\one(c)$\dnsp.\oss

\mypar{Lemma.}{bf-pairing} 
\emph{The groups\qss $\mathbold{K}\zero(c)$\qss and\qss
$K\one(c\to S)$\qss are the orthogonal complements of each other with
respect to the pairing\qss
$\langle\halfff \bullet\fff,\pff \bullet \fff\rangle_{\dff c}$\nnsp.\oss}

\proof\qss 
Let\qss $\theta\dff\in\dff H\zero(c)$\dnsp.\oss 
Then\qss $\theta$\qss 
is\dss a\dss linear combination of the form\vspace*{3pt} 
\[
\quad
\theta
\off =\off
\sum\nolimits_{\qff C}\qff n\fff(C)\dff o\fff(C)\dff,
\]

\vspace{-9pt}
where\dss $C$\dss runs over all components of\dss $c$\dss and\dss
$n\fff(C)\qff \in\qff \zzz$\nnsp.\oss 
If\dss $P$\dss is a component of\dss $Q$\nnsp,\oss then\vspace*{3pt} 
\[
\quad
\langle\pff \theta\fff,\pff\fff \hclass{\partial P} \qff\rangle_{\dff c}
\off =\off
\sum\nolimits_{\qff C}\qff n\fff(C)
\]

\vspace*{-9pt}
where\dss $C$\dss runs over all components of\dss $\partial P$\nnsp.\oss
On the other hand,\oss $\theta\qff \in\qff K\zero(c\to Q)$\qss
if and only of all such sums are equal to\dss $\oold$\nnsp.\oss
Since\qss $K\one(c\to Q)$\qss is generated by the fundamental classes\dss $\hclass{\partial P}$\dss
for components\dss $P$\dss of\dss $Q$\nnsp,\oss
it follows that\qss $K\zero(c\to Q)$\qss 
and\qss $K\one(c\to Q)$\qss
are the orthogonal complements of each other\halfff.\oss
A completely similar argument proves that\qss $K\zero(c\to \ccomp Q)$\qss 
and\qss $K\one(c\to \ccomp Q)$\qss
are the orthogonal complements of each other\halfff.\oss

By the Mayer--Vietoris equivalence the kernel\qss
$K\one(c\to S)$\qss is equal to the sum of the kernels\qss
$K\one(c\to Q)$\qss and\qss $K\one(c\to \ccomp Q)$\dnsp.\oss
By the definition\qss
$\mathbold{K}\zero(c)$\qss is equal to the intersection of the kernels\qss
$K\zero(c\to Q)$\qss and\qss $K\one(c\to \ccomp Q)$\dnsp.\oss
By combining these facts with the previous paragraph,\oss
we see that\qss $\mathbold{K}\zero(c)$\qss 
and\qss $K\one(c\to S)$\qss
are the orthogonal complements of each other\halfff.\oss  \eproof

\mypar{Lemma.}{duality-with-boundary}
\emph{A cycle\qss $\delta$\qss in\qss $Q$\qss
is\dss homologous\dss to\dss a\dss cycle\dss in\qss $c$\qss
if\dss and\dss only\dss if\qff\oss
$\langle\qff \delta\fff,\pff \tau \qff\rangle
\off =\off
\oold$\oss
for\dss all\dss cycles\qss $\tau$\qss in\qss $Q$\nnsp.\oss}

\proof\qss
The\qss \emph{``only\dss if''}\qss part is obvious.\oss
Let us prove the\qss \emph{``if''}\qss part\halfff.\oss
The exactness of the part
\[
\quad
H\one(c)\ttoo H\one(Q)\ttoo H\one(Q\fff,\pff c)
\]
of the homology sequence of the pair\qss $(Q\fff,\pff c)$\qss
implies that\qss $\delta$\qss is homologous to a cycle in\qss $c$\qss
if and only if the homology class of\qss $\delta$\qss in the relative
homology group\qss $H\one(Q\fff,\pff c)$\qss is equal to\dss $\oold$\nnsp.\oss
By the Poincar\'{e}--Lefschetz duality\qss
$a\qff \in\qff H\one(Q\fff,\pff c)$\qss 
is\dss equal\dss to\dss $\oold$\dss if\dss and\dss only\dss if\qss
$\langle\qff a\fff,\pff b \qff\rangle\qff =\qff \oold$\qss
for all\qss $b\qff \in\qff H\one(Q)$\dnsp.\oss
It remains to note that
\[
\quad
\langle\qff a\fff,\pff b \qff\rangle
\off =\off
\langle\qff \delta\fff,\pff \tau \qff\rangle
\]
if\dss $a$\dss is the homology class of\qss $\delta$\qss and\oss 
$b\off =\off \hclass{\tau}$\oss 
for some\qss $\tau\qff \in\qff Z\one(Q)$\dnsp.\oss  \eproof

\mypar{Theorem.}{weakly-torelli}
\emph{If\qss $\varphi$\qss is a weakly Torelli diffeomorphism,\oss
then for every\qss $\alpha\qff \in\qff \mathbold{Z}\one(Q\fff,\pff c)$\qss
the cycle\qss
$\varphi_{\dff *}(\alpha)\qff -\qff \alpha$\qss 
is homologous in\dss $Q$\dss to a cycle in\dss $c$\nnsp.\oss}

\proof\qss
Suppose
that\qss $\alpha\qff \in\qff \mathbold{Z}\one(Q\fff,\pff c)$\dnsp.\qff\oss
Let\qss $\tau\qff \in\qff Z\one(Q)$\qss be an arbitrary cycle.\qff\oss
Let
\[
\quad
\delta
\off =\off 
\varphi_{\fff *\fff}(\alpha)\qff -\qff \alpha
\hspace*{1.5em}\mbox{ and }\hspace*{1.5em}
\sigma
\off =\off  
\varphi_{\fff *\fff}(\tau)\qff -\qff \tau\dff.
\]
Since\qss $\alpha\qff \in\qff \mathbold{Z}\one(Q\fff,\pff c)$\dnsp,\oss
there exists\qss
$\beta\qff \in\qff Z\one(\ccomp Q\fff,\pff c)$\qss
such that\oss $\partial \alpha\off =\off \partial \beta$\oss
and\dss hence\qss
\[
\quad
\gamma
\off =\off 
\alpha\qff -\qff \beta
\]
is a cycle.\off\oss

\textbf{1.}\oss
As the first step of the proof\halfff,\oss let us prove that\vspace*{3pt}
\begin{equation}
\label{sum-is-zero}
\quad
\langle\qff \delta\fff,\pff \tau \qff\rangle
\qff +\qff
\langle\qff \gamma\fff,\pff \sigma \qff\rangle
\qff +\qff
\langle\qff \delta\fff,\pff \sigma \qff\rangle
\off =\off
\oold\dff.
\end{equation}

\vspace{-9pt}
Let\qss $\psi\dff \colon\dff S\toto S$\qss is the extension of\qss $\varphi$\qss by the identity.\qff\oss
Then\qss\vspace*{3pt}
\[
\quad
\psi_{\fff *\fff}(\gamma)\qff -\qff \gamma
\off =\off
\psi_{\fff *\fff}(\alpha\qff +\qff \beta)\qff -\qff (\alpha\qff +\qff \beta)
\off =\off \varphi_{\dff *\fff}(\alpha)\qff -\qff \alpha
\off =\off
\delta\dff.
\]

\vspace*{-9pt}
Since\qss $\psi$\qss is a diffeomorphism,\pss 
$\psi$\qss preserves the intersection pairing and hence\vspace*{3pt}
\[
\quad
\langle\qff \gamma\fff,\pff \tau \qff\rangle
\off =\off
\langle\qff \psi_{\fff *\fff}(\gamma)\fff,\pff \psi_{\dff *\fff}(\tau) \qff\rangle
\off 
=\off
\langle\qff \psi_{\fff *\fff}(\gamma)\fff,\pff \varphi_{\dff *\fff}(\tau) \qff\rangle\dff.
\]

\vspace*{-9pt}
Obviously,\oss 
$\psi_{\fff *\fff}(\gamma)
\off =\off 
\gamma\qff +\qff \delta$\oss
and\oss 
$\varphi_{\fff *\fff}(\tau)
\off =\off \tau\qff +\qff \sigma$\oss 
and\dss hence\vspace*{3pt}
\[
\quad
\langle\qff \psi_{\fff *\fff}(\gamma)\fff,\pff \varphi_{\dff *\fff}(\tau) \qff\rangle
\off =\off
\langle\qff \gamma\qff +\qff \delta\fff,\pff \tau\qff +\qff \sigma \qff\rangle
\]

\vspace*{-36pt}
\[
\quad
\phantom{\langle\qff \psi_{\fff *\fff}(\gamma)\fff,\pff \varphi_{\dff *\fff}(\tau) \qff\rangle
\off }
=\off
\langle\qff \gamma\fff,\pff \tau \qff\rangle
\qff +\qff
\langle\qff \delta\fff,\pff \tau \qff\rangle
\qff +\qff
\langle\qff \gamma\fff,\pff \sigma \qff\rangle
\qff +\qff
\langle\qff \delta\fff,\pff \sigma \qff\rangle\dff.
\]

\vspace{-9pt}
By combining the last two chains of equalities we see that\vspace*{3pt}
\[
\quad
\langle\qff \gamma\fff,\pff \tau \qff\rangle
\off =\off
\langle\qff \gamma\fff,\pff \tau \qff\rangle
\qff +\qff
\langle\qff \delta\fff,\pff \tau \qff\rangle
\qff +\qff
\langle\qff \gamma\fff,\pff \sigma \qff\rangle
\qff +\qff
\langle\qff \delta\fff,\pff \sigma \qff\rangle\dff.
\]

\vspace{-9pt}
Cancelling the terms\qss 
$\langle\qff \gamma\fff,\pff \tau \qff\rangle$\qss
completes the proof\dss of\qss (\ref{sum-is-zero}).\oss

\textbf{2.}\oss
The next step is to prove that\vspace*{3pt}
\begin{equation}
\label{delta-sigma}
\quad
\langle\qff \delta\fff,\pff \sigma \qff\rangle
\off =\off
\oold\dff.
\end{equation}

\vspace*{-9pt}
Since\dss $\delta$\dss is a cycle in\dss $Q$\dss and\dss $\rho$\dss
is a cycle in\dss $c$\nnsp,\oss
the intersection number\oss
$\langle\qff \delta\fff,\pff \rho \qff\rangle
\off =\off
\oold$\nnsp.\oss
On the other hand,\oss
$\langle\qff \delta\fff,\pff \sigma \qff\rangle
\off =\off
\langle\qff \delta\fff,\pff \rho \qff\rangle$\oss
because\dss $\sigma$\dss and\dss $\rho$\dss are homologous.\oss
This implies\qss (\ref{delta-sigma}).\oss

\textbf{3.}\oss
Now we will use the weak Torelli property in order to prove that\vspace*{3pt}
\begin{equation}
\label{gamma-sigma}
\quad
\langle\qff \gamma\fff,\pff \sigma \qff\rangle
\off =\off
\oold\dff.
\end{equation}

\vspace*{-9pt}
Since\dss $\varphi$\dss is a weak Torelli diffeomorphism,\oss
the cycle\dss $\sigma$\dss is homologous to\dss $\oold$\dss in\dss $S$\nnsp.\oss
By the Mayer--Vietoris equivalence this implies that\oss
$\sigma$\oss
is homologous in\dss $Q$\dss to a cycle\qss $\rho\qff \in\qff Z\one(c)$\qss
such that\dss $\rho$\dss is homologous to\dss $\oold$\dss in\dss $\ccomp Q$\nnsp.\qff\oss
Let\qss $r\qff \in\qff H\one(c)$\qss be the homology class of\dss $\rho$\dss in\dss $c$\nnsp.\oss 
Then\qss 
$r\qff \in\qff K\one(c\to \ccomp Q)$\qss 
and\dss hence\qss
$r\qff \in\qff K\one(c\to S)$\dnsp.\oss

The fact that\dss $\sigma$\dss and\dss $\rho$\dss are homologous also implies that\vspace*{3pt}
\[
\quad
\langle\qff \gamma\fff,\pff \sigma \qff\rangle
\off =\off
\langle\qff \gamma\fff,\pff \rho \qff\rangle
\off =\off
\langle\qff\fff \hclass{\gamma}\fff,\pff \hclass{\rho}_{\dff S} \qff\rangle\dff.
\]

\vspace*{-9pt}
Obviously,\oss the inclusion\dss
$\iota$\dss maps\oss 
$r\off =\off \hclass{\rho\dff}_{\dff c}$\oss
to\qss
$\hclass{\rho\dff}_{\dff S}$\qss and 
hence\qss (\ref{mv-pairing})\qss implies that\vspace*{3pt}
\[
\quad
\langle\qff\fff \hclass{\gamma}\fff,\pff \hclass{\rho}_{\dff S} \qff\rangle
\off =\off
\langle\qff\fff \mathcal{D}\dff \hclass{\gamma}\fff,\pff 
r \fff\qff\rangle_{\dff c}
\]

\vspace*{-9pt}
But\oss $\mathcal{D}\dff \hclass{\gamma}\off =\off \hclass{\partial \alpha}$\oss
by the definition of\qss $\mathcal{D}$\qss
and\dss hence\vspace*{3pt}
\[
\quad
\langle\qff\fff \mathcal{D}\dff \hclass{\gamma}\fff,\pff 
r \fff\qff\rangle_{\dff c}
\off =\off
\langle\qff\fff \hclass{\partial \alpha}\fff,\pff 
r \fff\qff\rangle_{\dff c}
\]

\vspace*{-9pt}
Since\oss 
$\hclass{\partial \alpha}\qff \in\qff \mathbold{K}\zero(c)$\oss
and\oss
$r\qff \in\qff K\one(c\to S)$\dnsp,\qff\oss
Lemma\qss \ref{bf-pairing}\qss implies that\vspace*{3pt}
\[
\quad
\langle\qff\fff \hclass{\partial \alpha}\fff,\pff 
r \fff\qff\rangle_{\dff c}
\off =\off
\oold\dff.
\]

\vspace*{-9pt}
The last four displayed equalities together immediately imply\qss (\ref{gamma-sigma}).

\textbf{4.}\oss
Finally,\oss the equalities\qss
(\ref{sum-is-zero}),\oss
(\ref{delta-sigma}),\oss
and\qss
(\ref{gamma-sigma})\qss
together imply that\oss
$\langle\qff \delta\fff,\pff \tau \qff\rangle
\off =\off
\oold$\nnsp.\oss
Since\dss $\tau$\dss is an arbitrary cycle in\dss $Q$\nnsp,\oss
Lemma\qss \ref{duality-with-boundary}\qss implies that\oss 
$\delta\off =\off \varphi_{\fff *\fff}(\alpha)\qff -\qff \alpha$\oss
is homologous to a cycle in\dss $c$\nnsp.\oss
The theorem follows.\oss  \eproof

\myitpar{The difference classes.}
Suppose that\dss $\varphi$\dss is a weakly Torelli diffeomorphism.\qff\oss
Let us identify the group\oss
$\mathbold{H}\one(c)
\off =\off
H\one(c)\bigl/\dff K\one(c\to S)$\oss
with the image of the map\oss 
$\dis
\iota\dff \colon\dff 
H\one(c)\toto H\one(S)$\dnsp.\oss

Let\qss
$\alpha\qff \in\qff \mathbold{Z}\one(Q\fff,\pff c)$\dnsp.\oss
By Theorem\qss \ref{weakly-torelli}\qss the cycle\qss
$\varphi_{\fff *\fff}(\alpha)\qff -\qff \alpha$\qss
is homologous in\dss $Q$\dss to a cycle\dss $\gamma$\dss in\dss $c$\nnsp.\oss
The identification of\qss $\mathbold{H}\one(c)$\qss
with the image of\dss $\iota$\dss turns the image of\qss
$\hclass{\gamma}
\qff \in\qff
H\one(c)$\oss 
in the group\qss
$\mathbold{H}\one(c)$\qss
into the homology class\qss $\hclass{\gamma}_{\dff S}$\qss
of the cycle\dss $\gamma$\dss in\dss $S$\nnsp.\oss

Since\dss $\gamma$\dss is homologous to\oss
$\varphi_{\fff *\fff}(\alpha)\qff -\qff \alpha$\nnsp,\oss
the homology class\qss $\hclass{\gamma}_{\dff S}$\qss
is equal to the homology class\qss 
$\hclass{\varphi_{\fff *\fff}(\alpha)\qff -\qff \alpha}_{\dff S}$\qss
of the cycle\qss $\varphi_{\fff *\fff}(\alpha)\qff -\qff \alpha$\qss in\dss $S$\nnsp.\oss
In particular\halfff,\oss
the image of\oss
$\hclass{\gamma}
\qff \in\qff
H\one(c)$\oss 
in the group\qss
$\mathbold{H}\one(c)$\qss 
depends only on\dss $\alpha$\nnsp,\oss
and if\dss $\alpha$\dss is a cycle in\dss $Q$\nnsp,\oss
then this image is equal to\dss $\oold$\dss by the definition
of weakly Torelli maps.\oss

The image of\oss
$\hclass{\gamma}
\qff \in\qff
H\one(c)$\oss 
in the group\qss
$\mathbold{H}\one(c)$\qss 
is called the\qss \emph{difference class}\qss of\dss $\alpha$\dss 
and is denoted by\qss $\Delta_{\dff \varphi}\fff(\alpha)$\dnsp.\oss 
Obviously,\oss the map\oss $\alpha\pff \longmapsto\pff \Delta_{\dff \varphi}\fff(\alpha)$\oss
is a homomorphism\vspace*{3pt}
\[
\quad
\Delta_{\dff \varphi}\qff \colon\qff
\mathbold{Z}\one(Q\fff,\pff c)\ttoo \mathbold{H}\one(c)\dff.
\]

\vspace*{-9pt}
It is called the\qss \emph{$\Delta$\dnsp-difference map}\qss of\dss $\varphi$\nnsp.\oss
The difference class\qss $\Delta_{\dff \varphi}\dff(\alpha)$\qss
depends only on\dss $\partial \alpha$\nnsp.\oss
Indeed,\oss
if\qss
$\alpha\fff,\pff \beta\qff \in\qff \mathbold{Z}\one(Q\fff,\pff c)$\qss
and\oss $\partial \alpha\off =\off \partial \beta$\nnsp,\oss
then\qss
$\alpha\qff -\qff \beta$\qss
is a cycle in\dss $Q$\dss and\dss hence\vspace*{3pt}
\[
\quad
\Delta_{\dff \varphi}\fff(\alpha)\qff -\qff \Delta_{\dff \varphi}\dff(\beta)
\off =\off
\Delta_{\dff \varphi}\fff(\alpha\qff -\qff \beta)
\off =\off
\oold\dff.
\]

\vspace*{-9pt}
If\qss
$\delta\dff \in\dff C\one(c)$\dnsp,\oss
then\pss
$\varphi_{\fff *\fff}(\delta)\pff =\pff \delta$\pss
and hence\pss
$\varphi_{\fff *\fff}(\delta)\qff -\qff \delta
\pff =\pff
\old{0}
\qff \in\qff
Z\one(c)$\pss
because the diffeomorphism\dss $\varphi$\dss is fixed on\dss $c$\nnsp.\oss
It follows that\oss 
$\Delta_{\dff \varphi}\fff(\delta)
\off =\off
\oold$\oss
for every\qss
$\delta\dff \in\dff C\one(c)$\dnsp.\oss
By combining this observation with the fact that\qss 
$\Delta_{\dff \varphi}\fff(\alpha)$\qss
depends only on\dss $\partial \alpha$\nnsp,\oss
we see that\halfff,\oss moreover\halfff,\oss
$\Delta_{\dff \varphi}\fff(\alpha)$\qss
depends only on the homology class\qss 
$\hclass{\partial \alpha}\qff \in\qff H\zero(c)$\nnsp.\oss

\myitpar{The\dss $\delta$\dnsp-difference map.}
Let\oss
$b\qff \colon\qff
Z\one(Q\fff,\pff c)\toto H\zero(c)$\oss be the homomorphism\oss
$\alpha\qff \longmapsto\qff \hclass{\partial \alpha}$\nnsp.\oss
By the definition the group\qss
$\mathbold{Z}\one(Q\fff,\pff c)$\qss
is the preimage of\qss
$\mathbold{K}\zero(c)$\qss
under the map\dss $b$\nnsp.\oss
It follows that\dss $b$\dss induces a surjective homomorphism\oss
$\dis
\mathbold{b}\qff \colon\qff
\mathbold{Z}\one(Q\fff,\pff c)\ttoo \mathbold{K}\zero(c)$\dnsp.\qff\oss

If\dss $\varphi$\dss is a weakly Torelli diffeomorphism,\oss
then the difference classes\qss $\Delta_{\dff \varphi}\dff(\alpha)$\qss
depend only on\qss 
$\hclass{\partial \alpha}\qff \in\qff H\zero(c)$\qss
and hence the surjectivity of\dss $\mathbold{b}$\dss 
implies that
there exists a unique homomorphism\oss 
$\delta_{\fff \varphi}$\oss
such that the triangle\vspace*{6pt}
\[
\quad
\begin{tikzcd}[column sep=tiny, row sep=huge]
\mathbold{Z}\one(Q\fff,\pff c)\off 
\arrow[rr, "{\displaystyle \mathbold{b}}"] 
\arrow[dr, "{\displaystyle \mathbold{\Delta}_{\dff \varphi}}"' near start]  
& & \off\mathbold{K}\zero(c) 
\arrow[dl, dashed, "{\displaystyle \delta_{\fff \varphi}}" near start] \\
& \mathbold{H}\one(c)  & 
\end{tikzcd}
\]

\vspace*{-6pt}
is commutative.\oss
The map\dss $\delta_{\fff \varphi}$\dss is called\qss
\emph{$\delta$\dnsp-difference map}\pss of\dss $\varphi$\nnsp.\oss

\mypar{Lemma.}{weakly-torelli-composition} \emph{If\pss $\varphi\fff,\pff \psi$\pss 
are weakly Torelli diffeomorphisms,\oss
then\qss $\psi\circ\fff \varphi$\qss 
is also a weakly Torelli diffeomorphism,\qff\oss and\qff\oss
$\dis
\delta_{\trf \psi\dff\circ\dff \varphi}
\off =\off 
\delta_{\trf \psi}\qff +\qff \delta_{\fff \varphi}$\nsp.}

\proof\qss Let\qss $\alpha\qff \in\qff \mathbold{Z}\one(Q\fff,\pff c)$\dnsp.\oss
Then\qss
$\varphi_{\fff *\fff}(\alpha)\qff \in\qff \mathbold{Z}\one(Q\fff,\pff c)$\qss
because\dss $\varphi$\dss is fixed on\dss $c$\nnsp.\oss 
Obviously,\vspace*{1.5pt}
\begin{equation}
\label{sum-composition}
\quad
\psi_{\fff *\fff}\bigl( \varphi_{\fff *\fff}(\alpha)\fff\bigr)\qff -\qff \alpha
\off =\off
\Bigl(\dff 
\psi_{\fff *\fff}\bigl( 
\varphi_{\fff *\fff}(\alpha)\fff\bigr)\qff -\qff \varphi_{\fff *\fff}\dff(\alpha) 
\dff\Bigr)
\pff +\pff
\Bigl(\dff 
\varphi_{\fff *\fff}(\alpha)\qff -\qff \alpha 
\dff\Bigr)\dff.
\end{equation}

\vspace{-10.5pt}
If\halfff,\oss moreover\halfff,\oss $\alpha\qff \in\qff Z\one(Q)$\nnsp,\oss
the\dss $\varphi_{\fff *\fff}(\alpha)\qff \in\qff Z\one(Q)$\dss also 
and hence both\qss
$\varphi_{\fff *\fff}(\alpha)\qff -\qff \alpha$\qss
and\qss
$\psi_{\fff *\fff}\bigl( \varphi_*\dff(\alpha)\fff\bigr)\qff -\qff \varphi_{\fff *}\dff(\alpha)$\qss
are homologous to\dss $\oold$\dss in\dss $S$\nnsp.\oss
Together with\qss (\ref{sum-composition})\qss this implies that the cycle\qss
$\psi_{\fff *\fff}\bigl( \varphi_{\fff *\fff}(\alpha)\fff\bigr)\qff -\qff \alpha$\qss
is homologous to\dss $\oold$\dss in\dss $S$\nnsp.\oss
It follows that\qss $\psi\circ\fff \varphi$\qss is weakly Torelli.

In general,\oss
Theorem\qss \ref{weakly-torelli}\qss implies that the cycles\qss
$\varphi_{\fff *}(\alpha)\qff -\qff \alpha$\qss and\qss
$\psi_{\fff *}\bigl( \varphi_*\dff(\alpha)\fff\bigr)\qff -\qff \varphi_{\fff *}\dff(\alpha)$\qss
are homologous in\dss $Q$\dss to some cycles\qss
$\delta\fff,\pff \varepsilon\qff \in\qff C\one(c)$\qss respectively.\oss
In view of\qss (\ref{sum-composition})\qss this implies that\qss
$\psi_{\fff *\fff}\bigl( \varphi_{\fff *\fff}(\alpha) \fff\bigr)\qff -\qff \alpha$\qss
is homologous in\dss $Q$\dss to\qss 
$\varepsilon\qff +\qff \delta\qff \in\qff Z\one(c)$\dnsp.\oss

Let\oss 
$\theta
\off =\off 
\hclass{\partial \alpha}
\off =\off 
\hclass{\partial\dff \varphi_{\fff *\fff}(\alpha)}$\dnsp.\oss
The images in\qss $\mathbold{H}\one(c)$\qss of\qss 
$\hclass{\delta}$\nnsp,\qss 
$\hclass{\varepsilon}$\qss 
and\qss 
$\hclass{\varepsilon\qff +\qff \delta}$\qss 
are equal to\vspace*{1.5pt}
\[
\quad
\delta_{\fff \varphi}(\theta),
\hspace*{1.2em}
\delta_{\trf \psi\halfff}(\theta),
\hspace*{1.2em}\mbox{ and }\hspace*{1.2em}
\delta_{\trf \psi\dff\circ\dff \varphi}(\theta)
\]

\vspace*{-10.5pt}
respectively.\oss
Since the last one is equal to the sum of the first two,\oss\vspace*{1.5pt}
\[
\quad
\delta_{\trf \psi\dff\circ\dff \varphi}(\theta)
\off =\off 
\delta_{\trf \psi\halfff}(\theta)\qff +\qff \delta_{\fff \varphi}(\theta)\dff.
\] 

\vspace*{-10.5pt}
Since every chain\qss $\theta\qff \in\qff \mathbold{K}\zero(c)$\qss is equal to\qss
$\hclass{\partial \alpha}$\qss for some\qss 
$\alpha\qff \in\qff \mathbold{Z}\one(Q\fff,\pff c)$\dnsp,\oss
the last displayed equality implies that\oss
$\dis\delta_{\trf \psi\dff\circ\dff \varphi}
\off =\off 
\delta_{\trf \psi}\qff +\qff \delta_{\fff \varphi}$\nsp.\oss  \eproof

\mypar{Lemma.}{weakly-torelli-inverse}
\emph{If\dss $\varphi$\dss is a weakly Torelli diffeomorphism,\oss 
then\qss $\varphi^{\dff -\dff \oldstylenums{1}}$\qss
is also a weakly Torelli diffeomorphism and\oss
$\dis
\delta_{\dff \varphi^{\dff -\dff \oldstylenums{1}}}
\off =\off 
-\qff \delta_{\dff \varphi}$\nsp.\oss}

\proof\qss
Let\qss $\alpha\qff \in\qff \mathbold{Z}\one(Q\fff,\pff c)$\dnsp.\oss
Then\oss
$\dis 
\beta
\off =\off 
\varphi^{\dff -\dff \oldstylenums{1}\dff}_{\fff *\fff}(\alpha)
\qff \in\qff \mathbold{Z}\one(Q\fff,\pff c)$\qss and\oss\vspace*{1.5pt} 
\begin{equation}
\label{minus-cycle}
\quad
\varphi^{\dff -\dff \oldstylenums{1}\dff}_{\fff *\fff}(\alpha) -\qff \alpha
\off =\off
\beta\qff -\qff \varphi_{\fff *\fff}(\beta)
\end{equation}

\vspace*{-10.5pt}
If\halfff,\oss moreover\halfff,\oss $\alpha\qff \in\qff Z\one(Q)$\nnsp,\oss
then\qss $\beta\qff \in\qff Z\one(Q)$\qss also 
and hence\qss
$\varphi_{\fff *\fff}(\beta)\qff -\qff \beta$\qss
is homologous to\dss $\oold$\dss in\dss $S$\nnsp.\oss
In view of\qss (\ref{minus-cycle})\qss this implies that the cycle\qss
$\varphi^{\dff -\dff \oldstylenums{1}\dff}_{\fff *\fff}(\alpha) -\qff \alpha$\qss
is homologous to\dss $\oold$\dss in\dss $S$\nnsp.\oss
It follows that\dss $\varphi^{\dff -\dff \oldstylenums{1}}$\dss is a weakly Torelli diffeomorphism.

In general,\oss
Theorem\qss \ref{weakly-torelli}\qss implies that the cycle\qss
$\varphi_{\fff *\fff}(\beta)\qff -\qff \beta$\qss
is homologous in\dss $Q$\dss to a cycle\qss $\delta\qff \in\qff Z\one(c)$\dnsp.\oss
In view\qss of\qss (\ref{minus-cycle})\qss 
$\varphi^{\dff -\dff \oldstylenums{1}\dff}_{\fff *\fff}(\alpha) -\qff \alpha$\qss 
is homologous in\dss $Q$\dss to\qss $-\qff \delta$\nnsp.\oss

Let\oss 
$\theta
\off =\off 
\hclass{\partial \alpha}
\off =\off 
\hclass{\partial \beta}$\dnsp.\oss
Then\qss
$\delta_{\dff \varphi}(\theta)$\qss
is equal to the image of\qss $\hclass{\delta}$\qss in\qss $\mathbold{H}\one(c)$\qss
and\qss
$\delta_{\dff \varphi^{\dff -\dff \oldstylenums{1}}}(\theta)$\qss
is equal to the image of\qss $\hclass{-\qff \delta}$\qss in\qss $\mathbold{H}\one(c)$\dnsp.\oss
Therefore\oss
$\delta_{\dff \varphi}(\theta)
\off =\off 
-\qff \delta_{\dff \varphi^{\dff -\dff \oldstylenums{1}}}(\theta)$\dnsp.\oss
Since every chain\qss $\theta\qff \in\qff \mathbold{K}\zero(c)$\qss is equal to\qss
$\hclass{\partial \alpha}$\qss for some\qss 
$\alpha\qff \in\qff \mathbold{Z}\one(Q\fff,\pff c)$\dnsp,\oss
the lemma follows.\oss  \eproof

\mysection{Extensions\qss by\qss the\qss identity\qss 
and\qss the\qss symmetry\qss property}{extensions-identity}

\vspace*{6pt}
\mypar{Theorem.}{extension-of-weakly-torelli} \emph{Let\pss 
$\varphi\dff \colon\dff Q\toto Q$\pss 
be a weakly Torelli diffeomorphism and let\pss
$\psi\dff \colon\dff S\toto S$\oss 
be the extension of\pss $\varphi$\pss by the identity map of\oss $\ccomp Q$\nnsp.\qff\oss 
Then}\oss\vspace*{4pt}
\[
\quad
\psi_{\fff *\fff}(\fff a\fff)\qff -\qff a
\off =\off
\delta_{\fff \varphi}\dff \circ\pff \mathcal{D}\dff(\halfff a\halfff)
\] 

\vspace*{-9pt}
\emph{for all\oss
$a\qff \in\qff H\one(S)$\dnsp,\oss
where the target\oss $\mathbold{H}\one(c)$\qss
of\oss $\delta_{\fff \varphi}$\pss
is identified with the image of the inclusion homomorphism\qss
$\iota\dff \colon\dff H\one(c)\toto H\one(S)$\dnsp.\oss}

\proof\oss Let\qss $a\qff \in\qff H\one(S)$\snsp.\oss
By the Mayer--Vietoris equivalence,\oss
$a\off =\off \hclass{\alpha\qff -\qff \beta}$\oss
for some chains\oss $\alpha\qff \in\qff C\one(Q)$\oss
and\oss 
$\beta\qff \in\qff C\one(\ccomp Q)$\oss
such that\qss
$\alpha\qff -\qff \beta$\qss is a cycle.\oss
Then\oss
$\dis
\mathcal{D}\dff(\halfff a\halfff)\off =\off \hclass{\fff \partial\alpha \fff}$\nnsp.\oss

Since\qss $\alpha\qff -\qff \beta$\qss is a cycle,\oss
$\displaystyle
\partial\alpha
\off =\off 
\partial \beta
\qff \in\qff C\zero(c)$\oss
and hence\oss
$\dis
\partial\alpha
\off =\off 
\partial \beta$\oss
bounds both in\dss $Q$\dss and in\dss $\ccomp Q$\nnsp.\qff\oss
It\dss follows that\pss 
$\hclass{\partial \alpha}
\qff \in\qff \mathbold{K}\zero(c)$\pss
and\pss
$\alpha\dff\in\dff \mathbold{Z}_{\dff 1\fff}(Q\fff,\pff c)$\nnsp.\oss 
Therefore\qss
$\Delta_{\dff \varphi}\dff(\alpha)$\qss is defined.\oss

After the identification of\qss $\mathbold{H}\one(c)$\qss
with the image of\qss $\iota\dff \colon\dff H\one(c)\toto H\one(S)$\qss 
the difference class\oss 
$\Delta_{\dff \varphi}\dff(\alpha)$\qss
turns into the homology class\oss
$\hclass{\varphi_{\fff *\fff}(\alpha)\qff -\qff \alpha}_{\dff S}\qff \in\qff H\one(S)$\oss
of\qss $\varphi_{\fff *\fff}(\alpha)\qff -\qff \alpha$\qss in\dss $S$\nnsp.\oss

On the other hand,\oss
$\psi_{\fff *\fff}(\fff a\fff)\qff -\qff a$\qss is equal to 
the homology class in\dss $S$\dss of the cycle\vspace*{3pt}
\[
\quad
\psi_{\fff *\fff}(\alpha\qff -\qff \beta)\qff -\qff (\alpha\qff -\qff \beta)
\off =\off
\psi_{\fff *\fff}(\alpha)\qff -\qff \psi_{\fff *\fff}(\beta)
\qff -\qff 
\alpha\qff +\qff \beta 
\]

\vspace*{-36pt}
\[
\quad
\phantom{\psi_{\fff *\fff}(\alpha\qff -\qff \beta)\qff -\qff (\alpha\qff -\qff \beta)
\off }
=\off
\varphi_{\fff *\fff}(\alpha)\qff -\qff \beta\qff -\qff \alpha\qff +\qff \beta
\off =\off
\varphi_{\fff *\fff}(\alpha)\qff -\qff \alpha\dff.
\]

\vspace*{-9pt}
It follows that\oss
$\dis
\psi_{\fff *\fff}(\fff a\fff)\qff -\qff a
\off =\off
\Delta_{\dff \varphi}\fff(\alpha)$\dnsp.\qff\oss
But\oss 
$\Delta_{\dff \varphi}\fff(\alpha)
\off =\off
\delta_{\fff \varphi\fff}(\dff \hclass{\fff \partial \alpha \fff} \dff)$\oss
by the definition of\qss $\delta_{\fff \varphi\fff}$\nnsp.\oss
By combining the last two equalities with\oss
$\dis
\mathcal{D}\dff(\halfff a\halfff)\off =\off \hclass{\fff \partial\alpha \fff}$\oss
we see that\vspace*{3pt}
\[
\quad
\psi_{\fff *\fff}(\fff a\fff)\qff -\qff a
\off =\off
\Delta_{\dff \varphi}\fff(\alpha)
\off =\off
\delta_{\fff \varphi\fff}(\dff \hclass{\fff \partial \alpha \fff} \dff)
\off =\off
\delta_{\fff \varphi}\dff \circ\pff \mathcal{D}\dff(\halfff a\halfff)\dff.
\]

\vspace*{-9pt}
The theorem follows.\oss  \eproof

\mypar{Theorem.}{extension-by-identity} 
\emph{Let\pss 
$\varphi\dff \colon\dff Q\toto Q$\pss 
be a weakly Torelli diffeomorphism and let\pss
$\psi\dff \colon\dff S\toto S$\oss 
be the extension of\pss $\varphi$\pss by the identity map of\oss $\ccomp Q$\nnsp.\qff\oss 
Then\qss $\psi$\qss is a Torelli diffeomorphism 
if and only if\oss
$\delta_{\fff \varphi}\off =\off \old{0}$\nnsp.\oss}\vspace*{3pt}

\proof\qss 
If\oss 
$\delta_{\fff \varphi}\off =\off \old{0}$\snsp,\oss
then Theorem\qss \ref{extension-of-weakly-torelli}\qss implies that the induced map\qss 
$\psi_{\fff *\fff}\dff \colon\dff H\one(S)\toto H\one(S)$\qss is equal to the identity
and hence\qss $\psi$\qss is a Torelli diffeomorphism.\oss
Conversely,\qff\oss
if\qss
$\psi$\qss is a Torelli diffeomorphism,\oss
then\qss $\varphi$\qss is a weakly Torelli diffeomorphism
by the remarks after the definition of the latter\halfff.\oss  \eproof

\myitpar{The induced pairing.} 
Lemma\qss \ref{bf-pairing}\qss 
implies that\pss (\ref{zero-one-intersection})\pss induces a non-degenerate pairing\vspace*{2pt} 
\[
\quad
\langle\halfff \bullet\fff,\pff \bullet \fff\rangle_{\dff c}\qff\colon\qff
\mathbold{K}\zero(c)\dff \times \mathbold{H}\one(c)
\qff \ttoo\qff \zzz\dff.
\]

\vspace*{-10pt}
The identification of\qss $\mathbold{H}\one(c)$\qss
with the image of\oss 
$\dis
\iota\dff \colon\dff 
H\one(c)\toto H\one(S)$\oss
turns\qss (\ref{mv-pairing})\qss into\oss\vspace*{2pt}
\begin{equation}
\label{mv-induced-pairing}
\quad
\langle\dff a\halfff\fff,\pff\dff b\dff\rangle
\off\off =\off\off
\langle\qff \mathcal{D}\fff(a)\fff,\pff\fff b\dff\rangle_{\dff c}\dff,
\end{equation}

\vspace*{-10pt}
where\qss $b\dff\in\dff \mathbold{H}\one(c)$\dnsp,\pss
$a\dff\in\dff H\one(S)$\qss
and\dss hence\qss
$\mathcal{D}\dff(a)\qff \in\qff \mathbold{K}\zero(c)$\qss
by the remarks after the definition of
the Mayer--Vietoris boundary map\qss $\mathcal{D}$\nnsp.\oss
A\dss homomorphism\vspace*{2pt}
\[
\quad
\delta\dff \colon\dff
\mathbold{K}\zero(c)\ttoo \mathbold{H}\one(c)
\]

\vspace*{-10pt}
is said to be\pss \emph{symmetric}\oss if\oss 
$\dis
\langle\qff a\halfff\fff,\off \delta\fff(b) \qff\rangle_{\dff c}
\off =\off
\langle\qff b\halfff,\off \delta\fff(a) \qff\rangle_{\dff c}
$\oss 
for\dss all\oss $a\fff,\pff b\qff\in\qff  \mathbold{K}\zero(c)$\dnsp.\oss

\mypar{Theorem.}{symmetry-theorem} \emph{If\pss $\varphi$\qss
is a weakly Torelli diffeomorphism\halfff,\oss 
then\pss $\delta_{\fff \varphi}$\pss 
is symmetric.\oss}

\proof\qss 
Let\pss
$\psi\dff \colon\dff S\toto S$\oss 
be the extension of\pss $\varphi$\pss by the identity.\qff\oss 
By Theorem\qss \ref{extension-of-weakly-torelli}\qss\vspace*{2pt} 
\[
\quad
\psi_{\fff *\fff}(\fff a\fff)\qff -\qff a
\off =\off
\delta_{\fff \varphi}\dff \circ\pff \mathcal{D}\dff(\halfff a\halfff)
\]

\vspace*{-10pt}
for all\oss $a\qff \in\qff H\one(S)$\dnsp.\oss
Let\oss $a\fff,\pff b\qff \in\qff  H\one(S)$\dnsp.\oss
Since\qss $\psi_{\fff *\fff}$\qss preserves the intersection pairing\halfff,\oss\vspace*{2pt}
\[
\quad
\langle\qff a\fff,\pff b \qff\rangle
\off =\off
\langle\qff a\qff +\qff \delta_{\fff \varphi}\dff \circ\pff \mathcal{D}\dff(\halfff a\halfff)\fff,\off\qff 
b\qff +\qff \delta_{\fff \varphi}\dff \circ\pff \mathcal{D}\dff(\halfff b\halfff) \qff\rangle\dff.
\] 

\vspace*{-10pt}
The intersection pairing\qss
$\langle\halfff \bullet\fff,\pff \bullet \fff\rangle$\qss 
is equal to\dss $\oold$\dss on the image of\oss 
$H\one(c)\toto H\one(S)$\dnsp,\oss and hence\vspace*{2pt}
\[
\quad
\langle\qff \delta_{\fff \varphi}\dff \circ\pff \mathcal{D}\dff(\halfff a\halfff)\fff,\qff\off 
\delta_{\fff \varphi}\dff \circ\pff \mathcal{D}\dff(\halfff b\halfff) \qff\rangle
\off =\off 
\oold\dff.
\]

\vspace*{-10pt}
It follows that\oss 
$\langle\qff a\fff,\pff b \qff\rangle
\off\off =\off\off
\langle\qff a\fff,\pff b \qff\rangle
\qff +\qff
\langle\qff a\fff,\off \delta_{\fff \varphi}\dff \circ\pff \mathcal{D}\dff(\halfff b\halfff) \qff\rangle
\qff +\qff
\langle\qff \delta_{\fff \varphi}\dff \circ\pff \mathcal{D}\dff(\halfff a\halfff)\fff,\off b \qff\rangle
$\dnsp,\qff\oss 
and\dss hence\vspace*{2pt} 
\[
\quad
\langle\qff a\fff,\off \delta_{\fff \varphi}\dff \circ\pff \mathcal{D}\dff(\halfff b\halfff) \qff\rangle
\qff +\qff
\langle\qff \delta_{\fff \varphi}\dff \circ\pff \mathcal{D}\dff(\halfff a\halfff)\fff,\off b \qff\rangle
\off =\off
\oold\dff.
\]

\vspace*{-10pt}
Since the intersection pairing is skew-symmetric,\oss this,\oss in\dss turn,\oss implies that\vspace*{2pt}
\[
\quad
\langle\qff a\fff,\off \delta_{\fff \varphi}\dff \circ\pff \mathcal{D}\dff(\halfff b\halfff) \qff\rangle
\off\off =\off\off
\langle\qff b\fff,\off \delta_{\fff \varphi}\dff \circ\pff \mathcal{D}\dff(\halfff a\halfff) \qff\rangle\dff.
\]

\vspace*{-10pt}
By combining this identity with\qss (\ref{mv-induced-pairing}),\oss
we get\oss\vspace*{2pt}
\[
\quad
\langle\qff \mathcal{D}\fff(a)\fff,\off 
\delta_{\fff \varphi}\dff \circ\pff \mathcal{D}\dff(\halfff b\halfff) \qff\rangle_{\dff c}
\off\off =\off\off
\langle\qff \mathcal{D}\fff(b)\fff,\off 
\delta_{\fff \varphi}\dff \circ\pff \mathcal{D}\dff(\halfff a\halfff) \qff\rangle_{\dff c}\dff.
\]

\vspace*{-10pt}
By the remarks after the definition of\qss $\mathcal{D}$\qss the image 
of\qss $\mathcal{D}$\qss is equal to\qss
$\mathbold{K}_{\dff 0\fff}(c)$\dnsp.\qff\oss
It follows that\qss 
$\delta_{\fff \varphi}$\qss is symmetric.\oss   \eproof

\mysection{General\qss extensions}{general-extensions}

\vspace*{\medskipamount}
\myitpar{The connectivity assumption.}
For the rest of the paper we will assume that the subsurface\dss $Q$\dss is connected.\oss
Of course,\oss this is the most important case,\oss
and the general case,\oss at least in principle,\oss
can be approached by dealing with each component separately.\oss
But one should keep in mind that in Torelli topology 
glueing the information about different components is not automatic.\oss
The glueing problem is one of the topics of the sequel\qss \cite{i-et-iii}\qss
of this paper\halfff.\oss

The main features of the case of connected\dss $Q$\dss are the following\halfff.\oss
The kernel\qss ${K\zero(c\to \ccomp Q)}$\qss is contained in the
kernel\qss $K\zero(c\to Q)$\dnsp.\oss
The kernel\qss ${K\one(c\to Q)}$\qss 
is generated by the fundamental class\dss $\hclass{\partial Q}$\dss
and hence is contained in\qss $K\one(c\to \ccomp Q)$\dnsp.\oss
It follows that\vspace*{3pt}
\[
\quad
\mathbold{K}\zero(c)\qff =\qff K\zero(c\to \ccomp Q)
\hspace*{1.5em}\mbox{ and }\hspace*{1.5em}
K\one(c\to S)\qss =\qff K\one(c\to \ccomp Q)\dff.
\]

\vspace*{-9pt}
Hence\oss 
$\mathbold{H}\one(c)
\off =\off
H\one(c)\bigl/\dff K\one(c\to \ccomp Q)$\dnsp.\oss

\myitpar{Completely reducible homomorphisms.} 
There is a canonical direct sum decomposition\vspace*{3pt}
\begin{equation*}
\quad
H_{\dff *\fff}(c)\off 
=\off \bigoplus\nolimits_{\trf P} \dff H_{\dff *\fff}(\partial P),
\end{equation*}

\vspace*{-9pt}
where\qss $P$\qss runs over all components of\qss $\ccomp Q$\dnsp.\off\oss
Since the subsurface\dss $Q$\dss is connected and hence\qss
$\mathbold{K}\zero(c)\qss =\qff K\zero(c\to \ccomp Q)$\dnsp,\oss
this decomposition leads to the direct sum decomposition\vspace*{3pt}
\begin{equation*}
\quad
\mathbold{K}\zero(c)\off 
=\off \bigoplus\nolimits_{\trf P} \dff \mathbold{K}\zero(\partial P)\dff,
\end{equation*}

\vspace*{-9pt}
where\oss 
$\displaystyle
\mathbold{K}\zero(\partial P)
\off =\off 
K\zero(\partial P\to P)$\dnsp,\oss
and to the direct sum decomposition\vspace*{3pt}
\begin{equation*}
\quad
\mathbold{H}\one(c)
\off =\off 
\bigoplus\nolimits_{\trf P} \dff \mathbold{H}\one(\partial P)\dff,
\end{equation*} 

\vspace*{-9pt}
where\oss
$\displaystyle
\mathbold{H}\one(\partial P)
\off =\off
H\one(\partial P\nsp:\fff P\halfff)
\off =\off
H\one(\partial P)\bigl/\dff K\one(\partial P\to P)$\dnsp.\off\oss
A homomorphism\oss\vspace*{3pt}
\[
\quad
\delta\dff\colon\dff \mathbold{K}\zero(c)\ttoo \mathbold{H}\one(c)
\]

\vspace*{-9pt}
is called\qss \emph{completely reducible}\qss if\qss $\delta$\qss
respects the above direct sum decompositions,\qff\oss 
i.e.\qss if\pss $\delta$\qss
maps\qss
$\mathbold{K}\zero(\partial P)$\qss into\qss
$\mathbold{H}\one(\partial P)$\qss
for every component\qss $P$\qss of\qss $\ccomp Q$\dnsp.\off\oss
If\qss $\delta$\qss is completely reducible,\oss
then\qss $\delta$\qss admits a direct sum decomposition\vspace*{3pt}
\begin{equation*}
\label{delta-decomposition}
\quad
\delta\off =\off  \bigoplus\nolimits_{\trf P} \dff \delta_{\fff P}
\end{equation*}

\vspace*{-9pt} 
into induced homomorphisms\oss
$\displaystyle
\delta_{\fff P}\dff\colon\dff 
\mathbold{K}\zero(\partial P)\ttoo 
\mathbold{H}\one(\partial P)$\dnsp.\oss

\mypar{Theorem.}{torelli-is-completely-reducible} 
\emph{Of\pss $\varphi\dff \colon\dff Q\toto Q$\pss can be
extended to a Torelli diffeomorphism\pss $S\toto S$\dnsp,\oss
then\pss $\varphi$\pss is a weakly Torelli diffeomorphism and\pss 
$\delta_{\fff \varphi}$\pss 
is completely reducible.}

\proof\qss 
By the remarks after the definition of weakly Torelli diffeomorphisms\qss $\varphi$\qss 
is a weakly Torelli diffeomorphism.\qff\oss
Let\dss $R$\dss be an arbitrary component of\dss $\ccomp Q$\dnsp,\qff\oss
and let\qss $\alpha\dff\in\dss \mathbold{Z}\one(Q\fff,\pff c)$\dnsp.\oss
It is sufficient to prove that\oss
$\displaystyle
\hclass{\partial\alpha}
\qff \in\qff 
\mathbold{K}\zero(\halfff\partial R\fff)$\oss
implies\oss
$\displaystyle
\Delta_{\dff \varphi}\fff(\alpha)
\qff \in\qff 
\mathbold{H}\one(\halfff\partial R\fff)$\dnsp.\oss

Let us present\dss $\partial \alpha$\dss as the sum\qss $\theta\qff +\qff \rho$\nnsp,\oss
where\qss $\theta\qff \in\qff C\zero(\partial R\fff)$\qss and\qss
$\rho\qff \in\qff C\zero(c\dff \smallsetminus\dff \partial R\fff)$\dnsp.\oss
Then\qss $\hclass{\partial \alpha}\qff =\qff \hclass{\theta}$\qss
and\qss $\hclass{\rho}\qff =\qff \oold$\nnsp.\oss
It follows that\dss $\theta$\dss bounds in\dss $R$\nnsp.\oss
Since\dss $Q$\dss is connected,\oss $\theta$\dss bounds in\dss $Q$\dss also,\oss
i.e.\oss $\theta\qff =\qff \partial \alpha'$\qss for some\qss
$\alpha'\qff \in\qff Z\one(Q\fff,\pff c)$\dnsp.\oss
Then\qss $\hclass{\partial \alpha'\fff}\qff =\qff \hclass{\partial \alpha}$\qss
and\qss $\partial \alpha'\qff \in\qff C\zero(\partial R\fff)$\dnsp.\oss
On the other hand,\oss 
$\Delta_{\dff \varphi}\fff(\alpha)$\qss 
depends only on\qss $\hclass{\alpha}$\qss
and hence\oss 
$\Delta_{\dff \varphi}\fff(\alpha'\fff)
\off =\off
\Delta_{\dff \varphi}\fff(\alpha)$\dnsp.\oss

Therefore we may assume that\qss $\partial \alpha\qff \in\qff C\zero(\partial R\fff)$\dnsp.\oss
Then the assumption\oss
$\displaystyle
\hclass{\partial\alpha}
\qff \in\qff 
\mathbold{K}\zero(\halfff\partial R\fff)$\oss
implies that there exists\oss
$\beta\dff\in\dff C\one(\halfff R\fff)$\oss
such that\oss
$\partial\alpha\off =\off \partial\beta$\oss
and\dss hence\qss\vspace*{3pt}
\[
\quad
\gamma
\off =\off 
\alpha\qff -\qff \beta
\]

\vspace*{-9pt}
is a cycle.\off\oss
Let\dss $\psi$\dss be a Torelli diffeomorphism of\dss $S$\dss 
extending\dss $\varphi$\snsp.\qff\oss
Then\qss $\psi_{\fff *\fff}(\gamma)$\qss
is ho\-mol\-o\-gous to\dss $\gamma$\dss in\dss $S$\nnsp.\oss
By the Mayer--Vietoris equivalence\vspace*{3pt}
\[
\quad
\psi_{\fff *\fff}(\gamma)\qff -\qff \gamma
\off =\off
\partial\fff(A\qff +\qff B)
\]

\vspace*{-9pt}
for some chains\qss 
$A\qff \in\qff C\two(Q)$\dnsp,\oss
$B\qff \in\qff C\two(\ccomp Q)$\dnsp.\oss 
Since\dss $\psi$\dss is an extension of\dss $\varphi$\nnsp,\oss
it follows that\vspace*{3pt}
\begin{equation}
\label{two-sides-e}
\quad
\varphi_{\fff *\fff}(\alpha)\qff -\qff \alpha\qff -\qff \partial A
\off =\off
\psi_{\fff *\fff}(\beta)\qff -\qff \beta\qff -\qff \partial B
\qff \in\qff C\one(c)\dff.
\end{equation}

\vspace*{-9pt}
Let us represent\qss $B$\qss as a sum\qss
$B\qff =\qff C\qff +\qff D$\qss
with\qss $C\qff \in\qff C\two(R\fff)$\qss
and\qss $D\dff \in\dff C\two(\ccomp Q\dff \smallsetminus\dff R\fff)$\dnsp.\oss
Then\vspace*{3pt}
\[
\quad
\varphi_{\fff *\fff}(\alpha)\qff -\qff \alpha\qff -\qff \partial\fff(A\qff +\qff D)
\off =\off
\psi_{\fff *\fff}(\beta)\qff -\qff \beta\qff -\qff \partial C
\qff \in\qff 
C\one(\partial R\fff)\dff.
\]

\vspace*{-9pt}
Let\oss
$\delta
\off =\off
\psi_{\fff *\fff}(\beta)\qff -\qff \beta\qff -\qff \partial C$\nnsp.\qff\oss
Then\dss $\delta$\dss is a cycle in\dss $\partial R$\dss
and the cycle\qss 
$\varphi_{\fff *\fff}(\alpha)\qff -\qff \alpha$\qss
is homologous in\qss $Q$\qss to\qss
$\delta\qff +\qff \partial D$\nnsp.\qff\oss
By the definition,\oss $\Delta_{\dff \varphi}\fff(\alpha)$\qss 
is equal to the image of the homology class\qss
$\hclass{\delta\qff +\qff \partial D}
\qff \in\qff
H\one(c)$\qss
in\dss the\dss group\qss $\mathbold{H}\one(c)$\dnsp.\oss

Since,\oss obviously,\oss 
$\hclass{\partial D}\qff \in\qff K\one(c\to S)$\dnsp,\oss
the image of\qss
$\hclass{\delta\qff +\qff \partial D}$\qss
in\qss
$\mathbold{H}\one(c)$\qss
is equal to the image of\qss $\hclass{\delta}$\qss
and hence belongs to\qss
$\mathbold{H}\one(\partial R\fff)$\dnsp.\oss
It follows that\oss
$\displaystyle
\Delta_{\dff \varphi}\fff(\alpha)
\qff \in\qff 
\mathbold{H}\one(\halfff\partial R\fff)$\dnsp.\oss  \eproof

\myitpar{Matrix presentations.} Let\qss $P$\qss be a component of\qss $\ccomp Q$\dnsp.\qff\oss
Let\qff\oss 
$C\zero\fff,\off\off C\one\fff,\off\off \ldots\fff,\off\off C_{\fff n}$\qff\oss
be the components of\qss $\partial P$\oss
oriented as components of\qss $\partial Q$\dnsp.\off\oss
For\oss 
$i
\off =\off 
\oldstylenums{1}\fff,\pff \oldstylenums{2}\fff,\pff \ldots\fff,\pff n$\oss let
\[
\quad
o_{\fff i}\off =\off o_{C_{\scriptstyle \dff i}}\qff -\qff o_{C_{\scriptstyle \dff \oold}}
\]
and let\qss $c_{\dff i}$ be the image of\qss $\hclass{C_{\fff i}}$\qss in\qss
$\mathbold{H}\one(\partial P\halfff)$\dnsp.\off\oss
Then\pff\oss 
\[
\quad
o\one\fff,\off\off o\two\fff,\off\off \ldots\fff,\off\off o_{\fff n}
\]  
is a basis of\pss $\mathbold{K}\zero(\partial P)$\pss and\oss 
\[
\quad
c\one\fff,\off\off c\two\fff,\off\off \ldots\fff,\off\off c_{\dff n}
\]
is a basis of\pss $\mathbold{H}\one(\partial P)$\dnsp.\oss

Suppose that\qss $\delta$\qss is completely reducible.\oss
Every
$\displaystyle
\delta\fff(o_{\fff i})$
belongs to\qss
$\mathbold{H}\one(\partial P)$\nsp,\oss 
and hence
\[
\quad
\delta\fff(o_{\fff i})\off 
=\off \sum\nolimits_{\qff j\qff =\qff \oldstylenums{1}}^{\dff n}\dff m_{\dff i\fff j}\trf c_{\dff j}
\]
for some integers\qss $m_{\qff i\dff j}$\nsp.\oss 
The matrix\qss 
$(\fff m_{\qff i\dff j}\fff)_{\fff \oldstylenums{1}\qff \leq\qff i\dff,\qff j\qff \leq\qff n}$\qss
is the matrix of the component\qss $\delta_{\dff P}$\dnsp. 

The symmetry condition
holds for\qss $\delta_{\dff P}$\qss  
if and only if it holds for every pair\qss
$o_{\fff i}\dff,\off o_{\fff j}$\qss of elements of our basis of\qss 
$\mathbold{K}\zero(\partial P)$\dnsp,\qff\oss
i.e.\qss if and only if
\[
\quad
\langle\qff o_{\fff i}\halfff\fff,\off \delta(o_{\fff j}) \qff\rangle
\off =\off
\langle\qff o_{\fff j}\halfff,\off \delta(o_{\fff i}) \qff\rangle
\]
for all\qss $i\fff,\pff j$\qss between\dss $\oldstylenums{1}$\dss and\dss $n$\nsp.\oss
But\qss 
$\langle\qff o_{\fff i}\halfff\fff,\off \delta(o_{\fff j}) \qff\rangle
\off =\off 
m_{\qff j\dff i}$\oss
and\oss
$\langle\qff o_{\fff j}\halfff\fff,\off \delta(o_{\fff i}) \qff\rangle
\off =\off 
m_{\qff i\dff j}$\dnsp.\oss
Therefore,\oss $\delta_{\dff P}$\qss is symmetric if and only if\qss 
$m_{\qff i\dff j}
\off =\off  
m_{\qff j\dff i}$\qss for all\qss 
$i\fff,\pff j$\qss 
between\dss $\oldstylenums{1}$\dss and\dss $n$\nsp,\oss
i.e.\qss if and only if the corresponding matrix\qss 
$(\fff m_{\qff i\dff j}\fff)_{\fff \oldstylenums{1}\qff \leq\qff i\fff,\qff j\qff \leq\qff n}$\qss 
is symmetric.\oss

\myitpar{Symmetric matrices.} Symmetric integer\qss $n\times n$\qss matrices
form a free abelian group\qss $SM\fff(n)$\qss of rank\qss 
$n(n\qff +\qff \oldstylenums{1})/\oldstylenums{2}$\nsp.\oss
A standard basis of\qss $SM\fff(n)$\qss consists of the following matrices.\oss
For a pair\oss
$k\fff,\pff l
\pff \in\pff 
\{\qff \oldstylenums{1}\fff,\off \oldstylenums{2}\fff,\off \ldots\fff,\off n \qff\}$\nsp,\qff\oss
let\qss $e\dff(k\dff, l)$\qss be\qss $n\times n$\qss matrix with entries
\[
\quad
e\dff(\fff k\dff, l\fff)_{\qff i\dff j}\off =\off \oldstylenums{1}
\hspace*{1.5em}\mbox{ if }\hspace*{1.5em} 
\{\qff i\fff,\pff j \qff\}
\off =\off 
\{\qff k\fff,\pff l \qff\}
\hspace*{1.5em}\mbox{ and }\hspace*{1.5em}
e\dff(\fff k\dff, l)_{\qff i\dff j}\off =\off \oldstylenums{0}
\hspace*{1.5em}\mbox{ otherwise. }\hspace*{1.5em}
\]
The matrices\qss $e\fff(k\dff, l)$\qss with\qss $k\dff\leq\dff l$\qss 
obviously form a basis of\qss $SM\fff(n)$\dnsp.\oss

For a subset\qss $\mathfrak{a}$\qss of\oss
$\displaystyle
\{\qff \oldstylenums{1}\fff,\off \oldstylenums{2}\fff,\off \ldots\fff,\off n \qff\}$\nsp,\qff\oss
let\qss $m\dff(\mathfrak{a})$\qss be\qss $n\times n$\qss matrix with the entries
\[
\quad
m\dff(\mathfrak{a})_{\qff i\dff j}
\off =\off 
\oldstylenums{1}
\hspace*{1.5em}\mbox{ if }\hspace*{1.5em} i\fff,\pff j\qff \in\qff \mathfrak{a}\dff,
\hspace*{1.5em}\mbox{ and }\hspace*{1.5em} 
m\dff(\mathfrak{a})_{\qff i\dff j}
\off =\off 
\oold
\hspace*{1.5em}\mbox{ otherwise. }
\]
For\qss $k\qff \leq\qff l$\nsp,\qff\oss let\oss 
$\displaystyle
\mathfrak{a}_{\qff k\fff,\pff l}
\off =\off 
\{\qff k\fff,\off k\qff +\qff \oldstylenums{1}\fff,\off \ldots\fff,\off l \qff\}$\dnsp,\oss
and\dss let\oss
$m\dff(\fff k\fff,\pff l\fff)
\off =\off 
m\dff(\dff \mathfrak{a}_{\qff k\fff,\pff l}\fff)$\dnsp.\off\oss

\mypar{Lemma.}{other-basis} \emph{The matrices\qss $m\fff(\fff k\fff,\pff l\fff)$\qss 
with\qss $k\qff \leq\qff l$\qss form a basis of\pss $SM\fff(n)$\dnsp.\oss}

\proof\qss Since\qss $SM\fff(n)$\qss is a free abelian group of rank\qss 
$n(n\qff +\qff \oldstylenums{1})/\oldstylenums{2}$\qss 
and there are\qss $n(n\qff +\qff \oldstylenums{1})/\oldstylenums{1}$\qss matrices\qss
$m\dff(\fff k\fff,\pff l\fff)$\qss with\qss $k\dff\leq\dff l$\snsp,\oss
it is sufficient to prove that these matrices generate the abelian group\qss $SM\fff(n)$\dnsp.\oss

Therefore it is sufficient to express the matrices\qss
$e\dff(k\dff, l)$\qss with\qss $k\qff \leq\qff l$\qss 
in terms of the matrices\qss $m\dff(\fff k\fff,\pff l\fff)$\qss 
with\qss $k\qff \leq\qff l$\nnsp.\off\oss
But\oss\vspace*{3pt} 
\[
\quad
e\dff(k\dff, k)
\off =\off
m\dff(\fff k\fff,\pff k\fff)\fff,
\] 

\vspace*{-36pt}
\[
\quad
e\dff(k\dff, k\qff +\qff 1)
\off =\off 
m\dff(\fff k\fff,\pff k\qff +\qff 1\fff)
\qff -\qff 
m\dff(\fff k\fff,\pff k\fff)
\qff -\qff 
m\dff(\fff k\qff +\qff 1\fff,\pff k\qff +\qff 1\fff)\fff,
\] 

\vspace*{-9pt}
and\dss if\oss $k\qff +\qff 1\qff <\qff l$\nsp,\qff\oss then\oss\vspace*{3pt}
\[
\quad
e\fff(k\dff, l)\off 
=\off m\fff(k\fff,\pff l)\qff 
-\qff m\fff(k\qff +\qff 1\fff,\pff l)\qff 
-\qff m\fff(k\fff,\pff l\qff -\qff 1)\qff 
+\qff m\fff(k\qff +\qff 1\fff,\pff l\qff -\qff 1)\fff.
\] 

\vspace*{-9pt}
It follows that the matrices\qss $m\dff(\fff k\fff,\pff l\fff)$\qss
with\qss $k\qff \leq\qff l$\qss generate\qss $SM\fff(n)$\dnsp,\oss
and hence form a basis of\qss $SM\fff(n)$\dnsp.\oss  \eproof

\myitpar{Peripheral twists and bi-twists.} Let\qss $P$\qss
be a component of\qss $\ccomp Q$\dnsp,\oss
and let\qss $U$ be the union of
several components of\qss $\partial P$\dnsp.\oss
Let us orient\qss $U$\qss as a part of
the boundary of\qss $Q$\dnsp.\oss

There are circles\qss $C$\qss in\qss $Q$\qss 
bounding together with\qss $U$\qss
a subsurface of\qss $Q$\dnsp.\oss 
We will call twist diffeomorphisms about such circles the\pss 
\emph{peripheral twists diffeomorphisms about}\pss $U$\dnsp.\oss

Let\qss $C$\qss be a circle in\qss $Q$\qss bounding together with\qss $U$\qss
a subsurface of\qss $Q$\dnsp,\oss 
and let\qss $D$\qss be a circle in\qss $\ccomp Q$\qss bounding together with\qss $U$\qss
a subsurface of\qss $\ccomp Q$\dnsp.\oss
The union of these two subsurfaces is a 
subsurface of\qss $S$\qss with boundary\qss $C\dff\cup\dff D$\dnsp.\off\oss
Let\oss $T_{\dff C}$\nsp,\oss $T_{\dff D}$\oss  
be twists diffeomorphisms of\qss $S$\qss about\pss $C$\dnsp,\oss $D$\pss 
fixed on\qss $\ccomp Q$\nsp,\oss $\ccomp P$\qss respectively.\qff\oss

The diffeomorphism\oss 
$\displaystyle
\psi
\off =\off 
T_{\dff C}\dff \circ\dff T_{\dff D}^{\qff -\dff \oldstylenums{1}}$\oss
leaves\dss $Q$\dss invariant and is equal to the identity on\dss $c$\dnsp.\oss
The isotopy class\oss 
$t_{\dff C}\dff \circ\dff t_{\dff D}^{\qff -\dff \oldstylenums{1}}$\oss
of\qss $\psi$\qss is a Dehn-Johnson twist 
and hence belongs to\dss $\tors$\dnsp.\oss
Therefore\qss $\psi$\qss is a Torelli diffeomorphism.\oss

We will call such diffeomorphisms\qss $\psi$\qss
\emph{peripheral bi-twist diffeomorphisms}\pss about\qss $U$\dnsp.\oss

The diffeomorphism\oss 
$\psi
\off =\off 
T_{\dff C}\dff \circ\dff T_{\dff D}^{\qff -\dff \oldstylenums{1}}$\oss 
induces a weak Torelli diffeomorphism\dss $\varphi$\dss of\dss $Q$\dnsp.\oss  
It is a twist diffeomorphism of\qss $Q$\qss about\qss $C$\dnsp,\oss
and its extension\qss $\chi$\qss by the identity to a diffeomorphism of\qss $S$\qss
is equal to\qss $T_{\dff C}$\nsp.\oss

\mypar{Lemma.}{twist-torelli-map}\oss \emph{Let\oss 
$\displaystyle
\hclass{U}\qff\in\qff \mathbold{H}\one(c)$\qss 
be the image of the fundamental class of\oss $U$\nsp.\oss 
Then\oss
\[
\quad
\delta_{\fff \varphi}\fff(a)
\off =\off 
\langle\qff a\halfff\fff,\pff \hclass{U} \pff\rangle_{\dff c}\qff\hclass{U}
\]
for all}\oss $a\dff\in\dff \mathbold{K}\zero(c)$\dnsp.\oss 

\proof\qss The orientation of\qss $U$\qss 
is the same as the orientation of\qss $U$\qss as a part of the boundary
of the subsurface of\dss $S$\dss bounded by\qss $C\dff\cup\dff U$\dnsp.\oss
Let us consider\qss $C$\qss with the orientation opposite to 
its orientation as a part of this boundary.\oss
Then\oss\vspace*{3pt}
\[
\quad
\hclass{C}\off =\off \hclass{U}
\] 

\vspace*{-9pt}
after the identification of\qss $\mathbold{H}\one(c)$\qss
with the image of\qss $H\one(c)$\qss in\qss $H\one(S)$\dnsp.\oss

Let\oss $a\qff \in\qff H\one(S)$\dnsp.\oss
By\dss Theorem\qss \ref{extension-of-weakly-torelli}\oss\vspace*{3pt}
\[
\quad
\delta_{\fff \varphi}\dff \circ\pff \mathcal{D}\dff(\fff a\fff)
\off =\off
\chi_{\fff *\fff}(\fff a\fff)\qff -\qff a\dff.
\]

\vspace*{-9pt}
Since\dss $\chi$\dss is a twist diffeomorphism of\dss $S$\dss 
about\dss $C$\snsp,\off\oss\vspace*{3pt}
\[
\quad
\chi_{\fff *\fff}(\fff a\fff)\qff -\qff a
\off =\off
\langle\qff a\dff,\pff\dff \hclass{C} \qff\rangle \qff\hclass{C}
\off =\off
\langle\qff a\dff,\pff\dff \hclass{U} \qff\rangle \qff\hclass{U}\dff.
\]

\vspace*{-9pt}
Since\oss
$\displaystyle
\langle\qff a\dff,\pff\dff \hclass{U} \qff\rangle
\off =\off
\langle\qff \mathcal{D}\dff(\fff a\fff)\fff,\pff\fff \hclass{U} \qff\rangle_{\dff c}$\oss
by the formula\oss (\ref{mv-induced-pairing}),\oss
it follows that\oss\vspace*{3pt}
\[
\quad
\delta_{\fff \varphi}\dff \circ\pff \mathcal{D}\dff(\fff a\fff)
\off =\off
\langle\qff \mathcal{D}\dff(\fff a\fff)\fff,\pff\fff \hclass{U} \qff\rangle_{\dff c} \qff \hclass{U}\dff.
\]

\vspace*{-9pt}
Since\qss $\mathcal{D}$\qss is a map onto\oss $\mathbold{K}\zero(c)$\dnsp,\oss
the lemma follows.\oss  \eproof

\mypar{Theorem.}{realization-of-delta}\oss 
\emph{For every symmetric completely reducible map}\oss\vspace*{3pt}
\[
\quad
\delta\dff \colon\qff
\mathbold{K}\zero(c)\ttoo \mathbold{H}\one(c)
\]

\vspace*{-9pt}
\emph{there exist a weakly Torelli diffeomorphism\pss $\varphi$\pss of\pss $Q$\pss
admitting an extension to a Torelli diffeomorphism of\pss $S$\pss 
and such that\oss 
$\displaystyle
\delta_{\fff \varphi}\off =\off \delta$\nsp.\oss}

\proof\qss 
Suppose that\trs $\delta$\trs is equal to\dss $\oold$\dss on all summands\sss of\pss
$\mathbold{K}\zero(c)$\pss except the summand\pss 
$\mathbold{K}\zero(\partial P\halfff)$\pss for some
component\dss $P$\dss of\dss $\ccomp Q$\dnsp.\qff\oss
Then\qss $\delta$\qss maps\pss
$\mathbold{K}\zero(\partial P\halfff)$\pss into\pss
$\mathbold{H}\one(\partial P\halfff)$\dnsp.  

We will use the notations from the discussion of matrix presentations
and symmetric matrices above.\oss
For a subset\qss\vspace*{3pt} 
\[
\quad
\mathfrak{A}
\off \subset\off
\{\qff \oldstylenums{1}\fff,\off\qff \oldstylenums{2}\fff,\off\qff \ldots\fff,\off\qff n\dff\}
\]

\vspace*{-9pt}
we will denote by\qss $U\dff(\mathfrak{A})$\qss be the union of circles\qss 
$C_{\dff i}$\qss with\pss 
$i\qff \in\qff \mathfrak{A}$\nsp.\oss
Let us choose for each such subset\qss $\mathfrak{A}$
a peripheral bi-twist diffeomorphism\qss $\psi\dff(\mathfrak{A})$\qss 
about\qss $U\dff(\mathfrak{A})$\dnsp.\oss
Let\qss $\varphi\fff(\mathfrak{A})$\qss be the diffeomorphism of\dss $Q$\dss
induced by\qss $\psi\dff(\mathfrak{A})$\dnsp.\oss
It is a peripheral twist diffeomorphism about\qss $U\dff(\mathfrak{A})$\dnsp.\oss

Lemma\qss \ref{twist-torelli-map}\qss implies that the\dss $\delta$\dnsp-difference map of\qss
$\varphi\dff(\mathfrak{A})$\qss is equal to\dss $\oold$\dss on all summands of\qss
$\mathbold{K}\zero(c)$\qss except of\qss 
$\mathbold{K}\zero(\partial P)$\dnsp,\oss
and the matrix of the map\vspace*{3pt}
\[
\quad
\mathbold{K}\zero(\partial P)\ttoo H\one(\partial P)_{\fff P}
\]

\vspace*{-9pt}
induced by this\dss $\delta$\dnsp-difference map is equal to\qss $m\dff(\mathfrak{A})$\dnsp.\oss
Since every symmetric matrix is an integer linear combination of matrices of this form,\oss
Lemmas\qss \ref{weakly-torelli-composition}\qss and\qss \ref{weakly-torelli-inverse}\qss 
imply that\qss $\delta$\qss
is equal to the\dss $\delta$\dnsp-difference map of
a product\dss $\varphi$\dss of diffeomorphisms of the form\qss 
$\varphi\dff(\mathfrak{A})$\dnsp.\oss
The corresponding product\dss $\psi$\dss of diffeomorphisms\qss 
$\psi\dff(\mathfrak{A})$\qss
is a Torelli diffeomorphism extending\dss $\varphi$\snsp.\oss
This proves the theorem when\qss $\delta$\qss is equal to\dss $\oold$\dss
on all summands of\qss $\mathbold{K}\zero(c)$\qss except one.\oss

Any completely reducible map\qss $\delta$\qss is equal to 
the sum over the components\qss $P$\qss of maps equal to\dss $\oold$\dss
on all summands\dss of\pss
$\mathbold{K}\zero(c)$\pss except\dss of\oss 
$\mathbold{K}\zero(\partial P\halfff)$\dnsp.\oss
Moreover\halfff,\oss 
if\qss $\delta$\qss is symmetric,\oss
then all summands of this sum are also symmetric.\oss
Therefore,\oss the general case of the theorem 
follows from the already proved case 
and Lemmas\qss \ref{weakly-torelli-composition}\qss 
and\qss \ref{weakly-torelli-inverse}.\oss  \eproof

\mypar{Theorem.}{char-ext} \emph{For a diffeomorphism\oss 
$\varphi\qff \colon\dff Q\toto Q$\oss 
fixed on\dss $c$\dss
the following two conditions are equivalent\halfff.}\vspace*{0pt}

({\fff}i{\fff}) $\varphi$\pss \emph{can be extended to a Torelli diffeomorphism of\pss $S$\dnsp.}

({\fff}ii{\fff}) $\varphi$\pss \emph{a weakly Torelli diffeomorphism and its\dss 
$\delta$\dnsp-difference map\qss $\delta_{\fff \varphi}$\qss is
completely reducible\halfff.}

\proof\qss By remarks after the definition of weakly Torelli diffeomorphisms 
and Theorem\qss \ref{torelli-is-completely-reducible}\qss
the condition\qss ({\fff}i{\fff})\qss implies\qss ({\fff}ii{\fff}).\oss

Conversely,\oss suppose that\dss $\varphi$\dss satisfies\qss ({\fff}ii{\fff}).\oss
Then by Theorem\qss \ref{symmetry-theorem}\qss the\dss 
$\delta$\dnsp-difference map\qss $\delta_{\fff \varphi}$\qss
is symmetric,\oss
and hence
by Theorem\qss \ref{realization-of-delta}\qss there exists
a peripheral diffeomorphism\qss $\varphi'$\qss of\dss $Q$\dss
admitting an extension\dss $\psi'$\dss to a Torelli diffeomorphism of\dss $S$\dss
and such that\qss  
\[
\quad
\delta_{\dff F'}
\off =\off 
\delta_{\dff F}\qff.
\]
Lemmas\qss \ref{weakly-torelli-composition}\qss and\qss \ref{weakly-torelli-inverse}\qss imply
that the map\oss\vspace*{3pt} 
\[
\quad
\varphi''
\off =\off 
(\dff \varphi'\dff)^{\dff -\dff \oldstylenums{1}}\qff \circ\qff \varphi
\]

\vspace*{-9pt} 
is peripheral and its\dss $\delta$\dnsp-difference map is\dss $\oold$\nnsp.\oss
Let\qss $\chi''\qff =\qff \varphi''\dff\backslash\fff S$\qss 
be the extension of\qss $\varphi''$\qss to\dss $S$\dss by the identity.\oss 
By Theorem\qss \ref{extension-by-identity}\qss $\chi''$\qss is a Torelli diffeomorphism.\oss
Let\oss\vspace*{3pt} 
\[
\quad
\psi
\off =\off 
\psi'\circ \chi''\dff.
\]

\vspace*{-9pt} 
By Lemma\qss \ref{weakly-torelli-composition}\pss $\psi$\qss is a Torelli diffeomorphism.\qff\oss
Since\qss $\psi'$\qss is an extension of\qss $\varphi'$\qss and\qss $\chi''$\qss
is an extension of\oss 
$(\dff \varphi'\dff)^{\dff -\dff  \oldstylenums{1}}\qff \circ\qff \varphi$\nnsp,\qff\oss
the diffeomorphism\dss $\psi$\dss is an extension of\dss $\varphi$\nnsp.\qff\oss
Therefore\qss $\varphi$\qss satisfies\qss ({\fff}i{\fff}).\oss  \eproof

\mysection{Multi-twists\qss about\qss the\qss boundary}{multi-twists-about-the-boundary}

\vspace*{\medskipamount}
\myitpar{Multi-twist diffeomorphism about the boundary.}
Suppose that for each a component\dss $C$\dss of\dss $c$\dss
an annulus\dss $A_{\fff C}$\dss in\dss $Q$\dss having\dss $C$\dss as one of 
its boundary components is fixed,\oss
and suppose that these annuli are pair-wise disjoint\halfff.\oss
Let\pss\vspace*{3pt} 
\[
\quad
\tau_{\dff C}\qff \colon\qff Q\ttoo Q
\] 

\vspace*{-9pt}
be the extension by the identity of a left twist diffeomorphism of\dss $A_{\fff C}$\snsp.\oss
A\qss \emph{multi-twist diffeomorphism of\pss $Q$\pss
about the boundary}\qss is a diffeomorphism of\dss $Q$\dss of the form\oss\vspace*{4pt} 
\begin{equation}
\label{boundary-multi-twist}
\quad
\tau\off =\off \prod\nolimits_{\dff C}\qff \tau_{\dff C}^{\qff m_{\dff C}}\qff,
\end{equation}

\vspace*{-8pt}
where\dss $C$\dss runs over the components of\qss $c\qff =\qff \partial Q$\qss
and\qss $m_{\dff c}$\qss are integers.\oss

The extension\qss $\tau\trf\backslash\fff S$\qss of such a diffeomorphism\dss 
$\tau$\dss by the identity
is a multi-twist diffeomorphism of\dss $S$\dss about\dss $c$\dss
and its isotopy class is equal to the Dehn multi-twist\oss\vspace*{4pt}
\begin{equation}
\label{boundary-dehn-multi-twist}
\quad
t
\off =\off 
\prod\nolimits_{\dff C}\qff t_{\dff C}^{\qff m_{\dff C}}
\off \in\off 
\mms\dff,
\end{equation} 

\vspace*{-8pt}
where each\qss $t_{\dff C}\qff \in\qff \mms$\qss 
is the left Dehn twist of\dss $S$\dss about\dss $C$\dnsp.\oss

\myitpar{Diagonal maps.} A homomorphism\oss
\[
\quad
D\dff\colon\dff H\zero(c)\ttoo \mathbold{H}\one(c)
\]
is said to be\qss \emph{diagonal}\qss if for every component\dss $C$\dss of\dss $c$\dss
it maps the generator\qss $o_C$\qss of\qss $H\zero(c)$\qss
from Section\qss \ref{general-extensions}\qss to an integer multiple of\qss
$\hclass{C}$\nnsp.\oss 
Every multi-twist diffeomorphism\sss $\tau$\sss about the boundary
tautologically defines a diagonal map\qss $D_{\dff \tau}$\nnsp.\qff\oss
Namely,\oss if\qss $\tau$\qss has the form\qss (\ref{boundary-multi-twist}),\qff\oss
then\qss $D_{\dff \tau}$\qss is given the rule
\[
\quad
D_{\dff \tau}
\dff \colon\dff 
o_C\off \longmapsto\off m_{\dff C}\dff\hclass{C}\dff.
\]

\vspace*{-6pt}
\mypar{Lemma.}{difference-multi-twist} \emph{Suppose that\qss 
$\tau$\qss
is a multi-twist diffeomorphism of\qss $Q$\qss
about the boundary.\oss 
Then\qss $\tau$\qss can be extended to 
a Torelli diffeomorphism of\qss $S$\qss and hence is a weakly Torelli diffeomorphism.\oss
If\oss $\tau$\qss is given by 
the formula\oss \textup{(\ref{boundary-multi-twist})},\oss
then\dss its\dss difference\sss map}\qss\vspace*{4pt} 
\[
\quad
\delta_{\dff \tau}
\qff \colon\qff 
\mathbold{K}\zero(c)\ttoo \mathbold{H}\one(c)
\]

\vspace*{-8pt} 
\emph{is equal to the restriction to\qss $\mathbold{K}\zero(c)$\qss of the diagonal map}\oss
$\displaystyle
D_{\dff \tau}
\qff \colon\qff 
H\zero(c)\ttoo H\one(c)$\dnsp.\oss

\proof\qss It is sufficient to consider the diffeomorphism\qss
(\ref{boundary-multi-twist}).\oss
For each component\dss $C$\dss of\dss $c$\dss let\dss
$B_{\fff C}$\dss be an annulus in\dss $\ccomp Q$\dss having\dss $C$\dss as 
one of its boundary components,\oss and let\qss 
$T_{\dff C}$\qss be the extension by the identity to\dss $S$\dss 
of a left twist diffeomorphism of\qss $B_{\fff C}$\nsp.\qff\oss
Let\oss\vspace*{3pt}
\[
\quad
T\off =\off \prod\nolimits_{\dff C}\qff T_{\dff C}^{\qff m_{\dff C}}\qff.
\]

\vspace*{-9pt}
Then\dss $T$\dss is equal to the identity on\dss $Q$\dnsp,\oss
and hence\qss $T^{\fff -\halfff \oldstylenums{1}}\circ\tau$\qss
extends\qss $\tau$\dnsp.\oss
On the other hand,\oss $T^{\fff -\halfff \oldstylenums{1}}\circ\tau$\qss is isotopic
to the identity,\oss and hence is a Torelli diffeomorphism.\oss
Therefore,\oss $\tau$\dss can be extended to a Torelli diffeomorphism
of\dss $S$\dss and hence is a weakly Torelli diffeomorphism.\oss

The action of the Dehn multi-twist\qss (\ref{boundary-dehn-multi-twist})\qss 
on\qss $H\one(S)$\qss is given by the formula\oss\vspace*{4pt}
\[
\quad
t_{\dff *\fff}(\fff a\fff)
\off =\off
a\qff +\qff
\sum\nolimits_{\dff C}\qff m_{\dff C}\trf 
\langle\qff a\dff,\pff \hclass{C} \qff\rangle\dff\hclass{C}\dff,
\]

\vspace*{-8pt}
where\qss $a\qff \in\qff H\one(S)$\dnsp.\oss
Theorem\qss \ref{extension-of-weakly-torelli}\qss implies that
$\displaystyle
t_{\dff *\fff}(\fff a\fff)\qff -\qff a
\off =\off
\delta_{\dff \tau}\trf \circ\pff \mathcal{D}\dff(a)$\dnsp,\oss
and hence\oss\vspace*{4pt}
\[
\quad
\delta_{\dff \tau}\trf \circ\pff \mathcal{D}\dff(a)
\off =\off
\sum\nolimits_{\dff C}\qff m_{\dff C}\qff 
\langle\qff a\qff,\pff \hclass{C} \pff\rangle \trf\hclass{C}\dff.
\]

\vspace*{-8pt}
Let us represent\qss $\mathcal{D}(a)$\qss
as a linear combination\vspace*{3pt}
\[
\quad
\mathcal{D}\dff(a)
\off =\off 
\sum\nolimits_{\dff C}\qff n_{\dff C}\qff o_C
\]

\vspace*{-9pt}
of generators\dss $o_C$\dss with integer coefficients\qss $n_{\dff C}$\snsp.\oss
Then\qss (\ref{zero-one-duality})\qss implies that\vspace*{3pt}
\[
\quad
n_{\dff D} 
\off =\off 
\sum\nolimits_{\dff C}\qff 
\langle\qff n_{\dff C}\qff o_C\qff,\pff\fff \hclass{D} \qff\rangle_{\dff c}
\off =\off 
\langle\qff \mathcal{D}\dff(a)\fff,\pff\fff \hclass{D} \qff\rangle_{\dff c}
\]

\vspace*{-9pt}
for every component\dss $D$\dss of\dss $c$\dnsp.\oss
On the other hand,\oss
$\dis
\langle\qff \mathcal{D}\dff(a)\fff,\pff\fff \hclass{D} \qff\rangle_{\dff c}
\off =\off
\langle\qff a\fff,\off \hclass{D} \qff\rangle$\oss
by\qss (\ref{mv-induced-pairing})\qss
and hence\oss
$n_{\dff D} 
\off =\off
\langle\qff a\fff,\off \hclass{D} \qff\rangle$\oss
for every component\dss $D$\dss of\dss $c$\dss and\vspace*{4pt}
\[
\quad
\mathcal{D}\dff(a)
\off =\off 
\sum\nolimits_{\dff D}\qff 
\langle\qff a\fff,\off \hclass{D} \qff\rangle\qff o_{\fff D}
\off =\off 
\sum\nolimits_{\dff C}\qff 
\langle\qff a\fff,\off \hclass{C} \qff\rangle\qff o_C\qff.
\]

\vspace*{-8pt} 
After applying\dss $D_{\dff \tau}$\dss to the last equality
we get\vspace*{4pt}
\[
\quad
D_{\dff \tau}\dff \circ\pff \mathcal{D}\dff(a)
\off =\off
\sum\nolimits_{\dff C}\qff m_{\dff C}\qff 
\langle\qff a\fff,\off \hclass{C} \qff\rangle \qff\hclass{C}\dff.
\]

\vspace*{-8pt}
By comparing this formula with the above formula for\pss 
$\delta_{\dff \tau}\dff \circ\pff \mathcal{D}\dff(a)$\nsp,\oss
we see that\vspace*{2pt}
\[
\quad
\delta_{\dff T}\dff \circ\pff \mathcal{D}
\off =\off
D_{\dff \tau}\dff \circ\pff \mathcal{D}\dff.
\]

\vspace{-10pt}
Since\qss 
$\mathcal{D}$\qss is a map onto\qss $\mathbold{K}\zero(c)$\dnsp,\pss
the lemma follows.\oss  \eproof

\mypar{Theorem.}{twist-correctable} 
\emph{For a diffeomorphism\qss $\varphi\dff \colon\qff Q\toto Q$\pss fixed on\dss $c$\dss
the following two conditions are equivalent.}

({\fff}i{\fff})
\emph{There exists a multi-twist diffeomorphism\dss 
$\tau\qff \colon\qff S\toto S$\dss about\dss $c$\dss 
such that\qss $\tau\qff \circ\qff (\varphi\fff\backslash S\halfff)$\qss
is a Torelli diffeomorphism.}

({\fff}ii{\fff}) 
\emph{\dnsp$\varphi$\qss a weakly Torelli diffeomorphism and\qss
$\delta_{\fff \varphi}$\qss is the restriction of a diagonal map.\oss}

\proof\qss Suppose that\qss ({\fff}i{\fff})\qss holds.\oss
We may assume that\qss $\tau$\qss is
the extension by the identity\qss $\rho\dff\backslash S$\qss
of a multi-twist diffeomorphism\dss $\rho$\dss of\dss $Q$\dss about the boundary.\oss
Then\vspace*{2pt} 
\[
\quad
(\fff\rho\circ \varphi\fff)\fff\backslash\fff S
\off =\off
\tau\dff\circ\fff (\fff \varphi\fff\backslash\fff S\halfff)
\]

\vspace*{-10pt}
is a Torelli diffeomorphism and hence Theorem\qss \ref{extension-by-identity}\qss
implies that\qss $\rho\circ \varphi$\qss is a weakly Torelli diffeomorphism and its
difference map is equal to\dss $\oold$\dnsp.\oss 
Since\dss $\rho$\dss is a weakly Torelli diffeomorphism,\oss
Lemmas\qss \ref{weakly-torelli-composition}\qss and\qss \ref{weakly-torelli-inverse}\qss imply
that\dss $\varphi$\dss is also a weakly Torelli diffeomorphism
and\oss $\delta_{\fff \varphi}\off =\off -\qff \delta_{\fff \rho}$\nsp.\oss
Since\qss $\delta_{\dff \rho}$\qss is equal to the restriction of a diagonal map
by Lemma\qss \ref{difference-multi-twist},\oss
it follows that the same is true for\qss $\delta_{\dff \varphi}$\nsp.\oss
This proves that\qss ({\fff}ii{\fff})\qss holds.\oss

Suppose that\qss ({\fff}ii{\fff})\qss holds.\oss
Then there is a multi-twist diffeomorphism\qss 
$\tau$\qss of\qss $Q$\qss about the boundary
such that\oss $\delta_{\fff \tau}\off =\off -\qff \delta_{\fff \varphi}$\nsp.\oss
In this case\qss $\tau\qff \circ\qff \varphi$\qss is a weakly Torelli diffeomorphism and its difference map
is equal to\dss $\oold$\dnsp.\oss
By Theorem\qss \ref{extension-by-identity}\qss the extension by the identity\qss 
$(\fff \tau\qff \circ\qff \varphi\fff)\backslash S$\qss
is a Torelli diffeomorphism of\qss $S$\dnsp.\oss
But\oss\vspace*{2pt} 
\[
\quad
(\fff \tau\qff \circ\qff \varphi\fff)\backslash S
\off =\off 
(\fff \tau\dff\backslash S\fff)\qff \circ\qff (\fff \varphi\fff\halfff\backslash S\fff)
\]

\vspace*{-10pt}
and\qss $\tau\dff\backslash S$\qss is a multi-twist about\qss $\partial Q$\dnsp.\oss
This proves that\qss ({\fff}i{\fff})\qss holds.\oss  \eproof

\mypar{Theorem.}{3-diagonal} \emph{If\oss $\partial P$\oss
consists of\oss $\leqslant\qff 3$\qss circles for every component\qss $P$\qss of\qss
$\ccomp Q$\nnsp,\oss 
then every completely reducible symmetric map\fff\qss
$\delta
\qff \colon\qff 
\mathbold{K}\zero(c)\ttoo \mathbold{H}\one(c)$\fff\qss 
is the restriction of a diagonal map.}

\proof\qss Suppose that\qss $\delta$\qss is equal to\dss $\oold$\qss on all summands\dss of\pss
$\mathbold{K}\zero(c)$\pss except the summand\pss 
$\mathbold{K}\zero(\partial P\halfff)$\pss for some
component\qss $P$\qss of\qss $\ccomp Q$\nnsp.\qff\oss
We will use the notations from the discussion of matrix presentations
and symmetric matrices from Section\qss \ref{general-extensions}.\oss

If\qss $\partial P$\qss is connected\halfff,\oss
then\oss
$\displaystyle
\mathbold{K}\zero(\partial P)
\off =\off
\mathbold{H}\one(\partial P)
\off =\off
\oold$\oss 
and the theorem is trivial.\oss

If\qss $\partial P$\qss consists of\dss $\oldstylenums{2}$\dss components\qss 
$C\zero\fff,\off C\one$\snsp,\oss 
then\qss $\mathbold{K}\zero(\partial P)$\qss 
is an infinite cyclic group with the generator\oss
$\displaystyle
o\one
\off =\off 
o_{C_{\scriptstyle \dff \oldstylenums{1}}}\qff -\qff o_{C_{\scriptstyle \dff \oldstylenums{0}}}$\nsp,\oss
and\qss $\mathbold{H}\one(\partial P)$\qss is an infinite cyclic group with the generator\oss
$\hclass{C\zero}\qff =\qff \hclass{C\one}$\snsp.\oss
Let\qss $\delta\dff(o\one)\qff =\qff m\dff \hclass{C\one}$\snsp.\oss 
If\qss $n\zero\fff,\pff n\one\qff \in\qff \zzz$\qss
and\qss $D$\qss is the\qss (unique)\qss diagonal map 
such that\oss
$\displaystyle
D\dff(o_C)\qff =\qff \oldstylenums{0}$\oss for\oss 
$\displaystyle
C\qff \neq\qff C\zero\fff,\pff C\one$\oss
and\oss\vspace*{3pt}
\[
\quad
D(o_{C_{\scriptstyle \dff \oldstylenums{0}}})\off =\off n\zero\dff \hclass{C\zero}\dff,
\quad\off
D(o_{C_{\scriptstyle 
\dff \oldstylenums{1}
}
})\off =\off 
n\dff \hclass{C\one}\dff,
\]

then\dss $\delta$\dss is equal to the restriction of\qss $D$\qss 
if and only if\qss $m\qff =\qff n\one\qff -\qff n\zero$\snsp.\oss
It follows that if\qss $\partial P$\qss
consists of\dss $\oldstylenums{2}$\dss components,\oss then $\delta$\qss is equal to the restriction of a diagonal map.\oss

If\qss $\partial P$\qss consists of\dss $\oldstylenums{3}$\dss components\qss 
$C\zero\fff,\off C\one\fff,\off C\two$\snsp,\oss 
then\qss $\mathbold{K}\zero(\partial P\halfff)$\qss is free abelian and\vspace*{2pt}
\[
\quad
o\one
\off =\pff 
o_{C_{\scriptstyle \dff \oldstylenums{1}}}\qff -\qff o_{C_{\scriptstyle \dff \oldstylenums{0}}}\dff,
\hspace*{1em}
o\two
\off =\off 
o_{C_{\scriptstyle \dff \oldstylenums{2}}}\qff -\qff o_{C_{\scriptstyle \dff \oldstylenums{0}}}\dff
\]

\vspace*{-10pt}
is a basis of\oss
$\mathbold{K}\zero(\partial P\halfff)$\dnsp.\oss
The group\qss $\mathbold{H}\one(\partial P\halfff)$\qss is free abelian
having\qss $\hclass{C\one}$\snsp,\oss $\hclass{C\two}$\qss as a basis.\oss
The induced map\oss
$\displaystyle
\delta_{\dff P}
\qff \colon\qff 
\mathbold{K}\zero(\partial P\halfff)\ttoo 
\mathbold{H}\one(\partial P\halfff)$\oss 
is symmetric and hence has the form\oss\vspace*{2pt}
\[
\quad
o\one
\off \longmapsto\off
m\one\dff \hclass{C\one}\qff +\qff m\dff \hclass{C\two}
\]

\vspace*{-37pt}
\[
\quad
o\two
\off \longmapsto\off
m\dff \hclass{C\one}\qff +\qff m\two\dff \hclass{C\two}\dff.
\]

\vspace*{-10pt}
Let\qss $D$\qss be the diagonal map such that\oss
$\displaystyle
D\dff(o_C)\qff =\qff \oold$\oss if\oss
with\oss $C\qff \neq\qff C\zero\fff,\pff C\one\fff,\pff C\two$\oss
and\oss\vspace*{2.5pt}
\[
\quad
D\dff(o_{C_{\scriptstyle \dff \oldstylenums{0}}})
\off =\off 
n\zero\dff \hclass{C\zero}\dff,
\quad\off
D\dff(o_{C_{\scriptstyle \dff \oldstylenums{1}}})
\off =\off 
n\one\dff \hclass{C\one}\dff,
\quad\off
D\dff(o_{C_{\scriptstyle \dff \oldstylenums{2}}})
\off =\off 
n\two\dff \hclass{C\two}\dff.
\]

\vspace*{-9.5pt}
As in Section\qss \ref{general-extensions},\oss we assume that circles\qss
$C\zero\fff,\off C\one\fff,\off C\two$\qss are oriented as components
of the boundary of\dss $Q$\dnsp.\qff\oss
Then\oss
$\displaystyle
\hclass{C\zero}\qff +\qff \hclass{C\one}\qff +\qff \hclass{C\two}
\off =\off
\oold$\snsp,\qff\oss and hence\oss\vspace*{3pt}
\[
\quad
D\dff(o\one)\off =\off n\one\dff \hclass{C\one}\qff -\qff n\zero\dff \hclass{C\zero}
\off =\off
(n\one\qff +\qff n\zero)\dff \hclass{C\one}\qff +\qff n\zero\dff \hclass{C\two}\dff, 
\]

\vspace*{-36pt}
\[
\quad
D\dff(o\two)\off =\off n\two\dff \hclass{C\two}\qff -\qff n\zero\dff \hclass{C\zero}
\off =\off
n\zero\dff \hclass{C\one}\qff +\qff (n\two\qff +\qff n\zero)\dff \hclass{C\two}\dff. 
\]

\vspace*{-9pt}
It follows that\qss $\delta$\qss is equal to the restriction of\qss $D$\qss if and only if\vspace*{2pt}
\[
\quad
m\one\off =\off n\one\qff +\qff n\zero\dff,
\quad\off
m\off =\off  n\zero\dff,
\quad\off
\quad
m\two\off =\off n\two\qff +\qff n\zero\dff.
\]

\vspace*{-10pt}
For any given\pss $m\fff,\pff m\one\fff,\pff m\two\qff\in\qff \zzz$\pss
this system of equation for\pss $n\zero\fff,\pff n\one\fff,\pff n\two\qff\in\qff \zzz$\pss
has a unique solution.\oss
It follows that\dss $\delta$\dss is the restriction of a diagonal map in this case also.

This completes the proof of the theorem when\qss $\delta$\qss is equal to\dss $\oold$\dss
on all summands of\qss $\mathbold{K}\zero(c)$\qss except one.\oss
In order to deduce from this case the general one,\oss
one needs only to repeat word by word the last paragraph of the proof
of Theorem\qss \ref{realization-of-delta}.\oss  \eproof

\mypar{Theorem.}{3-components-identity} \emph{Suppose that\oss $\partial P$\qss
consists of\oss $\leqslant\qff 3$\qss circles for every component\qss $P$\qss of\pss
$\ccomp Q$\nnsp.\oss 
If\qss 
$\varphi$\pss
can be extended to a Torelli diffeomorphism of\qss $S$\nnsp,\oss
then\qss $\tau\circ (\varphi\dff\backslash S\halfff)$\pss 
is a Torelli diffeomorphism for some multi-twist diffeomorphism\dss $\tau$\dss
of\pss $S$\qss about\dss $c$\nnsp.}

\proof\qss 
By Theorem\qss \ref{char-ext}\qss $\varphi$\qss is a weakly Torelli diffeomorphism
and\qss $\delta_{\fff \varphi}$\qss is completely reducible.\oss
By Theorem\qss \ref{symmetry-theorem}\qss $\delta_{\fff \varphi}$\qss is symmetric.\oss
Theorem\qss \ref{3-diagonal}\qss implies that\qss $\delta_{\fff \varphi}$\qss
is equal to the restriction of a diagonal map.\oss
It remains to apply Theorem\qss \ref{twist-correctable}.\oss  \eproof

\mysection{Torelli\qss groups\qss of\qss surfaces\qss with\qss boundary\fff?}{torelli-groups}

\vspace*{\bigskipamount}
What is the Torelli group of a connected subsurface\dss $Q$\dss
of a closed orientable surface\dss $S$\nsp?\oss
The first real insight into this question is due to A. Putman\oss \cite{p}.\oss
He suggested that a good candidate is the group\dss $\ppp(Q)$\qss of isotopy classes
of diffeomorphisms of\dss $Q$\dss
fixed on the boundary\dss $\partial Q$\dss and
such that their extensions by the identity to\dss $S$\dss 
act trivially on\dss $H\one(S)$\dnsp.\oss
The isotopies are also supposed to be fixed on\dss $\partial Q$\nnsp.\oss
Putman showed that\dss $\ppp(Q)$\dss depends only on the partition of\qss
$\partial Q$\qss into the boundaries\qss $\partial P$\qss of components\dss $P$\dss of
the complementary subsurface\qss $\ccomp Q$\dnsp.\oss

The present paper suggests that another good candidate is the group\dss $\ttt(Q)$\dss
of isotopy classes
of diffeomorphisms of\dss $Q$\dss
fixed on the boundary\dss $\partial Q$\dss and
and admitting an extension to a diffeomorphism of\dss $S$\dss
acting trivially on\dss $H\one(S)$\dnsp.\oss
The isotopies are again supposed to be fixed on\dss $\partial Q$\nnsp.\qff\oss
Like\dss $\ppp(Q)$\dnsp,\oss the group\qss $\ttt(Q)$\qss depends only on the partition of\qss
$\partial Q$\qss into the boundaries\qss $\partial P$\qss of components\dss $P$\dss of
the complementary subsurface\qss $\ccomp Q$\dnsp.\footnotemark

\footnotetext{The use of the same notation\dss $\ttt(\bullet)$\dss as the one used by D. Johnson
for the Torelli groups of closed surfaces should not be construed as
the desire to dismiss other candidates.}

By the very definition,\oss the groups\qss $\ttt(Q)$\qss have better
restriction properties\qss ({\fff\halfff}being closed under restrictions is one
of two properties put forward by Putman as the motivation of his definition).\oss
Suppose that\dss $\psi$\dss is a diffeomorphism of\dss $S$\dss
leaving\dss $Q$\dss invariant and fixed on\qss $\partial Q$\dnsp,\oss
and let\dss $\varphi$\dss be the induced diffeomorphism\qss $Q\toto Q$\dnsp.\oss
If the isotopy class of\dss $\psi$\dss belongs to\dss $\ttt(S)$\dnsp,\oss
then the isotopy class of\dss $\varphi$\dss belongs to\dss $\ttt(Q)$\dnsp,\oss
but not necessarily belongs to\dss $\ppp(Q)$\dnsp.\oss

The groups\qss $\ppp(Q)$\qss and\qss $\ttt(Q)$\qss 
are related by the short exact sequence\vspace*{-3pt}
\[
\hspace*{0.5em}
\begin{tikzcd}[column sep=scriptsize]
\mathbold{\oldstylenums{1}} 
\arrow[r]   
& \ppp(Q) \arrow[r] 
& \ttt(Q) \arrow[r, "{\displaystyle \delta}"]
& \mathbold{D}\dff(c) \arrow[r] 
& \mathbold{\oldstylenums{1}}  
\end{tikzcd}
\]

\vspace*{-12pt}
where\qss $\mathbold{D}\dff(c)$\qss is the group of completely reducible
symmetric map\oss
$\displaystyle
\mathbold{K}\zero(c)\ttoo \mathbold{H}\one(c)$\dnsp,\oss
and the\qss \emph{difference homomorphism}\qss $\delta$\dss 
maps the isotopy class of a diffeomorphism\dss $\varphi$\dss
into its difference map\qss $\delta_{\fff \varphi}$\nsp.\oss
It seems that\dss $\ttt(Q)$\dss
should be considered together with the homomorphism\dss $\delta$\dss
and the pair\qss $(\ttt(Q)\fff,\pff \delta)$\qss is an even better
candidate for the title of the Torelli group of\dss $Q$\nnsp.\oss

Suppose that\dss $S$\dss is partitioned into the union
of several subsurfaces such that their interiors are disjoint and every
two of them intersect by several common boundary components.\oss
The elements of groups\dss $\ppp(Q)$\dss for subsurfaces\dss $Q$\dss of
this partition can be pasted together into an element of the Torelli group\dss $\tors$\dnsp,\oss
but many other elements can be also pasted.\oss
If the elements of groups\dss $\mmod(Q)$\dss
can be pasted,\oss then they belong to the groups\qss $\ttt(Q)$\dnsp,\oss
but not all collections of elements of\qss $\ttt(Q)$\qss can be pasted,\oss
as shown in the sequel\qss \cite{i-et-iii}\qss of this paper\halfff.\oss
Still,\oss the groups\qss $\ttt(Q)$\qss together with the\dss $\delta$\dnsp-difference maps
provide a good framework for analyzing the pasting problem
and the classification of abelian
subgroups of\qss $\tors$\dnsp.\oss

\myappend{A\qss converse\qss to\qss Theorem\qss \ref{weakly-torelli}}{weakly-torelli-converse}

\vspace*{6pt}
\myappar{Theorem.}{weakly-torelli-converse}
\emph{If for every\qss $\alpha\qff \in\qff \mathbold{Z}\one(Q\fff,\pff c)$\qss
the cycle\qss
$\varphi_{\dff *}(\alpha)\qff -\qff \alpha$\qss 
is homologous in\dss $Q$\dss to a cycle in\dss $c$\nnsp,\oss
then\qss $\varphi$\qss is a weakly Torelli diffeomorphism.\oss}

\proof\qss
Suppose that\dss $\tau$\dss is a cycle in\dss $Q$\nnsp.\oss
Let\qss $\alpha\qff \in\qff \mathbold{Z}\one(Q\fff,\pff c)$\dnsp.\qff\oss
Let\vspace*{3pt}
\[
\quad
\sigma
\off =\off  
\varphi_{\fff *\fff}(\tau)\qff -\qff \tau
\hspace*{1.5em}\mbox{ and }\hspace*{1.5em}
\delta
\off =\off 
\varphi_{\fff *\fff}(\alpha)\qff -\qff \alpha\dff.
\]

\vspace*{-9pt}
Suppose that\qss
$\beta\qff \in\qff Z\one(\ccomp Q\fff,\pff c)$\qss
and\oss $\partial \alpha\off =\off \partial \beta$\nnsp.\oss
Then\qss\vspace*{3pt}
\[
\quad
\gamma
\off =\off 
\alpha\qff -\qff \beta
\]

\vspace*{-9pt}
is a cycle.\oss
The first step of the proof of Theorem\qss \ref{weakly-torelli}\qss
does not use the assumption that\dss $\varphi$\dss is a weakly Torelli diffeomorphism
and implies that\vspace*{3pt}
\begin{equation*}
\label{sum-is-zero-appendix}
\quad
\langle\qff \delta\fff,\pff \tau \qff\rangle
\qff +\qff
\langle\qff \gamma\fff,\pff \sigma \qff\rangle
\qff +\qff
\langle\qff \delta\fff,\pff \sigma \qff\rangle
\off =\off
\oold\dff.
\end{equation*}

\vspace{-9pt}
By the assumption\dss $\delta$\dss is homologous to a cycle in\dss $c$\nnsp.\oss
Since\qss $\tau\fff,\pff \sigma$\qss are cycles in\dss $Q$\nnsp,\oss
it follows that\oss
$\langle\qff \delta\fff,\pff \tau \qff\rangle
\off =\off
\langle\qff \delta\fff,\pff \sigma \qff\rangle
\off =\off
\oold$\nnsp.\oss
In view of the last displayed equality,\oss
this implies that\oss
$\langle\qff \gamma\fff,\pff \sigma \qff\rangle
\off =\off
\oold$\nnsp.\oss
Let\qss
$a$\dss be an arbitrary element of\dss $H\one(S)$\dnsp.\oss
By the Mayer--Vietoris equivalence the homology class\dss 
$a$\dss can be represented by
a cycle of the form\qss $\gamma\qff =\qff \alpha\qff -\qff \beta$\qss
and hence\oss\vspace*{3pt}
\[
\quad
\langle\qff a\fff,\pff \hclass{\sigma}_{\dff S} \qff\rangle
\off =\off 
\langle\qff \gamma\fff,\pff \sigma \qff \rangle
\off =\off
\oold\
\]

\vspace{-9pt}
Since\qss $a\qff \in\qff H\one(S)$\qss
is arbitrary,\oss
it follows that\qss $\hclass{\sigma}_{\dff S}\qff =\qff \oold$\nnsp.\oss
In other words,\oss the cycle\qss
$\varphi_{\fff *\fff}(\tau)\qff -\qff \tau$\qss
is a boundary in\dss $S$\dss for every\qss $\tau\qff \in\qff Z\one(Q)$\qss
and hence\dss $\varphi$\dss is weakly Torelli.\oss  \eproof

\myitpar{Remark.}
Let us assume only that\qss
$\varphi_{\dff *}(\alpha)\qff -\qff \alpha$\qss 
is homologous in\dss $Q$\dss to a cycle in\dss $c$\dss
for every\qss \emph{cycle}\qss $\alpha$\qss in\dss $S$\nnsp.\oss
As in Section\qss \ref{weakly-torelli-section},\oss
one can define the difference class\qss $\Delta_{\dff \varphi}\dff(\alpha)$\qss
of a cycle\dss $\alpha$\dss in\dss $Q$\dss
as the image in\dss $\mathbold{H}\one(c)$\dss of the homology class\qss
${\hclass{\gamma}\qff~\in\qff~H\one(c)}$\qss of any cycle\dss
$\gamma\qff \in\qff Z\one(c)$\qss 
homologous to\qss 
$\varphi_{\dff *}(\alpha)\qff -\qff \alpha$\qss in\dss $Q$\nnsp.\oss
Clearly,\oss the difference class\qss $\Delta_{\dff \varphi}\dff(\alpha)$\qss
depends only on the homology class\qss
$\hclass{\alpha}\qff \in\qff H\one(Q)$\dnsp,\oss
and hence these difference classes define a homomorphism\vspace*{3pt}
\[
\quad
d_{\dff \varphi}
\qff \colon\qff
H\one(Q)\ttoo \mathbold{H}\one(c)
\]

\vspace*{-9pt}
It turns out that every homomorphism\qss
$H\one(Q)\toto \mathbold{H}\one(c)$\qss
equal to zero on the image of\dss
$H\one(c)$\dss in\dss $H\one(Q)$\dss
can be realized as\qss
$d_{\dff \varphi}$\qss
for some diffeomorphism\dss $\varphi$\dss of\dss $Q$\dss fixed on\dss $c$\dss
and satisfying the above assumption.\oss
This follows from a construction of\qss D.\dss Johnson.\oss
See\qss \cite{j},\oss Lemma\qss 5\qss and\dss Appendix\qss I.\oss
In particular\halfff,\oss
such a diffeomorphism\dss $\varphi$\dss
does not need to be weakly Torelli.

\begin{flushright}

July\qss \old{8},\oss \old{2016},\oss
and\oss
April\qss \old{4},\oss \old{2017}\qss (the present version)
 
https\halfff:/\!/\hspace*{-0.06em}nikolaivivanov.com

E-mail\halfff:\oss nikolai.v.ivanov{\fff}@{\dff}icloud.com

\end{flushright}

\end{document}